\documentclass[11pt,leqno]{article}
\usepackage{amsfonts}
\usepackage{amsmath}
\usepackage{amssymb}
\usepackage{theorem}
\usepackage{latexsym}
\usepackage{mathrsfs}

\newtheorem{theorem}{\bf Theorem}
\newtheorem{lemma}[theorem]{\bf Lemma}

\newtheorem{remark}[theorem]{\bf Remark}

\def\thebibliography#1{\section*{References\markboth
 {REFERENCES}{REFERENCES}}\list
 {[\arabic{enumi}]}{\settowidth\labelwidth{[#1]}\leftmargin\labelwidth
 \advance\leftmargin\labelsep
 \usecounter{enumi}}
 \def\newblock{\hskip .11em plus .33em minus .07em}
 \sloppy
 \sfcode`\.=1000\relax}

\begin{document}
\vspace*{0ex}
\begin{center}
{\Large\bf
A mathematical justification of the Isobe-Kakinuma model \\[0.5ex]
for water waves with and without bottom topography 
}
\end{center}

\begin{center}
Tatsuo Iguchi
\end{center}

\begin{abstract}
We consider the Isobe--Kakinuma model for water waves in both cases of the flat and the variable bottoms. 
The Isobe--Kakinuma model is a system of Euler--Lagrange equations for an approximate Lagrangian 
which is derived from Luke's Lagrangian for water waves by approximating the velocity potential 
in the Lagrangian appropriately. 
The Isobe--Kakinuma model consists of $(N+1)$ second order and a first order partial differential equations, 
where $N$ is a nonnegative integer. 
We justify rigorously the Isobe--Kakinuma model as a higher order shallow water approximation in the strongly 
nonlinear regime by giving an error estimate between the solutions of the Isobe--Kakinuma model and of the 
full water wave problem in terms of the small nondimensional parameter $\delta$, 
which is the ratio of the mean depth to the typical wavelength. 
It turns out that the error is of order $O(\delta^{4N+2})$ in the case of the flat bottom and of order 
$O(\delta^{4[N/2]+2})$ in the case of variable bottoms. 
\end{abstract}

\section{Introduction}
\label{section:intro}
In this paper we consider the motion of a water filled in $(n+1)$-dimensional Euclidean space 
together with the water surface. 
The water wave problem is mathematically formulated as a free boundary problem 
for an irrotational flow of an inviscid and incompressible fluid under the gravitational field. 
Let $t$ be the time, $x = (x_1,\ldots,x_n)$ the horizontal spatial coordinates, 
and $z$ the vertical spatial coordinate. 
We assume that the water surface and the bottom are represented as $z = \eta(x,t)$ and $z = - h + b(x)$, 
respectively, where $\eta = \eta(x,t)$ is the surface elevation, $h$ is the mean depth, 
and $b = b(x)$ represents the bottom topography. 
As was shown by J. C. Luke \cite{Luke1967}, the water wave problem has a variational structure. 
His Lagrangian density is of the form 
\begin{equation}\label{intro:Luke's Lagrangian}
\mathscr{L}_{\rm Luke}(\Phi,\eta) = \int_{-h+b(x)}^{\eta(x,t)}\biggl(\partial_t\Phi(x,z,t)
 +\frac12|\nabla_X\Phi(x,z,t)|^2+gz\biggr){\rm d}z,
\end{equation}
where $\Phi = \Phi(x,z,t)$ is the velocity potential, 
$\nabla_X=(\nabla,\partial_z)=(\partial_{x_1},\ldots,\partial_{x_n},\partial_z)$, and 
$g$ is the gravitational constant. 
He showed that the corresponding Euler--Lagrange equation is exactly the basic equations for water waves. 
Concerning the water wave problem, we refer to H. Lamb \cite{Lamb1993}, J. J. Stoker \cite{Stoker1992}, 
and D. Lannes \cite{Lannes2013-2}.

M. Isobe \cite{Isobe1994, Isobe1994-2} and T. Kakinuma \cite{Kakinuma2000, Kakinuma2001, Kakinuma2003} 
approximated the velocity potential $\Phi$ in Luke's Lagrangian by 
\[
\Phi^{\mbox{\rm\tiny app}}(x,z,t) = \sum_{i=0}^N\Psi_i(z;b)\phi_i(x,t),
\]
where $\{\Psi_i\}$ is an appropriate function system in the vertical coordinate $z$ and may depend on 
the bottom topography $b$ and $(\phi_0,\phi_1,\ldots,\phi_N)$ are unknown variables, 
and derived an approximate Lagrangian density 
$\mathscr{L}^{\mbox{\rm\tiny app}}(\phi_0,\phi_1,\ldots,\phi_N,\eta)
=\mathscr{L}_{\rm Luke}(\Phi^{\mbox{\rm\tiny app}},\eta)$. 
The Isobe--Kakinuma model is the corresponding Euler--Lagrange equation for the 
approximated Lagrangian. 
We have to choose the function system $\{\Psi_i\}$ carefully in order that the Isobe-Kakinuma model 
would be a good approximation for the full water wave problem. 
One of the choices is obtained by the bases of a Taylor series of the velocity potential $\Phi(x,z,t)$ 
with respect to the vertical spatial coordinate $z$ around the bottom. 
Such an expansion has been already used by J. Boussinesq \cite{Boussinesq1872} in the case of the flat bottom. 
From this point of view, one of the natural choices of the function system is given by 
\[
\Psi_i(z;b) 
 = 
\begin{cases}
 (z + h)^{2i} & \mbox{in the case of the flat bottom}, \\
 (z + h - b(x))^i & \mbox{in the case of the variable bottom}.
\end{cases}
\]
Here we note that the later choice is valid also for the case of the flat bottom. 
However, it turns out that the terms of odd degree do not play any important role in such a case 
so that the former choice economizes the computational resources in the numerical computations. 
In this paper, to treat the both cases at the same time, we adopt the approximation 
\begin{equation}\label{intro:approximation}
\Phi^{\mbox{\rm\tiny app}}(x,z,t) = \sum_{i=0}^N(z+h-b(x))^{p_i}\phi_i(x,t),
\end{equation}
where $p_0,p_1,\ldots,p_N$ are nonnegative integers satisfying $0=p_0<p_1<\cdots<p_N$. 
Then, the corresponding Isobe--Kakinuma model has the form 
\begin{equation}\label{intro:IK model}
\left\{
 \begin{array}{l}
  \displaystyle
  H^{p_i} \partial_t \eta + \sum_{j=0}^N \biggl\{ \nabla \cdot \biggl(
   \frac{1}{p_i+p_j+1} H^{p_i+p_j+1} \nabla\phi_j
   - \frac{p_j}{p_i+p_j} H^{p_i+p_j} \phi_j \nabla b \biggr) \\
  \displaystyle\phantom{ H^{p_i} \partial_t \eta + \sum_{j=0}^N \biggl\{ }
   + \frac{p_i}{p_i+p_j} H^{p_i+p_j} \nabla b \cdot \nabla\phi_j
   - \frac{p_ip_j}{p_i+p_j-1} H^{p_i+p_j-1} (1 + |\nabla b|^2) \phi_j \biggr\} = 0 \\
  \makebox[27em]{}\mbox{for}\quad i=0,1,\ldots,N, \\
  \displaystyle
  \sum_{j=0}^N H^{p_j} \partial_t \phi_j + g\eta 
   + \frac12 \biggl\{ \biggl| \sum_{j=0}^N ( H^{p_j}\nabla\phi_j - p_j H^{p_j-1}\phi_j\nabla b ) \biggr|^2 
   + \biggl( \sum_{j=0}^N p_j H^{p_j-1} \phi_j \biggr)^2 \biggr\} = 0,
 \end{array}
\right.
\end{equation}
where $H = H(x,t) = h + \eta(x,t) - b(x)$ is the depth of the water. 
Here and in what follows we use the notational convention $0/0=0$. 
This is the Isobe--Kakinuma model that we are going to consider in this paper. 
This system consists of $(N+1)$ evolution equations for $\eta$ and only one evolution equation 
for $(N+1)$ unknowns $(\phi_0,\phi_1,\ldots,\phi_N)$, so that this is a overdetermined and 
underdetermined composite system. 
For more details to this model, we refer R. Nemoto and T. Iguchi \cite{NemotoIguchi2017}.

One of the interesting features of the model is its linear dispersion relation. 
Let $c_{IK}(\xi)$ and $c_{WW}(\xi)$ be the phase speed of the plane wave solution related to the wave 
vector $\xi \in \mathbf{R}^n$ of the linearized Isobe--Kakinuma model and the linearized water wave problem 
around the rest state in the case of the flat bottom, respectively. 
Then, under the choice $p_i=2i$, $(c_{IK}(\xi))^2$ becomes $[2N/2N]$ Pad\'e approximant of $(c_{WW}(\xi))^2$. 
More precisely, it holds that 
\[
\biggl|\biggl(\frac{c_{WW}(\xi)}{\sqrt{gh}}\biggr)^2-\biggl(\frac{c_{IK}(\xi)}{\sqrt{gh}}\biggr)^2\biggr|
\leq C(h|\xi|)^{4N+2}
\]
with a positive constant $C$ depending only on $N$. 
Under the choice $p_i=i$, we do not have such a beautiful result as above, 
but we still have following nice estimate 
\[
\biggl|\biggl(\frac{c_{WW}(\xi)}{\sqrt{gh}}\biggr)^2-\biggl(\frac{c_{IK}(\xi)}{\sqrt{gh}}\biggr)^2\biggr|
\leq C(h|\xi|)^{4[N/2]+2}. 
\]
For the details, we refer to R. Nemoto and T. Iguchi \cite{NemotoIguchi2017}. 
Since $h|\xi|$ is essentially the ratio of the mean depth to the wave length, these estimates anticipate 
that the Isobe--Kakinuma model would be an approximation to the full water wave problem with an error of 
order $O(\delta^{4N+2})$ in the case of the flat bottom with the choice $p_i=2i$, 
and of order $O(\delta^{4[N/2]+2})$ in the case of variable bottoms with the choice $p_i=i$, 
where $\delta$ is a nondimensional parameter defined as the aspect ratio. 
In this paper, we will show that this is correct even for the nonlinear problem with variable bottoms. 
In a particular case, that is, in the case $N=1$ and $p_1 = 2$ with the flat bottom, this was shown by 
T. Iguchi \cite{Iguchi2017}. 
Therefore, this paper is a generalization of his results.

In order to compare this Isobe--Kakinuma model with the full water wave problem more precisely 
in the shallow water regime, we need to rewrite \eqref{intro:IK model} in an appropriate nondimensional form. 
Let $\lambda$ be the typical wave length and introduce a nondimensional parameter 
$\delta$ by the aspect ratio $\delta=h/\lambda$, which measures the shallowness of the water. 
We rescale the independent and the dependent variables by 
\[
  x=\lambda\tilde{x}, \quad z=h\tilde{z}, \quad t=\frac{\lambda}{\sqrt{gh}}\tilde{t}, 
  \quad \eta=h\tilde{\eta}, \quad \phi_i=\frac{\lambda\sqrt{gh}}{\lambda^{p_i}}\tilde{\phi}_i. 
\]
Here we note that these rescaling of dependent variables are related to the 
strongly nonlinear regime of the wave. 
Plugging these into \eqref{intro:IK model} and dropping the tilde sign in the notation we obtain the 
Isobe--Kakinuma model in the nondimensional form. 
It follows directly from the equations that if the solution and its derivatives are uniformly bounded 
with respect to the small parameter $\delta$, then $\phi_i$ is of order $O(\delta)$ in the case where 
$p_i$ is an odd integer. 
Therefore, it is more convenient to rescale $\phi_i$ again by 
\[
\phi_i = \delta^{p_i - 2[p_i/2]}\tilde{\phi}_i,
\]
where $[p_i/2]$ is the integer part of $p_i/2$. 
Note that $p_i - 2[p_i/2]=0$ if $p_i$ is even, $=1$ if $p_i$ is odd. 
Then, the Isobe--Kakinuma model in nondimensional form has the form 
\begin{equation}\label{intro:ndIK model}
\left\{
 \begin{array}{l}
  \displaystyle
  H^{p_i} \partial_t \eta + \sum_{j=0}^N \delta^{2(p_j - [p_j/2])} \biggl\{ \nabla \cdot \biggl(
   \frac{1}{p_i+p_j+1} H^{p_i+p_j+1} \nabla\phi_j
   - \frac{p_j}{p_i+p_j} H^{p_i+p_j} \phi_j \nabla b \biggr) \\
  \displaystyle\phantom{ H^{p_i}\partial_t \eta + \sum_{j=0}^N }
   + \frac{p_i}{p_i+p_j} H^{p_i+p_j} \nabla b \cdot \nabla\phi_j
   - \frac{p_ip_j}{p_i+p_j-1} H^{p_i+p_j-1} (\delta^{-2} + |\nabla b|^2) \phi_j \biggr\} = 0 \\
  \makebox[27em]{}\mbox{for}\quad i=0,1,\ldots,N, \\
  \displaystyle
  \sum_{j=0}^N \delta^{2(p_j - [p_j/2])} H^{p_j} \partial_t \phi_j + \eta 
   + \frac12 \biggl\{ \biggl| \sum_{j=0}^N \delta^{2(p_j - [p_j/2])}
    ( H^{p_j}\nabla\phi_j - p_j H^{p_j-1}\phi_j\nabla b ) \biggr|^2 \\
  \displaystyle\phantom{\sum_{j=0}^N \delta^{2(p_j - [p_j/2])} H^{p_j} \partial_t \phi_j
     + \eta + \frac12 \biggl\{ }
   + \delta^2\biggl( \sum_{j=0}^N \delta^{2(p_j - [p_j/2] - 1)} p_j H^{p_j-1} \phi_j \biggr)^2 \biggr\} = 0,
 \end{array}
\right.
\end{equation}
where $H = H(x,t) = 1 + \eta(x,t) + b(x)$. 
We consider the initial value problem to this Isobe--Kakinuma model \eqref{intro:ndIK model} 
under the initial conditions 
\begin{equation}\label{intro:initial conditions}
(\eta,\phi_0,\ldots,\phi_N)=(\eta_{(0)},\phi_{0(0)},\ldots,\phi_{N(0)}) \quad\makebox[3em]{at} t=0.
\end{equation}
Solvability of the initial value problem \eqref{intro:ndIK model}--\eqref{intro:initial conditions} 
was first given by Y. Murakami and T. Iguchi \cite{MurakamiIguchi2015} in a particular case and then by 
R. Nemoto and T. Iguchi \cite{NemotoIguchi2017} in the general case under physically reasonable 
conditions on the initial data.

On the other hand, the initial value problem to the full water wave problem 
in Zakharov--Craig--Sulem formulation in the nondimensional form is written as 
\begin{equation}\label{intro:WW}
\left\{
 \begin{array}{l}
  \partial_t\eta-\Lambda(\eta,b,\delta)\phi=0, \\
  \displaystyle
  \partial_t\phi+\eta+\frac12|\nabla\phi|^2
   -\delta^2\frac{(\Lambda(\eta,b,\delta)\phi+\nabla\eta\cdot\nabla\phi)^2}{2(1+\delta^2|\nabla\eta|^2)}
   =0,
 \end{array}
\right.
\end{equation}
\begin{equation}\label{intro:IC of WW}
(\eta,\phi)=(\eta_{(0)},\phi_{(0)})  \quad\makebox[3em]{at} t=0,
\end{equation}
where $\phi=\phi(x,t)$ is the trace of the velocity potential $\Phi$ on the water surface, that is, 
$\phi(x,t) = \Phi(x,\eta(x,t),t)$ 
and $\Lambda(\eta,b,\delta)$ is the Dirichlet-to-Neumann map for Laplace's equation. 
More precisely, the linear operator $\Lambda(\eta,b,\delta)$ depending nonlinearly on 
the surface elevation $\eta$, the bottom topography $b$, and the parameter $\delta$ is defined by 
\begin{equation}\label{intro:DN}
\Lambda(\eta,b,\delta)\phi = (\delta^{-2}\partial_z\Phi - \nabla\eta\cdot\nabla\Phi)|_{z=\eta(x,t)},
\end{equation}
where $\Phi$ is a unique solution to the boundary value problem for Laplace's equation 
\begin{equation}\label{intro:BBP}
\left\{
 \begin{array}{lll}
  \Delta\Phi + \delta^{-2} \partial_z^2\Phi = 0 & \mbox{in} & -1 + b(x) < z < \eta(x,t), \\
  \Phi = \phi & \mbox{on} & z = \eta(x,t), \\
  \delta^{-2} \partial_z \Phi - \nabla b \cdot \nabla\Phi = 0 & \mbox{on} & z = -1 + b(x).
 \end{array}
\right.
\end{equation}
For the details of this formulation to the water wave problem, we refer to V. E. Zakharov \cite{Zakharov1968}, 
W. Craig and C. Sulem \cite{CraigSulem1993}, and D. Lannes \cite{Lannes2013-2}. 
In this paper we will give an error estimate between the solutions of the initial value problems to the 
Isobe--Kakinuma model \eqref{intro:ndIK model}--\eqref{intro:initial conditions} and to the full water 
wave problem \eqref{intro:WW}--\eqref{intro:IC of WW} under appropriate conditions on the initial data.

There are huge literatures devoted to modelization for the full water wave problem and 
many approximate models were proposed and analyzed, especially, in weekly nonlinear regimes. 
Even in the strongly nonlinear regime, there are several model equations. 
Among them, the most famous model is the shallow water equations, 
which is also called Saint-Venant equations. 
The equations in the full water wave problem \eqref{intro:WW} can be expanded with respect to 
$\delta^2$ and the shallow water equations are derived in the limit $\delta \to +0$, so that 
the shallow water equations are the approximation of the full water wave problem with an error $O(\delta^2)$. 
This approximation of the equations leads to the approximation of the solution in the same order of the error as 
\[
|\eta^{\mbox{\tiny WW}}(x,t)-\eta^{\mbox{\tiny SW}}(x,t)| \lesssim \delta^2
\]
on some time interval independent of $\delta$, where $\eta^{\mbox{\tiny WW}}$ and $\eta^{\mbox{\tiny SW}}$ 
are the solutions to the full water wave problem and to the shallow water equations, respectively. 
The other famous model in the strongly nonlinear regime is the Green--Naghdi equations, 
which are derived by introducing the vertically averaged horizontal velocity field and by retaining 
the terms of order $\delta^2$ in the expansion of the equations. 
Therefore, the Green--Naghdi equations are the approximation of the full water wave problem 
with an error $O(\delta^4)$. 
This approximation of the equations leads again to the approximation of the solution in the same order 
of the error as 
\[
|\eta^{\mbox{\tiny WW}}(x,t)-\eta^{\mbox{\tiny GN}}(x,t)| \lesssim \delta^4
\]
on some time interval independent of $\delta$, where $\eta^{\mbox{\tiny GN}}$ is a solution to 
the Green--Naghdi equations. 
Concerning these and related results, we refer to Y. A. Li \cite{Li2006}, T. Iguchi \cite{Iguchi2009, Iguchi2011}, 
B. Alvarez-Samaniego and D. Lannes \cite{Alvarez-SamaniegoLannes2008}, and 
H. Fujiwara and T. Iguchi \cite{FujiwaraIguchi2015}.

Compared to these approximations, the Isobe--Kakinuma model \eqref{intro:ndIK model} is not the approximation 
of the equations so that it is not straight forward to derive a precise error estimate, even in the formal level. 
To analyze the shallow water approximation, we will further restrict ourselves to the following two cases: 
\begin{itemize}
\item[(H1)]
$p_i = 2i$ $(i=0,1,\ldots,N)$ with the flat bottom, that is, $b(x) \equiv 0$
\item[(H2)]
$p_i = i$ $(i=0,1,\ldots,N)$ with general bottom topographies
\end{itemize}
We rewrite the approximation \eqref{intro:approximation} in the nondimensional form as 
\begin{equation}\label{intro:nd approximation}
\Phi^{\mbox{\rm\tiny app}}(x,z,t) = \sum_{i=0}^N \delta^{2(p_i - [p_i/2])}(z + 1 - b(x))^{p_i} \phi_i(x,t),
\end{equation}
while as was shown by J. Boussinesq \cite{Boussinesq1872}, in the case of the flat bottom 
(which corresponds to the case (H1)) the velocity potential $\Phi$ 
which satisfies the boundary value problem \eqref{intro:BBP} can be expanded in a Taylor series as 
\[
\Phi(x,z,t) = \sum_{i=0}^{\infty} \delta^{2i}(z + 1)^{2i} \frac{(-\Delta)^i\phi_0(x,t)}{(2i)!},
\]
where $\phi_0$ is the trace of the velocity potential $\Phi$ on the bottom. 
Although $\phi_i$ is not equal to $\frac{1}{(2i)!}(-\Delta)^i\phi_0$, 
we may regard \eqref{intro:nd approximation} to an approximation with an error of order $O(\delta^{2N+2})$ 
in the case (H1) and of order $O(\delta^{2[N/2]+2})$ in the case (H2). 
Therefore, one may expect that the Isobe--Kakinuma model would be an approximation with an error of 
these orders. 
However, surprisingly we shall see in this paper that the precise error is much smaller than these 
orders and is given by 
\begin{equation}\label{intro:justification}
|\eta^{\mbox{\tiny WW}}(x,t) - \eta^{\mbox{\tiny IK}}(x,t)| \lesssim 
\begin{cases}
\delta^{4N+2} & \mbox{in the case (H1)}, \\
\delta^{4[N/2]+2} & \mbox{in the case (H2)},
\end{cases}
\end{equation}
as was expected by the analysis of linear dispersion relations, 
where $\eta^{\mbox{\tiny IK}}$ is the solution of the Isobe--Kakinuma model. 
As mentioned before, in the case $N=1$ and $p_1 = 2$ with flat bottom, 
this error estimate was shown by T. Iguchi \cite{Iguchi2017}.

As another higher order shallow water approximation in the strongly nonlinear regime, 
extended Green--Naghdi equations were proposed by Y. Matsuno \cite{Matsuno2015, Matsuno2016}. 
His $\delta^{2N}$ model is an approximation of the full water wave equations with an error of 
order $\delta^{2N+2}$ and contains $(2N+1)$th order derivative terms. 
As is well known, higher order derivative terms are troublesome in a numerical computation. 
Moreover, it is not so easy to write down explicitly the extended Green--Naghdi equations for large $N$. 
We remark also that until now there is no rigorous justification of his $\delta^{2N}$ model 
in the sense of approximation of the solutions as mentioned above. 
Compared to this model, the Isobe--Kakinuma model does not contain any higher order derivative terms. 
This is one of strong advantages of the Isobe--Kakinuma model.

The contents of this paper are as follows. 
In Section \ref{section:result} we present our main results in this paper, that is, an existence of 
the solution of the initial value problem to the Isobe--Kakinuma model 
\eqref{intro:ndIK model}--\eqref{intro:initial conditions} on some time interval independent of the 
parameter $\delta \in (0,1]$, the consistency of the Isobe--Kakinuma model with the water wave 
equations \eqref{intro:WW}, and the rigorous justification of the model by establishing 
an error estimate of the solutions such as \eqref{intro:justification}. 
In Section \ref{section:t-detivative} we derive estimates for the time derivatives of the solution 
to \eqref{intro:ndIK model} 
and related partial differential operators of elliptic type with particular care on the dependence 
on the parameter $\delta$. 
Since the hypersurface $t = 0$ in the space-time $\mathbf{R}^n \times \mathbf{R}$ is characteristic 
for the Isobe--Kakinuma model, these estimations are not straightforward. 
In Section \ref{section:estimate I} we prove the existence of the solution to 
\eqref{intro:ndIK model}--\eqref{intro:initial conditions} on some time interval independent of $\delta$. 
Here, we do not need the special choice of the indices $p_i$. 
In Section \ref{section:estimate II}, under the additional conditions (H1) or (H2) we prove 
uniform boundedness of the solution to \eqref{intro:ndIK model} and its derivatives. 
In Sections \ref{section:consistency I}--\ref{section:consistency III} we prove a consistency of the 
Isobe--Kakinuma model with the water wave equations. 
One of the key elements for obtaining the precise error estimates is to introduce a modified 
approximate velocity potential $\widetilde{\Phi}^{\mbox{\rm\tiny app}}$, which approximates the 
velocity potential $\Phi$ for the full water wave problem with an error of the order 
indicated in \eqref{intro:justification}. 
The approximation $\Phi^{\mbox{\rm\tiny app}}$ does not possess such a nice property. 
In Section \ref{section:justification} we derive an error estimate between the solutions to the 
Isobe--Kakinuma model and to the water wave equations by using the stability of 
the water wave equations.

\medskip
\noindent
{\bf Acknowledgement} \ 
This work was carried out when the author was visiting Universit\'e de Bordeaux 
on his sabbatical leave during the 2017 academic year. 
He is very grateful to the member of Institut de Math\'matiques de Bordeaux, 
especially, David Lannes for their kind hospitalities and for fruitful discussions. 
This work was partially supported by JSPS KAKENHI Grant Number JP17K18742 and JP17H02856.

\section{Main results}
\label{section:result}
\setcounter{equation}{0}
{\bf Notation}. \ 
We denote by $W^{m,p}(\mathbf{R}^n)$ the $L^p$ Sobolev space of order $m$ on $\mathbf{R}^n$. 
The norms of the Sobolev space $H^m=W^{m,2}(\mathbf{R}^n)$ and of a Banach space $B$
are denoted by $\|\cdot\|_m$ and by $\|\cdot\|_B$, respectively. 
The $L^2$-norm and the $L^2$-inner product are simply denoted by $\|\cdot\|$ and $(\cdot,\cdot)_{L^2}$, 
respectively. 
We put $\partial_t=\partial/\partial t$, $\partial_j=\partial/\partial x_j$, 
and $\partial_z=\partial/\partial z$. 
$[P,Q]=PQ-QP$ denotes the commutator. 
We put $J = (1-\Delta)^{1/2}$ and $J_{\delta} = (1-\delta^2\Delta)^{1/2}$. 
$[a]$ denotes the integer part of the real number $a$ and 
$a_1 \vee a_2 = \max\{a_1,a_2\}$. 
We denote the Kronecker delta by $\delta_{ij}$, that is, 
$\delta_{ij}=1$ if $i=j$ and $\delta_{ij}=0$ if $i\ne j$. 
We fix $t_0 > n/2$, which satisfies $t_0 -n/2 \ll 1$ if necessary. 
For a matrix $A$ we denote by $A^{\rm T}$ the transpose of $A$. 
For a vector $\mbox{\boldmath$\phi$}=(\phi_0,\phi_1,\ldots,\phi_N)^{\rm T}$ we denote the last $N$ 
components by $\mbox{\boldmath$\phi$}'=(\phi_1,\ldots,\phi_N)^{\rm T}$. 
We use the notational convention $0/0=0$. 
We denote by $C(a_1,a_2,\ldots)$ a positive constant depending on $a_1,a_2,\ldots$. 
$f \lesssim g$ means that there exists a non-essential positive constant $C$ such that 
$f \leq Cg$ holds.

\bigskip
To state our main results it is convenient to introduce rescaled variables 
$\mbox{\boldmath$\phi$}^\delta = (\phi_0^\delta,\ldots,\phi_N^\delta)^{\rm T}$ by 
\begin{equation}\label{result:rescale}
\phi_i^\delta = \delta^{2(p_i - [p_i/2])} \phi_i
\end{equation}
for $i=0,1,\ldots,N$. 
Then, the Isobe--Kakinuma model \eqref{intro:ndIK model} is written in these rescaled variables as 
\begin{equation}\label{result:IK model}
\left\{
 \begin{array}{l}
  \displaystyle
  H^{p_i} \partial_t \eta + \sum_{j=0}^N  \biggl\{ \nabla \cdot \biggl(
   \frac{1}{p_i+p_j+1} H^{p_i+p_j+1} \nabla\phi_j^\delta
   - \frac{p_j}{p_i+p_j} H^{p_i+p_j} \phi_j^\delta \nabla b \biggr) \\
  \displaystyle\phantom{ H^{p_i}\partial_t \eta + \sum_{j=0}^N }
   + \frac{p_i}{p_i+p_j} H^{p_i+p_j} \nabla b \cdot \nabla\phi_j^\delta
   - \frac{p_ip_j}{p_i+p_j-1} H^{p_i+p_j-1} (\delta^{-2} + |\nabla b|^2) \phi_j^\delta \biggr\} = 0 \\
  \makebox[27em]{}\mbox{for}\quad i=0,1,\ldots,N, \\
  \displaystyle
  \sum_{j=0}^N H^{p_j} \partial_t \phi_j^\delta + \eta 
   + \frac12 \biggl\{ \biggl| \sum_{j=0}^N 
    ( H^{p_j}\nabla\phi_j^\delta - p_j H^{p_j-1}\phi_j^\delta \nabla b ) \biggr|^2 
   + \delta^{-2} \biggl( \sum_{j=0}^N p_j H^{p_j-1} \phi_j^\delta \biggr)^2 \biggr\} = 0.
 \end{array}
\right.
\end{equation}
We denote the corresponding initial data by 
$\mbox{\boldmath$\phi$}_{(0)}^\delta = (\phi_{0(0)}^\delta,\ldots,\phi_{N(0)}^\delta)^{\rm T}$ with 
\begin{equation}\label{result:rescale id}
\phi_{i(0)}^\delta = \delta^{2(p_i - [p_i/2])} \phi_{i(0)}
\end{equation}
for $i=0,1,\ldots,N$ and put 
$\mbox{\boldmath$\phi$}_{(0)}^{\delta \prime} = (\phi_{1(0)}^\delta,\ldots,\phi_{N(0)}^\delta)^{\rm T}$, 
where $\mbox{\boldmath$\phi$}_{(0)} = (\phi_{0(0)},\ldots,\phi_{N(0)})^{\rm T}$ 
is the initial data for the original variables $\mbox{\boldmath$\phi$} = (\phi_0,\ldots,\phi_N)^{\rm T}$ 
in \eqref{intro:initial conditions}.

As explained in the previous section, the Isobe--Kakinuma model \eqref{result:IK model} is a 
overdetermined and underdetermined composite system and we have $(N + 1)$ evolution equations 
for only one unknown $\eta$, so that the initial value problem 
\eqref{intro:ndIK model}--\eqref{intro:initial conditions} is not solvable in general. 
In fact, if the problem has a solution $(\eta,\phi_0,\ldots,\phi_N)$, then by eliminating the 
time derivative $\partial_t\eta$ from the evolution equations we see that the solution has to 
satisfy the relations 
\begin{align}\label{result:compatibility}
& H^{p_i} \sum_{j=0}^N \nabla \cdot \biggl(
   \frac{1}{p_j+1} H^{p_j+1} \nabla\phi_j^\delta
   - \frac{p_j}{p_j} H^{p_j} \phi_j^\delta \nabla b \biggr) \nonumber \\
&= \sum_{j=0}^N \biggl\{ \nabla \cdot \biggl(
   \frac{1}{p_i+p_j+1} H^{p_i+p_j+1} \nabla\phi_j^\delta
   - \frac{p_j}{p_i+p_j} H^{p_i+p_j} \phi_j^\delta \nabla b \biggr)  \\
&\phantom{ =\sum_{j=0}^N\biggl\{ }
 \displaystyle
   + \frac{p_i}{p_i+p_j} H^{p_i+p_j} \nabla b \cdot \nabla\phi_j^\delta
   - \frac{p_ip_j}{p_i+p_j-1} H^{p_i+p_j-1} (\delta^{-2} + |\nabla b|^2)\phi_j^\delta \biggr\} \nonumber 
\end{align}
for $i=1,\ldots,N$. 
Therefore, as a necessary condition the initial date $(\eta_{(0)},\mbox{\boldmath$\phi$}_{(0)})$ and 
the bottom topography $b$ have to satisfy the relation \eqref{result:compatibility} for the existence of 
the solution.

Another important condition on the well-posedness of the Isobe--Kakinuma model is related to a 
generalized Rayleigh--Taylor sign condition for water wave problem and states the positivity of 
the function $a$ defined by 
\begin{align}\label{result:a}
a &= 1 + \sum_{j=1}^N p_j H^{p_j-1} \partial_t \phi_j^\delta \\
&\quad\;
 + \mbox{\boldmath$u$} \cdot \sum_{j=1}^N \Bigl( p_j H^{p_j-1} \nabla\phi_j^\delta
  - p_j(p_j-1) H^{p_j-2} \phi_j^\delta \nabla b \Bigr) 
 + w \sum_{j=1}^N p_j(p_j-1) H^{p_j-2} \phi_j^\delta. \nonumber
\end{align}
For the details of this function $a$, we refer to R. Nemoto and T. Iguchi \cite{NemotoIguchi2017}.

The following theorem is one of the main results in this paper and asserts the existence of 
the solution to the initial value problem \eqref{intro:ndIK model}--\eqref{intro:initial conditions} 
on a time interval independent of the small parameter $\delta$ with a uniform bounds of the solution 
$(\eta,\mbox{\boldmath$\phi$}^\delta)$ in the rescaled variables.

\begin{theorem}\label{result:theorem 1}
Let $c_0,M_0$ be positive constants and $m$ an integer such that $m>n/2+1$. 
There exist a time $T_1>0$ and a constant $M$ such that for any $\delta \in (0,1]$ if the initial 
data $(\eta_{(0)},\mbox{\boldmath$\phi$}_{(0)}^\delta)$ and $b$ satisfy the relation 
\eqref{result:compatibility} and 
\begin{equation}\label{result:uniform estimate 1}
\begin{cases}
 \|(\eta_{(0)},\nabla \mbox{\boldmath$\phi$}_{(0)}^\delta)\|_m
  + \delta^{-1}\|\mbox{\boldmath$\phi$}_{(0)}^{\delta \prime}\|_m + \|b\|_{W^{m+1,\infty}} \leq M_0, \\
 1 + \eta_{(0)}(x) - b(x) \geq c_0, \quad a(x,0)\geq c_0 
   \qquad\mbox{for}\quad x\in\mathbf{R}^n,
\end{cases}
\end{equation}
then the initial value problem \eqref{intro:ndIK model}--\eqref{intro:initial conditions} has a 
unique solution $(\eta,\phi_0,\ldots,\phi_N)$ on the time interval $[0,T_1]$. 
Moreover, the solution satisfies the uniform bound: 
\begin{equation}\label{result:uniform estimate 2}
\begin{cases}
 \|\eta(t)\|_m + \|\nabla \mbox{\boldmath$\phi$}^\delta (t)\|_m
  + \delta^{-1} \|\mbox{\boldmath$\phi$}^{\delta \prime} (t)\|_m
  + \delta^{-2} \|\mbox{\boldmath$\phi$}^{\delta \prime} (t)\|_{m-1} \\
 \quad
  + \|\partial_t \eta(t)\|_{m-1} + \|\partial_t \mbox{\boldmath$\phi$}^\delta (t)\|_m
  + \delta^{-1} \|\partial_t \mbox{\boldmath$\phi$}^{\delta \prime} (t)\|_{m-1}
  + \delta^{-2} \|\partial_t \mbox{\boldmath$\phi$}^{\delta \prime} (t)\|_{m-2} \\
 \quad
  + \|\partial_t^2 \eta(t)\|_{m-2} + \|\partial_t^2 \mbox{\boldmath$\phi$}^\delta (t)\|_{m-1}
  + \delta^{-1} \|\partial_t^2 \mbox{\boldmath$\phi$}^{\delta \prime} (t)\|_{m-2}
  \leq M, \\
 1 + \eta(x,t) - b(x) \geq c_0/2, \quad a(x,t)\geq c_0/2 
   \qquad\mbox{for}\quad x\in\mathbf{R}^n, 0 \leq t \leq T_1,
\end{cases}
\end{equation}
where $\mbox{\boldmath$\phi$}^\delta = (\phi_0^\delta,\ldots,\phi_N^\delta)^{\rm T}$ are the rescaled 
variables defnied by \eqref{result:rescale} and 
$\mbox{\boldmath$\phi$}^{\delta \prime} = (\phi_1^\delta,\ldots,\phi_N^\delta)^{\rm T}$.

Furthermore, if we assume in addition that {\rm (H1)} or {\rm (H2)}, then we have the uniform bound:
\[
\begin{cases}
 \|\phi_j(t)\|_{m-2j+1} \leq M & \mbox{in the case {\rm (H1)}}, \\
 \|\phi_j(t)\|_{m-2[(j+1)/2]+1} \leq M & \mbox{in the case {\rm (H2)}}
\end{cases}
\]
for $i=1,\ldots,N$ and $0 \leq t\leq T_1$. 
\end{theorem}

We proceed to show that the Isobe--Kakinuma model \eqref{intro:ndIK model} is consistent 
with the water wave equations \eqref{intro:WW} at order $O(\delta^{4N+2})$ in the case (H1) 
and at order $O(\delta^{4[N/2]+2})$ in the case (H2). 
In view of the facts that the unknown $\phi$ in Zakharov--Craig--Sulem formulation is the 
trace of the velocity potential $\Phi$ on the water surface and that the unknowns 
$\mbox{\boldmath$\phi$} = (\phi_0,\ldots,\phi_N)^{\rm T}$ for Isobe--Kakinuma model 
appear in the approximation \eqref{intro:nd approximation} of $\Phi$, 
these variables are related approximately by the formula 
\begin{equation}\label{result:phi}
\phi = \sum_{i=0}^N H^{p_i} \phi_i^\delta
 = \sum_{i=0}^N \delta^{2(p_i-[p_i/2])} H^{p_i} \phi_i.
\end{equation}

\begin{theorem}\label{result:theorem 2}
In addition to hypothesis of Theorem \ref{result:theorem 1} we assume that {\rm (H1)} or {\rm (H2)} 
and that $m \geq 4N+2$ and $m > n/2 + 2N +2$ in the case {\rm (H1)}, 
$m \geq 4[N/2]+2+\delta_{N1}$ and $m > n/2+2[N/2]+2$ in the case {\rm (H2)}. 
Let $(\eta,\phi_0,\ldots,\phi_N)$ be the solution obtained in Theorem \ref{result:theorem 1} and 
define $\phi$ by \eqref{result:phi}. 
Then, $(\eta,\phi)$ satisfy the water wave equations approximately as 
\begin{equation}\label{result:appWW}
\left\{
 \begin{array}{l}
  \partial_t \eta - \Lambda(\eta,b,\delta) \phi = \mathfrak{r}_1, \\
  \displaystyle
  \partial_t \phi + \eta + \frac12 |\nabla\phi|^2
   - \delta^2 \frac{(\Lambda(\eta,b,\delta)\phi+\nabla\eta\cdot\nabla\phi)^2}{2(1+\delta^2|\nabla\eta|^2)}
   = \mathfrak{r}_2.
 \end{array}
\right.
\end{equation}
Here, $(\mathfrak{r}_1,\mathfrak{r}_2)$ satisfy 
\begin{equation}
 \begin{cases}
  \|(\mathfrak{r}_1(t),\mathfrak{r}_2(t))\|_{m-4(N+1)} \leq C\delta^{4N+2}
   & \mbox{in the case {\rm (H1)}}, \\
  \|(\mathfrak{r}_1(t),\mathfrak{r}_2(t))\|_{m-4([N/2]+1)} \leq C\delta^{4[N/2]+2}
   & \mbox{in the case {\rm (H2)}},
 \end{cases}
\end{equation}
where $C$ is a positive constant independent of $\delta \in (0,1]$ and $t \in [0,T_1]$. 
\end{theorem}

The above theorem concerns the approximation of the equations. 
Next, we will be concerned with the approximation of the solution to give a 
rigorous justification of the Isobe--Kakinuma model. 
Here we recall the existence theorem for the initial value problem to the water wave equations 
\eqref{intro:WW}--\eqref{intro:IC of WW} obtained by T. Iguchi \cite{Iguchi2009} and 
B. Alvarez-Samaniego and D. Lannes \cite{Alvarez-SamaniegoLannes2008}. 
See also D. Lannes \cite{Lannes2013-2}.

\begin{theorem}\label{result:theorem 3}
Let $c_0,M_0>0$ and $m>n/2+1$. 
There exist a time $T_2>0$ and constants $C,\delta_*>0$ such that for any 
$\delta\in(0,\delta_2]$ if the initial data $(\eta_{(0)},\phi_{(0)})$ satisfy 
\[
\left\{
 \begin{array}{l}
  \|\eta_{(0)}\|_{m+3+1/2}+\|\nabla\phi_{(0)}\|_{m+3} \leq M_0, \\[0.5ex]
  1+\eta_{(0)}(x) \geq c_0 \qquad\mbox{for}\quad x\in\mathbf{R}^n,
 \end{array}
\right.
\]
then the initial value problem \eqref{intro:WW}--\eqref{intro:IC of WW} 
has a unique solution $(\eta,\phi)$ on the time interval $[0,T_2]$. 
Moreover, the solution satisfies the uniform bound: 
\[
\left\{
 \begin{array}{l}
  \|\eta(t)\|_{m+3}+\|\nabla\phi(t)\|_{m+2}
   +\|\partial_t\eta(t)\|_{m+2}+\|\partial_t\phi(t)\|_{m+2} \leq C, \\[0.5ex]
  1+\eta(x,t) \geq c_0/2, 
   \qquad\mbox{for}\quad x\in\mathbf{R}^n, \; 0\leq t\leq T_2. 
 \end{array}
\right.
\]
\end{theorem}

\begin{remark}
{\rm 
In the above theorem the constant $\delta_2$ is small. 
As in the case of Theorem \ref{result:theorem 1} we can reduce the restriction 
$0<\delta\leq\delta_1$ to $0\leq\delta\leq1$, if we impose the sign condition 
$a^{\mbox{\rm\tiny WW}}(x,0)\geq c_0$ on the initial data, where 
$a^{\mbox{\rm\tiny WW}}=1+\delta^2\partial_tZ+\delta^2\mbox{\boldmath $v$}\cdot\nabla Z$ with 
\[
\left\{
 \begin{array}{l}
  Z=(1+\delta^2|\nabla\eta|^2)^{-1}(\Lambda(\eta,\delta)\phi+\nabla\eta\cdot\nabla\phi), \\[0.5ex]
  \mbox{\boldmath $v$}=\nabla\phi-\delta^2Z\nabla\eta.
 \end{array}
\right.
\]
}
\end{remark}

In order that the solution to the Isobe--Kakinuma model 
\eqref{intro:IK model}--\eqref{intro:initial conditions} approximates the solution to the 
water wave problem \eqref{intro:WW}--\eqref{intro:IC of WW}, 
we need to prepare the initial data $\mbox{\boldmath$\phi$}_{(0)}^\delta$ for the Isobe--Kakinuma model 
appropriately in terms of the initial data $\eta_{(0)}$ and $\phi_{(0)}$ for the water wave equations. 
As a matter of fact, the necessary conditions \eqref{result:compatibility} and the approximate relation 
\eqref{result:phi} between $\phi$ and $\mbox{\boldmath$\phi$}^\delta$ determine uniquely the initial data 
$\mbox{\boldmath$\phi$}_{(0)}^\delta$ from $(\eta_{(0)},\phi_{(0)})$, and $b$ as guaranteed by 
Lemma \ref{t-derivative:compatibility} in the Section \ref{section:t-detivative}. 
Moreover, it follows from Lemma \ref{estimate I:lemma 4} in Section \ref{section:estimate I} that 
$\|a(\cdot,0) - 1\|_{m-1} \leq C\delta$ with a constant independent of $\delta$. 
Therefore, by taking $\delta_*$ sufficiently small if necessary, we have $a(x,0) \geq 1/2$, 
so that the conditions \eqref{result:uniform estimate 1} in Theorem \ref{result:theorem 1} will be 
satisfied and we can construct the solution to the Isobe--Kakinuma model. 
The next theorem gives a rigorous justification of the Isobe--Kakinuma model for the water wave problem 
as a higher order shallow water approximation.

\begin{theorem}\label{result:theorem 4}
Let $c_0,M_0$ be positive constants and $m$ an integer such that $m>n/2+1$, 
suppose that {\rm (H1)} or {\rm (H2)} holds, and put $T_*=\min\{T_1,T_2\}$, 
where $T_1$ and $T_2$ are those in Theorems {\rm \ref{result:theorem 1}} and {\rm \ref{result:theorem 3}}. 
Suppose also that $0<\delta\leq\delta_*$ and the initial data $(\eta_{(0)},\phi_{(0)})$ and $b$ satisfy 
\begin{equation}\label{result:condition for ID 2}
 \begin{cases}
  \|\eta_{(0)}\|_{m+4N+8} + \|\nabla\phi_{(0)}\|_{m+4N+7} \leq M_0
   & \mbox{in the case {\rm (H1)}}, \\
  \|\eta_{(0)}\|_{m+4[N/2]+8} + \|\nabla\phi_{(0)}\|_{m+4[N/2]+7} + \|b\|_{W^{m+4[N/2]+8,\infty}} \leq M_0
   & \mbox{in the case {\rm (H2)}}, \\
  1 + \eta_{(0)}(x) - b(x) \geq c_0 \qquad\mbox{for}\quad x\in\mathbf{R}^n. 
 \end{cases}
\end{equation}
Then, \eqref{result:compatibility} and \eqref{result:phi} determine uniquely the initial data 
$\mbox{\boldmath$\phi$}_{(0)}^\delta$ to the Isobe--Kakinuma model. 
Let $(\eta^{\mbox{\rm\tiny WW}},\phi^{\mbox{\rm\tiny WW}})$ be the solution to the initial value 
problem \eqref{intro:WW}--\eqref{intro:IC of WW} obtained in Theorem \ref{result:theorem 3} 
and $(\eta^{\mbox{\rm\tiny IK}},\mbox{\boldmath$\phi$}^\delta)$ the solution to the initial value problem 
\eqref{result:IK model}--\eqref{result:rescale id} obtained in Theorem \ref{result:theorem 1}, 
and define $\phi^{\mbox{\rm\tiny IK}}$ by \eqref{result:phi}. 
Then, for any $\delta\in(0,\delta_*]$ and $t\in [0,T_*]$ we have 
\begin{equation}\label{result:error estimate 2}
\|\eta^{\mbox{\rm\tiny WW}}(t)-\eta^{\mbox{\rm\tiny IK}}(t)\|_{m+2}
 +\|\phi^{\mbox{\rm\tiny WW}}(t)-\phi^{\mbox{\rm\tiny IK}}(t)\|_{m+2} \leq 
\begin{cases}
 C\delta^{4N+2} & \mbox{in the case {\rm (H1)}}, \\
 C\delta^{4[N/2]+2} & \mbox{in the case {\rm (H2)}},
\end{cases}
\end{equation}
where $C$ is a positive constant independent of $\delta$ and $t$. 
\end{theorem}

\begin{remark}
{\rm 
The error estimate \eqref{result:error estimate 2} together with the Sobolev imbedding theorem 
implies the pointwise error estimate \eqref{intro:justification}. 
}
\end{remark}

We will give the proof of Theorems \ref{result:theorem 1}, \ref{result:theorem 2}, and 
\ref{result:theorem 4} in Sections \ref{section:estimate I}--\ref{section:estimate II}, 
\ref{section:consistency I}--\ref{section:consistency III}, and 
\ref{section:justification}, respectively.

\section{Estimate of the time derivate and related operators}
\label{section:t-detivative}
\setcounter{equation}{0}
\setcounter{theorem}{0}
One of the difficulties for the analysis of the Isobe--Kakinuma model \eqref{intro:ndIK model} 
(equivalently \eqref{result:IK model}) lies in the fact that the hypersurface $t=0$ in 
the space-time $\mathbf{R}^n\times\mathbf{R}$ is characteristic for the model. 
In fact, the evolution equation for 
$\mbox{\boldmath$\phi$}^\delta = (\phi_0^\delta,\ldots,\phi_N^\delta)^{\rm T}$ 
is underdetermined so that we cannot express the time derivative 
$\partial_t \mbox{\boldmath$\phi$}^\delta$ 
in terms of the spatial derivatives directly from the equation. 
Nevertheless, we can express it implicitly along with the calculations in R. Nemoto and T. Iguchi 
\cite{NemotoIguchi2017}. 
Since we need to trace carefully the dependence of the small parameter $\delta$, we outline them.

We introduce second order differential operators $L_{ij}=L_{ij}(H,b,\delta)$ $(i,j=0,1,\ldots)$ depending on 
the water depth $H$, the bottom topography $b$, and the parameter $\delta$ by 
\begin{align}\label{t-derivative:L}
L_{ij}\psi_j
&= - \nabla \cdot \biggl(
   \frac{1}{p_i+p_j+1} H^{p_i+p_j+1} \nabla\psi_j
   - \frac{p_j}{p_i+p_j} H^{p_i+p_j} \psi_j \nabla b \biggr) \\[0.5ex]
&\quad\,
  - \frac{p_i}{p_i+p_j} H^{p_i+p_j} \nabla b \cdot \nabla\psi_j
   + \frac{p_ip_j}{p_i+p_j-1} H^{p_i+p_j-1} ( \delta^{-2} + |\nabla b|^2 ) \psi_j. \nonumber
\end{align}
Then, we have $L_{ij}^*=L_{ji}$, where $L_{ij}^*$ is the adjoint operator of $L_{ij}$ in $L^2(\mathbf{R}^n)$. 
We introduce also the functions $\mbox{\boldmath$u$}$ and $w$ by 
\begin{equation}\label{t-derivative:uw}
\mbox{\boldmath$u$} = \sum_{i=0}^N ( H^{p_i} \nabla\phi_i^\delta - p_i H^{p_i-1} \phi_i^\delta \nabla b ), \quad
w = \delta^{-2} \sum_{i=1}^N p_i H^{p_i-1} \phi_i^\delta.
\end{equation}
Since $\mbox{\boldmath$u$} = (\nabla\Phi^{\mbox{\rm\tiny app}})|_{z=\eta}$ and 
$\delta^2 w = (\partial_z \Phi^{\mbox{\rm\tiny app}})|_{z=\eta}$, 
where $\Phi^{\mbox{\rm\tiny app}}$ is the approximate velocity potential defined by 
\eqref{intro:nd approximation}, $\mbox{\boldmath$u$}$ and $\delta^2 w$ represent approximately 
the horizontal and the vertical components of the velocity field on the water surface, respectively. 
We note that both $\mbox{\boldmath$u$}$ and $w$ would be expected of order $O(1)$. 
Then, the Isobe--Kakinuma model \eqref{intro:ndIK model} and the necessary conditions 
\eqref{result:compatibility} can be written simply as 
\begin{equation}\label{t-derivative:IK model}
\left\{
 \begin{array}{l}
  \displaystyle
  H^{p_i} \partial_t \eta - \sum_{j=0}^N L_{ij} \phi_j^\delta = 0 \qquad\mbox{for}\quad i=0,1,\ldots,N, \\
  \displaystyle
  \sum_{j=0}^N H^{p_j} \partial_t \phi_j^\delta + \eta + \frac12( |\mbox{\boldmath$u$}|^2 + \delta^2 w^2 ) = 0
 \end{array}
\right.
\end{equation}
and 
\begin{equation}\label{t-derivative:compatibility}
\sum_{j=0}^N ( L_{ij} - H^{p_i} L_{0j} ) \phi_j^\delta = 0 \qquad\mbox{for}\quad i=1,\ldots,N,
\end{equation}
respectively. 
In view of these equations we introduce also linear operators $\mathscr{L}_i = \mathscr{L}_i(H,b,\delta)$ 
$(i = 0,1,\ldots,N)$ depending on the water depth $H$, the bottom topography $b$, and the parameter 
$\delta$, and acting on 
$\mbox{\boldmath$\varphi$} = (\varphi_0,\ldots,\varphi_N)^{\rm T}$ by 
\begin{equation}\label{t-derivative:scrLi}
\mathscr{L}_0 \mbox{\boldmath$\varphi$} = \sum_{j=0}^N H^{p_j} \varphi_j, \qquad
\mathscr{L}_i \mbox{\boldmath$\varphi$} = \sum_{j=0}^N ( L_{ij} - H^{p_i} L_{0j} ) \varphi_j
 \quad\mbox{for}\quad i=1,\ldots,N, 
\end{equation}
and put 
\begin{equation}\label{t-derivative:scrL}
\mathscr{L} \mbox{\boldmath$\varphi$} = (\mathscr{L}_0 \mbox{\boldmath$\varphi$},
\ldots,\mathscr{L}_N \mbox{\boldmath$\varphi$})^{\rm T}. 
\end{equation}
Then, the necessary conditions \eqref{t-derivative:compatibility} have the simple form 
\begin{equation}\label{t-derivative:compatibility 2}
\mathscr{L}_i \mbox{\boldmath$\phi$}^\delta = 0
 \quad\mbox{for}\quad i=1,\ldots,N.
\end{equation}
We note that the operators $\mathscr{L}_i$ for $i=1,\ldots,N$ can be written explicitly as 
\begin{align}\label{t-derivative:Li}
\mathscr{L}_i \mbox{\boldmath$\varphi$}
&= \sum_{j=0}^N \biggl\{
 - \biggl( \frac{1}{p_i+p_j+1} - \frac{1}{p_j+1} \biggr) H^{p_i+p_j+1} \Delta\varphi_j \\
&\qquad\quad
 + \biggl( \frac{p_j}{p_i+p_j} - \frac{p_j}{p_j} \biggr) H^{p_i+p_j} \nabla \cdot (\varphi_j\nabla b)
 \nonumber \\
&\qquad\quad
 - \frac{p_i}{p_i+p_j} H^{p_i+p_j} \nabla b \cdot \nabla\varphi_j
   + \frac{p_ip_j}{p_i+p_j-1} H^{p_i+p_j-1} ( \delta^{-2} + |\nabla b|^2 ) \varphi_j \biggr\}, 
   \nonumber
\end{align}
and that they do not include the term $\nabla H$. 
Therefore, differentiating \eqref{t-derivative:compatibility 2} with respect to the time $t$ and using the first 
equation in \eqref{t-derivative:IK model} with $i=0$ to eliminate $\partial_t \eta$, we obtain 
\begin{equation}\label{t-derivative:eq phi pre}
\mathscr{L}_i \partial_t \mbox{\boldmath$\phi$}_j^\delta = F_i
 \qquad\mbox{for}\quad i=1,\ldots,N,
\end{equation}
where 
\begin{equation}\label{t-derivative:Fi}
F_i = - \biggl( \biggl( \frac{\partial}{\partial H} \mathscr{L}_i \biggr) \mbox{\boldmath$\phi$}^\delta \biggr)
 \sum_{j=0}^N L_{0j} \phi_j^\delta 
\end{equation}
for $i=1,\ldots,N$. 
We note that $F_i$ does not contain any time derivatives. 
Then, it follows from \eqref{t-derivative:eq phi pre} and the second equation in 
\eqref{t-derivative:IK model} that 
\begin{equation}\label{t-derivative:eq phi}
\mathscr{L} \partial_t \mbox{\boldmath$\phi$}^\delta = \mbox{\boldmath$F$},
\end{equation}
where $\mbox{\boldmath$F$} = (F_0,\ldots,F_N)^{\rm T}$ and 
\begin{equation}\label{t-derivative:F0}
F_0 = - \eta - \frac12( |\mbox{\boldmath$u$}|^2 + \delta^2 w^2 ).
\end{equation}
Therefore, the time derivative of $\mbox{\boldmath$\phi$}^\delta$ can be represented implicitly 
in terms of the spatial derivatives by 
$\partial_t \mbox{\boldmath$\phi$}^\delta = \mathscr{L}^{-1} \mbox{\boldmath$F$}$.

To investigate the operator $\mathscr{L}^{-1}$, assuming $\mbox{\boldmath$F$}$ to be a given function 
we consider the equation 
\begin{equation}\label{t-derivative:elliptic 1}
\mathscr{L} \mbox{\boldmath$\varphi$} = \mbox{\boldmath$F$}.
\end{equation}
Let $\mbox{\boldmath$\varphi$}$ be a solution of this equation. 
It follows from the first component of \eqref{t-derivative:elliptic 1} that 
\begin{equation}\label{t-derivative:reduction}
\varphi_0 = F_0 - \sum_{j=1}^N H^{p_j} \varphi_j.
\end{equation}
Plugging this into the other components of \eqref{t-derivative:elliptic 1} we obtain 
\begin{equation}\label{t-derivative:elliptic 2}
P_i \mbox{\boldmath$\varphi$}' = F_i - (L_{i0} - H^{p_i}L_{00}) F_0 \quad\mbox{for}\quad i=1,\ldots,N,
\end{equation}
where $\mbox{\boldmath$\varphi$}'=(\varphi_1,\ldots,\varphi_N)^{\rm T}$ and 
$P_j = P_j(H,b,\delta)$ $(j=1,\ldots,N)$ are second order differential operators defined by 
\begin{equation}\label{linear:P}
P_i\mbox{\boldmath$\varphi$}'=
\sum_{j=1}^N \{(L_{ij}-H^{p_i}L_{0j})\varphi_j-(L_{i0}-H^{p_i}L_{00})(H^{p_j}\varphi_j) \}. 
\end{equation}
We further introduce the operator $P\mbox{\boldmath$\varphi$}' = (P_1\mbox{\boldmath$\varphi$}',\ldots,
P_N\mbox{\boldmath$\varphi$}')^{\rm T}$. 
Since $L_{ij}^*=L_{ji}$, we see easily that $P$ is symmetric in $L^2(\mathbf{R}^n)$. 
Moreover, we have the following lemma.

\begin{lemma}\label{t-derivative:lemma 1}
Let $c_0, c_1$ be positive constants. 
There exists a positive constant $C=C(c_0,c_1)$ depending only on $c_0$ and $c_1$ such that 
if $H,\nabla b\in L^{\infty}(\mathbf{R}^n)$ satisfy $H(x)\geq c_0$ and $|\nabla b(x)|\leq c_1$, 
then for any $\delta \in (0,1]$ we have 
\[
(P\mbox{\boldmath$\varphi$}',\mbox{\boldmath$\varphi$}')_{L^2}
\geq C^{-1}( \|\nabla \mbox{\boldmath$\varphi$}'\|^2 + \delta^{-2} \|\mbox{\boldmath$\varphi$}'\|^2 ). 
\]
\end{lemma}

\noindent
{\bf Proof}. \ 
Introducing $\varphi_0 = -\sum_{j=1}^N H^{p_j} \varphi_j$, we have 
\begin{align*}
(P \mbox{\boldmath$\varphi$}',\mbox{\boldmath$\varphi$}')_{L^2}
&= \sum_{i,j=0}^N (L_{ij}\varphi_j,\varphi_i)_{L^2} \\
&= \int_{\mathbf{R}^n} \! {\rm d}x \! \int_0^{H} \biggl\{
  \biggl| \sum_{i=0}^N (z^{p_i} \nabla\varphi_i - p_i z^{p_i-1} \varphi_i \nabla b) \biggr|^2
 + \delta^{-2} \biggl( \sum_{i=0}^N p_i z^{p_i-1} \varphi_i \biggr)^2 \biggr\} {\rm d}z,
\end{align*}
which gives the desired estimate. 
For the details, we refer to R. Nemoto and T. Iguchi \cite{NemotoIguchi2017}. 
\quad$\Box$

\bigskip
Once we obtain such a coercive estimate, by the standard theory of elliptic partial differential 
equations, we can obtain the following lemma.

\begin{lemma}\label{t-derivative:lemma 2}
Let $c_0,M$ be positive constants and $m$ an integer such that $m>n/2+1$. 
There exists a positive constant $C=C(c_0,M,m)$ such that if $\eta$ and $b$ satisfy 
\begin{equation}\label{t-derivative:assumption}
\left\{
 \begin{array}{l}
  \|\eta\|_m + \|b\|_{W^{m,\infty}} \leq M, \\[0.5ex]
  c_0 \leq H(x) = 1 + \eta(x) - b(x) \quad\mbox{for}\quad x\in\mathbf{R}^n,
 \end{array}
\right.
\end{equation}
then for $k = 0,\pm 1,\ldots,\pm (m - 1)$ and $\delta\in(0,1]$ we have 
\begin{equation}\label{t-derivative:elliptic estimate}
\|J_\delta P^{-1} \mbox{\boldmath$G$}'\|_k \leq C\delta^2 \|J_\delta^{-1}\mbox{\boldmath$G$}'\|_k.
\end{equation}
\end{lemma}

\begin{remark}\label{t-derivative:remark}
{\rm 
For the estimation to the time derivate $\partial_t \mbox{\boldmath$\phi$}^\delta$, 
it is sufficient to show the above estimate \eqref{t-derivative:elliptic estimate} 
in the Sobolev space with nonnegative indices. 
However, the estimate with negative indices plays an important role in deriving an error 
estimate between the solutions to the Isobe--Kakinuma model and to the full water wave problem. 
}
\end{remark}

\noindent
{\bf Proof}. \ 
Put $\mbox{\boldmath$\varphi$}' = P^{-1}\mbox{\boldmath$G$}'$. 
Noting that $\|\nabla \mbox{\boldmath$\varphi$}'\|^2 + \delta^{-2} \|\mbox{\boldmath$\varphi$}'\|^2$ 
is equivalent to $\delta^{-2}\|J_\delta \mbox{\boldmath$\varphi$}'\|^2$ uniformly with respect to 
$\delta \in (0,1]$, we see by Lemma \ref{t-derivative:lemma 1} that 
\[
\|J_\delta \mbox{\boldmath$\varphi$}'\|^2
\lesssim \delta^2 (P \mbox{\boldmath$\varphi$}',\mbox{\boldmath$\varphi$}')_{L^2}
= \delta^2 (\mbox{\boldmath$G$}',\mbox{\boldmath$\varphi$}')_{L^2}
\leq \delta^2 \|J_\delta^{-1} \mbox{\boldmath$G$}'\| \|J_\delta \mbox{\boldmath$\varphi$}'\|,
\]
which yields the estimate \eqref{t-derivative:elliptic estimate} in the case $k=0$.

Let $1 \leq k \leq m-1$ and $\alpha$ be a multi-index such that $|\alpha| \leq k$. 
Applying the differential operator $\partial^\alpha$ to the equation 
$P \mbox{\boldmath$\varphi$}' = \mbox{\boldmath$G$}'$, we have 
$P \partial^\alpha \mbox{\boldmath$\varphi$}'
 = \partial^\alpha \mbox{\boldmath$G$}' - [\partial^\alpha,P]\mbox{\boldmath$\varphi$}'$, 
so that 
\[
\|J_\delta \partial^\alpha \mbox{\boldmath$\varphi$}'\|
 \lesssim \delta^2 ( \|J_\delta^{-1} \partial^\alpha \mbox{\boldmath$G$}'\|
  + \|J_\delta^{-1} [\partial^\alpha,P]\mbox{\boldmath$\varphi$}'\| ).
\]
We evaluate the commutator $[\partial^\alpha,P]$ by writing down explicitly the operator $P$. 
Let $t_0 > n/2$ and remember the standard commutator estimate 
\[
\|[\partial^\alpha,u]v\|
 \lesssim 
\begin{cases}
  \|u\|_{W^{|\alpha|,\infty}} \|v\|_{|\alpha| - 1}, \\
  \|u\|_{|\alpha| \vee t_0 + 1 } \|v\|_{|\alpha| - 1}.
\end{cases}
\]
By expanding the commutator $[\partial^\alpha,u]v = \partial^\alpha(uv) - u\partial^\alpha v$, 
evaluating each terms separately, and using the calculus inequalities 
$\|uv\|_k \lesssim \|u\|_{|k| \vee t_0}\|u\|_k$ and 
$\|uv\|_k \lesssim \|u\|_{W^{|k|,\infty}}\|v\|_k$ for any integers $k$, 
we also have $\|[\partial^\alpha,u]v\| \lesssim \|u\|_{|\alpha| \vee t_0} \|v\|_{|\alpha|}$ and 
\begin{align*}
\|[\partial^\alpha,u]v\|_{-1}
 \lesssim 
\begin{cases}
  \|u\|_{W^{|\alpha|-1 \vee 1,\infty}} \|v\|_{|\alpha| - 1}, \\
  \|u\|_{|\alpha| - 1\vee 1 \vee t_0 } \|v\|_{|\alpha| - 1}.
\end{cases}
\end{align*}
In the following, we use these calculus inequalities without any comment. 
We also note that we need to handle a smooth function $f(H)$ of $H=1+\eta-b$. 
Under the conditions in \eqref{t-derivative:assumption}, $f(H)$ does not belong to 
$H^m$ nor $W^{m,\infty}$, in general. 
However, we can decompose it as $f(H) = f(1-b) + f_1(\eta,b)\eta$ with a smooth function $f_1$, 
and the first term belongs to $W^{m,\infty}$ and the second one to $H^m$. 
We will also use this fact without any comment. 
Noting 
\begin{equation}\label{t-derivative:multiplier}
\delta \|J_\delta^{-1} \nabla u\| \leq \|u\|, \quad
\delta \|J_\delta^{-1} u\| \leq \|u\|_{-1}, \quad
\|J_\delta^{-1} u\| \leq \|u\|,
\end{equation}
and using the above calculus inequalities, we see that 
\[
\delta^2 \|J_\delta^{-1} [\partial^\alpha,P]\mbox{\boldmath$\varphi$}'\|
\lesssim \delta \|\nabla \mbox{\boldmath$\varphi$}'\|_{k-1} + \|\mbox{\boldmath$\varphi$}'\|_{k-1}
\lesssim \|J_\delta \mbox{\boldmath$\varphi$}'\|_{k-1},
\]
so that 
$\|J_\delta \mbox{\boldmath$\varphi$}'\|_k
 \lesssim \delta^2 \|J_\delta^{-1} \mbox{\boldmath$G$}'\|_k
 + \|J_\delta \mbox{\boldmath$\varphi$}'\|_{k-1}$,
which yields the estimate \eqref{t-derivative:elliptic estimate} for positive $k$ by induction on $k$.

The estimate for negative $k$ comes from the standard duality argument. 
We note that the operator $P$ is symmetric in $L^2(\mathbf{R}^n)$ so is $P^{-1}$. 
Let $1 \leq k' \leq m$. 
Then, we see that 
\begin{align*}
|(J_\delta P^{-1}\mbox{\boldmath$G$}', \mbox{\boldmath$F$}')_{L^2}| 
&= |(J_\delta^{-1} \mbox{\boldmath$G$}', J_\delta P^{-1} J_\delta \mbox{\boldmath$F$}')_{L^2}| 
 \leq \|J_\delta^{-1} \mbox{\boldmath$G$}'\|_{-k'} \|J_\delta P^{-1} J_\delta \mbox{\boldmath$F$}'\|_{k'} \\
&\lesssim \delta^2 \|J_\delta^{-1} \mbox{\boldmath$G$}'\|_{-k'}
  \|J_\delta^{-1}(J_\delta \mbox{\boldmath$F$}')\|_{k'}
 = \delta^2 \|J_\delta^{-1} \mbox{\boldmath$G$}'\|_{-k'} \|\mbox{\boldmath$F$}'\|_{k'},
\end{align*}
which gives 
$\|J_\delta P^{-1}\mbox{\boldmath$G$}'\|_{-k'} 
 \lesssim \delta^2 \|J_\delta^{-1} \mbox{\boldmath$G$}'\|_{-k'}$. 
This gives the estimate \eqref{t-derivative:elliptic estimate} for negative $k$. 
\quad$\Box$

\bigskip
Thanks of Lemma \ref{t-derivative:lemma 2}, a unique existence of the solution $\mbox{\boldmath$\varphi$}$ 
to \eqref{t-derivative:elliptic 1} is guaranteed in appropriate function spaces. 
Concerning estimates of the solution, we have the following lemma.

\begin{lemma}\label{t-derivative:lemma 3}
Let $c_0,M$ be positive constants and $m$ an integer such that $m>n/2+1$. 
There exists a positive constant $C=C(c_0,M,m)$ such that if $\eta$ and $b$ satisfy the conditions 
in \eqref{t-derivative:assumption} and if $\mbox{\boldmath$\varphi$}$ is a solution of 
\eqref{t-derivative:elliptic 1}, then for $k = 0,\pm 1,\ldots,\pm (m - 1)$ and $\delta\in(0,1]$ we have 
\begin{equation}\label{t-derivative:elliptic estimate 2}
\begin{cases}
 \|\nabla \varphi_0\|_k + \delta^{-1} \|J_\delta \mbox{\boldmath$\varphi$}'\|_k 
  \leq C( \|\nabla F_0\|_k + \delta \|J_\delta^{-1} \mbox{\boldmath$F$}'\|_k ), \\
 \|\mbox{\boldmath$\varphi$}\|_{k+1} 
  \leq C( \|F_0\|_{k+1} + \delta \|J_\delta^{-1} \mbox{\boldmath$F$}'\|_k ),
\end{cases}
\end{equation}
where $\mbox{\boldmath$F$}' = (F_1,\ldots,F_N)^{\rm T}$.

If, in addition, $F_0 = 0$, then we have 
\begin{equation}\label{t-derivative:elliptic estimate 3}
\|\mbox{\boldmath$\varphi$}\|_k \leq C \delta^2\|\mbox{\boldmath$F$}'\|_k.
\end{equation}
\end{lemma}

\noindent
{\bf Proof}. \ 
Since $\mbox{\boldmath$\varphi$}'$ satisfy \eqref{t-derivative:elliptic 2}, 
it follows from Lemma \ref{t-derivative:lemma 2} that 
\[
\|J_\delta \mbox{\boldmath$\varphi$}'\|_k
\lesssim \delta^2\|J_\delta^{-1}\mbox{\boldmath$F$}'\|_k
 + \sum_{i=1}^N \delta^2\|J_\delta^{-1} (L_{i0} - H^{p_i}L_{00})F_0\|_k.
\]
Here, by writing down the operator $L_{i0} - H^{p_i}L_{00}$ explicitly and noting 
\eqref{t-derivative:multiplier}, we have 
\[
\delta^2\|J_\delta^{-1} (L_{i0} - H^{p_i}L_{00})F_0\|_k
\lesssim \delta \|\nabla F_0\|_k,
\]
so that 
$\|J_\delta \mbox{\boldmath$\varphi$}'\|_k
 \lesssim \delta^2\|J_\delta^{-1}\mbox{\boldmath$F$}'\|_k + \delta \|\nabla F_0\|_k$. 
Now, we estimate $\varphi_0$ by using \eqref{t-derivative:reduction} and obtain 
$\|\nabla\varphi_0\|_k \lesssim \|\nabla F_0\|_k + \|\mbox{\boldmath$\varphi$}'\|_{k+1}
 \lesssim \|\nabla F_0\|_k + \delta^{-1}\|J_\delta \mbox{\boldmath$\varphi$}'\|_k$. 
Therefore, we obtain the first estimate in \eqref{t-derivative:elliptic estimate 2}. 
Similarly, by \eqref{t-derivative:reduction} we also have 
$\|\varphi_0\|_{k+1} \lesssim \|F_0\|_{k+1} + \delta^{-1}\|J_\delta \mbox{\boldmath$\varphi$}'\|_k$. 
In view of $\|u\|_{k+1} \leq \delta^{-1}\|J_\delta u\|_k$, we obtain the second estimate in 
\eqref{t-derivative:elliptic estimate 2}.

If $F_0=0$, then it follows from the first estimate in \eqref{t-derivative:elliptic estimate 2} that 
$\|\mbox{\boldmath$\varphi$}'\|_k \leq \|J_\delta \mbox{\boldmath$\varphi$}'\|_k
 \lesssim \delta^2 \|J_\delta^{-1}\mbox{\boldmath$F$}'\|_k \lesssim \delta^2 \|\mbox{\boldmath$F$}'\|_k$. 
This together with \eqref{t-derivative:reduction} gives the estimate for $\varphi_0$. 
\quad$\Box$

\bigskip
Now, we are ready to give an estimate for the time derivative $\partial_t \mbox{\boldmath$\phi$}^\delta$. 
We introduce a mathematical energy $E_m(t)$ by 
\begin{equation}\label{t-derivative:energy}
E_m(t) = \|\eta(t)\|_m^2 + \|\nabla \mbox{\boldmath$\phi$}^\delta(t)\|_m^2
 + \delta^{-2} \|\mbox{\boldmath$\phi$}^{\delta\prime}(t)\|_m^2,
\end{equation}
where $\mbox{\boldmath$\phi$}^{\delta\prime} = (\phi_1^\delta,\ldots,\phi_N^\delta)^{\rm T}$.

\begin{lemma}\label{t-derivative:lemma 4}
Let $c_0,M$ be positive constants and $m$ an integer such that $m>n/2+1$. 
There exists a positive constant $C=C(c_0,M,m)$ such that if $(\eta,\mbox{\boldmath$\phi$}^\delta)$ is 
a solution to the Isobe--Kakinuma model \eqref{result:IK model} satisfying 
\begin{equation}\label{t-derivative:assumption 2}
\begin{cases}
 E_m(t) + \|b\|_{W^{m+1,\infty}} \leq M, \\
 c_0 \leq H(x,t) = 1 + \eta(x,t) - b(x) \quad\mbox{for}\quad x \in \mathbf{R}^n, \; 0 \leq t\leq T,
\end{cases}
\end{equation}
then we have 
$\|\partial_t \eta(t)\|_{m-1}^2 + \|\partial_t \mbox{\boldmath$\phi$}^\delta(t)\|_m^2
 + \delta^{-2} \|\partial_t \mbox{\boldmath$\phi$}^{\delta \prime}(t)\|_{m-1}^2 \leq C E_m(t)$
for $0 \leq t\leq T$. 
\end{lemma}

\noindent
{\bf Proof}. \ 
We remind that $\mbox{\boldmath$u$}$ and $w$ were defined by \eqref{t-derivative:uw}, so that 
we easily have $\|\mbox{\boldmath$u$}\|_m^2 + \delta^2 \|w\|_m^2 \lesssim E_m(t)$. 
We remind also that $\partial_t \mbox{\boldmath$\phi$}^\delta$ satisfies \eqref{t-derivative:eq phi} 
where $F_0$ and $\mbox{\boldmath$F$}' = (F_1,\ldots,F_N)^{\rm T}$ are defined by 
\eqref{t-derivative:F0} and \eqref{t-derivative:Fi}, respectively. 
By using the explicit expressions, we see that 
$\|F_0\|_m^2 + \delta^2 \|\mbox{\boldmath$F$}'\|_m^2 \lesssim E_m(t)$. 
Therefore, applying the second estimate in Lemma \ref{t-derivative:lemma 3} and noting 
$\|J_\delta^{-1} u\|_k \leq \|u\|_k$ we obtain 
$\|\partial_t \mbox{\boldmath$\phi$}^\delta\|_m^2 
\lesssim \|F_0\|_m^2 + \delta^2 \|\mbox{\boldmath$F$}'\|_{m-1}^2 \lesssim E_m(t)$. 
On the other hand, applying the first estimate in Lemma \ref{t-derivative:lemma 3} and noting 
$\|u\|_k \leq \|J_\delta u\|_k$ we obtain 
$\delta^{-2} \|\partial_t \mbox{\boldmath$\phi$}^{\delta \prime}\|_{m-1}^2 
\lesssim \|\nabla F_0\|_{m-1}^2 + \delta^2 \|\mbox{\boldmath$F$}'\|_{m-1}^2 \lesssim E_m(t)$. 
The estimate for $\partial_t \eta$ follows directly from the first equation in 
\eqref{result:IK model} with $i=0$. 
\quad$\Box$

\section{Uniform estimate of the solution I}
\label{section:estimate I}
\setcounter{equation}{0}
\setcounter{theorem}{0}
In this section we will prove the first half of Theorem \ref{result:theorem 1}, that is, the existence 
of the solution on a time interval independent of $\delta \in (0,1]$ and a uniform bound 
\eqref{result:uniform estimate 1} to the rescaled variables $(\eta,\mbox{\boldmath$\phi$}^\delta)$ by 
using an energy method.

Let $\alpha$ be a multi-index satisfying $1 \leq |\alpha| \leq m$. 
Applying $\partial^\alpha$ to the Isobe--Kakinuma model \eqref{result:IK model}, 
after a tedious but straightforward calculation, we obtain 
\begin{equation}\label{estimate I:IK model}
\left\{
 \begin{array}{l}
  \displaystyle
  H^{p_i} (( \partial_t + \mbox{\boldmath$u$} \cdot \nabla ) \partial^\alpha \eta)
     - \sum_{j=0}^N L_{ij} (\partial^\alpha \phi_j^\delta)  = - f_{i,\alpha}
   \quad\mbox{for}\quad i=0,1,\ldots,N, \\
  \displaystyle
  \sum_{j=0}^N H^{p_j} ( (\partial_t + \mbox{\boldmath$u$} \cdot \nabla) \partial^\alpha \phi_j^\delta )
   + a \partial^\alpha \eta = f_{N+1,\alpha},
 \end{array}
\right.
\end{equation}
where 
\begin{align}\label{estimate I:fi}
f_{i,\alpha}
&=  [\partial^\alpha,H^{p_i}] \partial_t \eta
 + ((\nabla \cdot (H^{p_i} \mbox{\boldmath$u$})) \partial^\alpha \eta \\
&\quad\;
 + \sum_{j=0}^N \biggl\{  \nabla \cdot \biggl\{
  \biggl( [\partial^\alpha,\frac{1}{p_i+p_j+1} H^{p_i+p_j+1}] - H^{p_i+p_j} (\partial^\alpha \eta) \biggr) 
   \nabla \phi_j^\delta \nonumber \\
&\phantom{ \quad\; + \sum_{j=0}^N \biggl\{ \nabla \cdot \biggl\{ }
 + \biggl( [\partial^\alpha,\frac{p_j}{p_i+p_j} H^{p_i+p_j}(\nabla b)]
   - p_j H^{p_i+p_j-1}(\nabla b) (\partial^\alpha \eta) \biggr) \phi_j^\delta \biggr\} \nonumber \\
&\qquad
 + \frac{p_j}{p_i+p_j} [\partial^\alpha, H^{p_i+p_j} \nabla b] \cdot \nabla\phi_j^\delta
 - \frac{p_ip_j}{p_i+p_j-1} [\partial^\alpha,  H^{p_i+p_j-1}(\delta^{-2}+|\nabla b|^2)] \phi_j^\delta
 \biggr\}, \nonumber
\end{align}
\begin{align}\label{estimate I:fN+1}
f_{N+1,\alpha}
&= - \sum_{j=1}^N\Bigl( [\partial^\alpha, H^{p_j}] - p_j H^{p_j-1} (\partial^\alpha \eta) \Bigr) 
   \partial_t \phi_j^\delta \\
&\quad\;
 - \frac12 ( \partial^\alpha ( |\mbox{\boldmath$u$}|^2 )
   - 2\mbox{\boldmath$u$} \cdot \partial^\alpha \mbox{\boldmath$u$} )
 - \frac12 \delta^2 ( \partial^\alpha ( w^2 ) - 2 w \partial^\alpha w ) \nonumber \\
&\quad\;
 - \mbox{\boldmath$u$} \cdot \sum_{j=1}^N \Bigl\{
  \Bigl( [\partial^\alpha, H^{p_j}] - p_j H^{p_j-1}(\partial^\alpha \eta) \Bigr) \nabla \phi_j^\delta 
  \nonumber \\
&\qquad\qquad
 - p_j \Bigl( [\partial^\alpha, H^{p_j-1}(\nabla b)] - (p_j-1)H^{p_j-2}(\nabla b)(\partial^\alpha \eta)
  \Bigr) \phi_j^\delta \Bigr\} \nonumber \\
&\quad\;
 - w \cdot \sum_{j=1}^N p_j \Bigl( [\partial^\alpha, H^{p_j-1}] - (p_j-1)H^{p_j-2}(\partial^\alpha \eta)
  \Bigr) \phi_j^\delta, \nonumber
\end{align}
and $a$ is related to a generalized Rayleigh--Taylor sign condition and is given by 
\begin{align}\label{estimate I:a}
a &= 1 + \sum_{j=1}^N p_j H^{p_j-1} \partial_t \phi_j^\delta \\
&\quad\;
 + \mbox{\boldmath$u$} \cdot \sum_{j=1}^N \Bigl( p_j H^{p_j-1} \nabla\phi_j^\delta
  - p_j(p_j-1) H^{p_j-2} \phi_j^\delta \nabla b \Bigr) 
 + w \sum_{j=1}^N p_j(p_j-1) H^{p_j-2} \phi_j^\delta. \nonumber
\end{align}
We can rewrite \eqref{estimate I:IK model} in a matrix form as 
\begin{equation}\label{estimate I:IK model mat}
\begin{pmatrix}
 0 & \mbox{\boldmath$l$}^{\rm T} \\
 -\mbox{\boldmath$l$} & O
\end{pmatrix}
( \partial_t + \mbox{\boldmath$u$} \cdot \nabla ) \partial^\alpha
\begin{pmatrix}
 \eta \\
 \mbox{\boldmath$\phi$}^\delta
\end{pmatrix}
+
\begin{pmatrix}
 a & \mbox{\boldmath$0$}^{\rm T} \\
 \mbox{\boldmath$0$} & L
\end{pmatrix}
\partial^\alpha
\begin{pmatrix}
 \eta \\
 \mbox{\boldmath$\phi$}^\delta
\end{pmatrix}
=
\begin{pmatrix}
 f_{N+1,\alpha} \\
 \mbox{\boldmath$f$}_\alpha
\end{pmatrix},
\end{equation}
where $\mbox{\boldmath$f$}_\alpha = (f_{0,\alpha},\ldots,f_{N,\alpha})^{\rm T}$, 
$L = (L_{ij})_{0 \leq i,j \leq N}$, and 
\begin{equation}\label{estimate I:l}
\mbox{\boldmath$l$} = \mbox{\boldmath$l$}(H) = (H^{p_0},\ldots,H^{p_N})^{\rm T}.
\end{equation}
Since $L_{ij}^* = L_{ji}$, the matrix operator $L$ acting on 
$\mbox{\boldmath$\varphi$} = (\varphi_0,\ldots,\varphi_N)^{\rm T}$ is symmetric in $L^2(\mathbf{R}^n)$. 
Moreover, we have already shown the positivity of $L$ in the proof of Lemma \ref{t-derivative:lemma 1}, 
that is, we have the following lemma.

\begin{lemma}\label{estimate I:lemma 1}
Let $c_1,C_0$ be positive constants. 
There exists a positive constant $C=C(c_1,C_0)$ such that 
if $H,\nabla b\in L^{\infty}(\mathbf{R}^n)$ satisfy $C_0^{-1} \leq H(x) \leq C_0$ and $|\nabla b(x)|\leq c_1$, 
then for any $\delta \in (0,1]$ we have 
\[
C^{-1}( \|\nabla \mbox{\boldmath$\varphi$}\|^2 + \delta^{-2} \|\mbox{\boldmath$\varphi$}'\|^2 )
\leq (L\mbox{\boldmath$\varphi$},\mbox{\boldmath$\varphi$})_{L^2}
\leq C ( \|\nabla \mbox{\boldmath$\varphi$}\|^2 + \delta^{-2} \|\mbox{\boldmath$\varphi$}'\|^2 )
\]
where $\mbox{\boldmath$\varphi$}' = (\varphi_1,\ldots,\varphi_N)^{\rm T}$. 
\end{lemma}

By making use of this nice structure of the equations, we derive an energy estimate 
which leads uniform bound of the solution in the rescaled variables. 
Before carrying out the estimate, we need to show that an appropriate norm of the right-hand side of 
\eqref{estimate I:IK model mat} would be evaluated by our energy function $E_m(t)$ 
uniformly with respect to $\delta \in (0,1]$. 
However, it contains a term $\delta^{-2}[\partial^\alpha, H^{p_i+p_j-1}] \phi_j^\delta$, 
which cannot be estimated directly because of the coefficient $\delta^{-2}$. 
Nevertheless, thanks of the commutator we can gain a regularity of order one. 
Using this fact and the necessary conditions \eqref{t-derivative:compatibility}, 
we can handle such a term. 
We remind that the necessary conditions can be written simply as 
$\mathscr{L}_i\mbox{\boldmath$\phi$}^\delta = 0$ for $i=1,\ldots,N$.

\begin{lemma}\label{estimate I:lemma 2}
Let $c_0,M$ be positive constants and $m$ an integer such that $m>n/2+1$. 
There exists a positive constant $C=C(c_0,M,m)$ such that if $\eta$ and $b$ satisfy 
\begin{equation}\label{estimate I:assumption}
\left\{
 \begin{array}{l}
  \|\eta\|_{m-1} + \|b\|_{W^{m+1,\infty}} \leq M, \\[0.5ex]
  c_0 \leq H(x) = 1 + \eta(x) - b(x) \quad\mbox{for}\quad x\in\mathbf{R}^n,
 \end{array}
\right.
\end{equation}
and if $\mbox{\boldmath$\varphi$}$ satisfies $\mathscr{L}_i\mbox{\boldmath$\varphi$} = F_i$ 
for $i=1,\ldots,N$, then for $k = 0,\pm 1,\ldots,\pm (m-1)$ and $\delta\in(0,1]$ we have 
\[
\delta^{-2} \|\mbox{\boldmath$\varphi$}^{\prime}\|_k
 \leq C ( \|\nabla \mbox{\boldmath$\varphi$}\|_{k+1} + \|\mbox{\boldmath$\varphi$}'\|_{k+1}
 + \|\mbox{\boldmath$F$}'\|_k ).
\]
\end{lemma}

\noindent
{\bf Proof}. \ 
In view of \eqref{t-derivative:Li}, 
we see that the equation $\mathscr{L}_i\mbox{\boldmath$\varphi$} = F_i$ is equivalent to 
\begin{align}\label{estimate I:Li}
\sum_{j=1}^N \frac{p_ip_j}{p_i+p_j-1} H^{p_j} \varphi_j
&= \frac{\delta^2}{1+\delta^2|\nabla b|} \sum_{j=0}^N H^{p_j+1} \biggl\{
 - \frac{p_i}{(p_i+p_j+1)(p_j+1)} H\Delta\varphi_j \\
&\qquad\qquad
 + \frac{p_ip_j}{(p_i+p_j)p_j} \nabla \cdot (\varphi_j \nabla b)
 + \frac{p_i}{(p_i+p_j)} \nabla b \cdot \nabla \varphi_j \biggr\} \nonumber \\
&\quad\;
 + \frac{\delta^2}{1+\delta^2|\nabla b|} H^{1-p_i} F_i \nonumber
\end{align}
for $i=1,\ldots,N$. 
Since $N \times N$ matrix $A_1' = \bigl( \frac{p_ip_j}{p_i+p_j-1} \bigr)_{1 \leq i,j \leq N}$ 
is nonsingular, the desired estimate comes from standard calculus inequalities. 
\quad$\Box$

\begin{lemma}\label{estimate I:lemma 3}
Let $c_0,M$ be positive constants and $m$ an integer such that $m>n/2+1$. 
There exists a positive constant $C=C(c_0,M,m)$ such that if $(\eta,\mbox{\boldmath$\phi$}^\delta)$ is 
a solution to the Isobe--Kakinuma model \eqref{result:IK model} satisfying the conditions in 
\eqref{t-derivative:assumption 2}, then we have 
\[
\begin{cases}
 \delta^{-4} ( \|\mbox{\boldmath$\phi$}^{\delta \prime}(t)\|_{m-1}^2
  + \|\partial_t \mbox{\boldmath$\phi$}^{\delta \prime}(t)\|_{m-2}^2 ) \leq C E_m(t), \\
 \|\mbox{\boldmath$f$}_\alpha\|^2 + \|f_{N+1,\alpha}\|_1^2 \leq C E_m(t).
\end{cases}
\] 
\end{lemma}

\noindent
{\bf Proof}. \ 
By Lemma \ref{estimate I:lemma 2}, we have 
$\delta^{-2} \|\mbox{\boldmath$\phi$}^{\delta \prime}\|_{m-1} 
 \lesssim \|\nabla \mbox{\boldmath$\phi$}^{\delta}\|_m
  + \|\mbox{\boldmath$\phi$}^{\delta \prime}\|_m
 \lesssim E_m(t)^{1/2}$. 
Let $\mbox{\boldmath$F$}' = (F_1,\ldots,F_N)^{\rm T}$ be defined by \eqref{t-derivative:Fi}. 
Then, we have 
$\|\mbox{\boldmath$F$}'\|_{m-1}
 \lesssim \|\nabla \mbox{\boldmath$\phi$}^{\delta}\|_m
  + \delta^{-2} \|\mbox{\boldmath$\phi$}^{\delta \prime}\|_{m-1}
 \lesssim  E_m(t)^{1/2}$. 
Since $\partial_t \mbox{\boldmath$\phi$}^\delta$ satisfies \eqref{t-derivative:eq phi}, 
by Lemmas \ref{estimate I:lemma 2} and \ref{t-derivative:lemma 4}, we get 
$\delta^{-2} \|\partial_t \mbox{\boldmath$\phi$}^{\delta \prime}\|_{m-2}
 \lesssim \|\partial_t \mbox{\boldmath$\phi$}^\delta\|_m
  + \|\mbox{\boldmath$F$}'\|_{m-2}
 \lesssim  E_m(t)^{1/2}$. 
Therefore, we obtain the first estimate of the lemma. 
Note that we also have 
$\|\partial_t \eta\|_{m-1} + \|\mbox{\boldmath$u$}\|_m + \delta \|w\|_m \lesssim E_m(t)^{1/2}$. 
Thus, by using the standard commutator estimate and an estimate for a symmetric commutator 
$\|\partial^\alpha(uv) - (\partial^\alpha u)v - u(\partial^\alpha v)\|_1
 \lesssim \|u\|_{|\alpha| \vee t_0+1} \|v\|_{|\alpha| \vee t_0+1}$, 
we obtain the second estimate of the lemma. 
\quad$\Box$

\bigskip
In our energy estimate, we need to handle the time derivative $\partial_t a$. 
Since the coefficient $a$ contains $\partial_t \mbox{\boldmath$\phi$}^{\delta \prime}$, 
we need to estimate the second order time derivative $\partial_t^2 \mbox{\boldmath$\phi$}^{\delta \prime}$.

\begin{lemma}\label{estimate I:lemma 4}
Let $c_0,M$ be positive constants and $m$ an integer such that $m>n/2+1$. 
There exists a positive constant $C=C(c_0,M,m)$ such that if $(\eta,\mbox{\boldmath$\phi$}^\delta)$ is 
a solution to the Isobe--Kakinuma model \eqref{result:IK model} satisfying the conditions in 
\eqref{t-derivative:assumption 2}, then we have 
\[
\begin{cases}
 \|\partial_t^2 \eta(t)\|_{m-2}^2 + \|\partial_t^2 \mbox{\boldmath$\phi$}^\delta (t)\|_{m-1}^2
  + \delta^{-2} \|\partial_t^2 \mbox{\boldmath$\phi$}^\delta (t)\|_{m-2}^2  \leq C E_m(t), \\
 \|a-1\|_m^2 +\delta^{-2} \|a-1\|_{m-1}^2 + \|\partial_t a\|_{m-1}^2 \leq C E_m(t).
\end{cases}
\]
\end{lemma}

\noindent
{\bf Proof}. \ 
Differentiating the first equation in \eqref{t-derivative:IK model} (equivalently \eqref{result:IK model}) 
with $i=0$ with respect to $t$, we have 
\[
\partial_t^2 \eta = \sum_{j=0}^N L_{0j} \partial_t \phi_j^\delta
 - \nabla \cdot ((\partial_t \eta)\mbox{\boldmath$u$}),
\]
which together with Lemma \ref{t-derivative:lemma 4} yields 
$\|\partial_t^2 \eta\|_{m-2}^2 \lesssim E_m(t)$. 
Note that the operator $\mathscr{L}$ depends on $H$ but not on $\nabla H$. 
Therefore, differentiating the second equation in \eqref{t-derivative:IK model} 
and the necessary conditions $\mathscr{L}_i \mbox{\boldmath$\phi$}^\delta = 0$ for $i=1,\ldots,N$ twice 
with respect to $t$, we have 
\[
\mathscr{L} \partial_t^2 \mbox{\boldmath$\phi$}^\delta = \mbox{\boldmath$F$}_1,
\]
where $\mbox{\boldmath$F$}_1 = (F_{0,1},\ldots,F_{N,1})^{\rm T}$, and
\[
\begin{cases}
 \displaystyle
 F_{0,1} = - \partial_t \eta - (\partial_t \eta)\sum_{j=1}^N p_j H^{p_j-1} \partial_t \phi_j^\delta
   - \mbox{\boldmath$u$} \cdot \partial_t \mbox{\boldmath$u$} - \delta^2 w \partial_t w, \\
 \displaystyle
 F_{i,1} = - (\partial_t^2 \eta)\biggl( \frac{\partial}{\partial H} \mathscr{L}_i \biggr)
   \mbox{\boldmath$\phi$}^\delta
  - (\partial_t \eta)^2\biggl( \frac{\partial^2}{\partial H^2} \mathscr{L}_i \biggr)
   \mbox{\boldmath$\phi$}^\delta
  - 2(\partial_t \eta)\biggl( \frac{\partial}{\partial H} \mathscr{L}_i \biggr)
   \partial_t \mbox{\boldmath$\phi$}^\delta
\end{cases}
\]
for $i=1,\ldots,N$. 
Here, we note that $\bigl( \frac{\partial}{\partial H} \bigr)^j \mathscr{L}_i$ is also a second 
order differential operators like $\mathscr{L}_i$. 
Therefore, by Lemmas \ref{t-derivative:lemma 4} and \ref{estimate I:lemma 3} we have 
$\|\partial_t \mbox{\boldmath$u$}\|_{m-1}^2 + \delta^2 \|\partial_t w\|_{m-1}^2
 \lesssim E_m(t)$ and 
$\|F_{0,1}\|_{m-1}^2 + \|\mbox{\boldmath$F$}_1'\|_{m-2}^2 \lesssim E_m(t)$, 
where $\mbox{\boldmath$F$}_1' = (F_{1,1},\ldots,F_{N,1})^{\rm T}$. 
Applying the both estimates in Lemma \ref{t-derivative:lemma 3}, we obtain 
$\|\partial_t^2 \mbox{\boldmath$\phi$}^\delta\|_{m-1}
  + \delta^{-1} \|\partial_t^2 \mbox{\boldmath$\phi$}^\delta\|_{m-2}
 \lesssim \|F_{0,1}\|_{m-1} + \delta \|J_\delta^{-1} \mbox{\boldmath$F$}_1'\|_{m-2}
 \lesssim E_m(t)$, so that the first estimate of the lemma is proved. 
In view of \eqref{estimate I:a}, the above estimates together with Lemmas \ref{t-derivative:lemma 4} 
and \ref{estimate I:lemma 3} yield the second estimate of the lemma. 
\quad$\Box$

\bigskip
Now, we are ready to give a proof of the first half of Theorem \ref{result:theorem 1}. 
Since the existence theorem has already been established by R. Nemoto and T. Iguchi \cite{NemotoIguchi2017} 
in the function spaces, it is sufficient to show \eqref{result:uniform estimate 1} for some time interval 
independent of $\delta \in (0,1]$. 
Moreover, in view of Lemmas \ref{t-derivative:lemma 4}, \ref{estimate I:lemma 3}, 
and \ref{estimate I:lemma 3}, it is sufficient to show that 
\begin{equation}\label{estimate I:estimate}
E_m(t) \leq M_1, \quad c_0/2 \leq H(x,t) \leq 2C_0,\quad c_0/2 \leq a(x,t) \leq 2C_0 
\end{equation}
for any $x\in\mathbf{R}^n$, $0\leq t\leq T$, and $0<\delta \leq 1$, where $C_0$ is chosen so that 
$H(x,0) \leq C_0$ and $a(x,0) \leq C_0$ and the constant $M_1$ and the time $T$ will be determined later. 
Note that by Lemma \ref{estimate I:lemma 4} such a constant $C_0$ exists under our assumption on the initial 
data and the bottom topography. 
In the following we simply write the constants depending only on $(c_0,C_0,M_0,m)$ by $C_1$ and the 
constants depending also on $M_1$ by $C_2$, which may change from line to line.

We remind that the solution satisfies \eqref{estimate I:IK model mat}. 
In view of this symmetric form of the equations, we introduce an energy function $\mathscr{E}_m(t)$ by 
\[
\mathscr{E}_m(t) = \sum_{|\alpha| \leq m} \{
 (a\partial^\alpha \eta(t),\partial^\alpha \eta(t))_{L^2}
 + (L\partial^\alpha \mbox{\boldmath$\phi$}^\delta(t),
  \partial^\alpha \mbox{\boldmath$\phi$}^\delta(t))_{L^2} \}.
\]
Now, suppose that the solution satisfies \eqref{estimate I:estimate}. 
Then, by Lemma \ref{estimate I:lemma 1} we see that 
\begin{equation}\label{estimate I:equivalence}
C_1^{-1} E_m(t) \leq \mathscr{E}_m(t) \leq C_1 E_m(t)
\end{equation}
for $0\leq t\leq T$. 
For $1 \leq |\alpha| \leq m$ we take the $L^2$-inner product of \eqref{estimate I:IK model mat} with 
$(\partial_t + \mbox{\boldmath$u$} \cdot \nabla) \partial^\alpha 
( \eta, \mbox{\boldmath$\phi$}^\delta )^{\rm T}$ and use the symmetry of the operator $L$ and 
integration by parts. 
For $|\alpha|=0$ we evaluate it directly. 
Then, we obtain 
\begin{align}\label{estimate I: energy estimate}
\frac{\rm d}{{\rm d}t} \mathscr{E}_m(t)
&= \sum_{|\alpha| \leq m} \{
 ((\partial_t a) \partial^\alpha \eta,\partial^\alpha \eta)_{L^2}
 + ([\partial_t,L]\partial^\alpha \mbox{\boldmath$\phi$}^\delta,
  \partial^\alpha \mbox{\boldmath$\phi$}^\delta)_{L^2} \} \\
&\quad\;
 + \sum_{1 \leq |\alpha| \leq m} \{
   ((\nabla \cdot (a\mbox{\boldmath$u$})) \partial^\alpha \eta,\partial^\alpha \eta)_{L^2}
  - 2 (L\partial^\alpha \mbox{\boldmath$\phi$}^\delta,
  (\mbox{\boldmath$u$} \cdot \nabla)\partial^\alpha \mbox{\boldmath$\phi$}^\delta)_{L^2} \nonumber \\
&\qquad\qquad\qquad
 + 2 (f_{N+1,\alpha},(\partial_t + \mbox{\boldmath$u$} \cdot \nabla) \partial^\alpha \eta)_{L^2}
 + 2 (\mbox{\boldmath$f$}_\alpha,
  (\partial_t + \mbox{\boldmath$u$} \cdot \nabla) \partial^\alpha \mbox{\boldmath$\phi$}^\delta)_{L^2} \}
  \nonumber \\
&\quad\;
  + 2 (a \eta,\partial_t \eta)_{L^2}
  + 2 (L \mbox{\boldmath$\phi$}^\delta, \partial_t \mbox{\boldmath$\phi$}^\delta)_{L^2}. \nonumber
\end{align}
To evaluate the term with the commutator $[\partial_t,L]$, it is sufficient to see that 
\begin{align*}
([\partial_t,L]\mbox{\boldmath$\varphi$},\mbox{\boldmath$\varphi$})_{L^2}
&= \sum_{i,j=0}^N \int_{\mathbf{R}^n} (\partial_t \eta)\{
 H^{p_i+p_j} \nabla \varphi_j \cdot \nabla \varphi_i
 - p_j H^{p_i+p_j-1} \varphi_j \nabla b \cdot \nabla \varphi_i \\
&\qquad
 - p_i H^{p_i+p_j-1} \varphi_i \nabla b \cdot \nabla \varphi_j 
 +p_ip_j H^{p_i+p_j-2} (\delta^{-2} + |\nabla b|^2) \varphi_j \varphi_i \} {\rm d}x,
\end{align*}
which yields 
$|([\partial_t,L]\mbox{\boldmath$\varphi$},\mbox{\boldmath$\varphi$})_{L^2}|
 \lesssim \|\nabla \mbox{\boldmath$\varphi$}\|^2 + \delta^{-2} \|\mbox{\boldmath$\varphi$}\|^2$. 
To evaluate the term 
$(L\partial^\alpha \mbox{\boldmath$\phi$}^\delta,
  (\mbox{\boldmath$u$} \cdot \nabla)\partial^\alpha \mbox{\boldmath$\phi$}^\delta)_{L^2}$, 
we decompose the operator $L$ into its principal term 
$L^{\rm pr} = (L_{ij}^{\rm pr})_{0 \leq i,j \leq N}$ and the remainder part 
$L^{\rm low} = (L_{ij}^{\rm low})_{0 \leq i,j \leq N}$, where 
\[
\begin{cases}
\displaystyle
L_{ij}^{\rm pr}\varphi_j
= - \nabla \cdot \biggl(
   \frac{1}{p_i+p_j+1} H^{p_i+p_j+1} \nabla\varphi_j \biggr) 
  +\delta^{-2} \frac{p_ip_j}{p_i+p_j-1} H^{p_i+p_j-1} \varphi_j, \\[2ex]
\displaystyle
L_{ij}^{\rm low}\varphi_j
= \nabla \cdot \biggl( \frac{p_j}{p_i+p_j} H^{p_i+p_j} \varphi_j \nabla b \biggr) 
  - \frac{p_i}{p_i+p_j} H^{p_i+p_j} \nabla b \cdot \nabla\varphi_j \\[2ex]
\displaystyle
\phantom{ L_{ij}^{\rm low}\varphi_j = }
   + \frac{p_ip_j}{p_i+p_j-1} H^{p_i+p_j-1} |\nabla b|^2 \varphi_j.
\end{cases}
\]
We can evaluate the term 
$(L^{\rm low} \partial^\alpha \mbox{\boldmath$\phi$}^\delta,
 (\mbox{\boldmath$u$} \cdot \nabla)\partial^\alpha \mbox{\boldmath$\phi$}^\delta)_{L^2}$ 
directly by the Cauchy--Schwarz inequality, whereas 
the term 
$(L^{\rm pr} \partial^\alpha \mbox{\boldmath$\phi$}^\delta,
 (\mbox{\boldmath$u$} \cdot \nabla)\partial^\alpha \mbox{\boldmath$\phi$}^\delta)_{L^2}$ 
is evaluated by the expression 
\begin{align*}
& (L^{\rm pr} \mbox{\boldmath$\varphi$},
 (\mbox{\boldmath$u$} \cdot \nabla) \mbox{\boldmath$\varphi$})_{L^2} \\
&= \sum_{i,j=0}^N \int_{\mathbf{R}^n} \biggl\{
 \frac{1}{p_i+p_j+1} \biggl(H^{p_i+p_j+1} \nabla \varphi_j \cdot
  [\nabla, \mbox{\boldmath$u$} \cdot \nabla ] \varphi_i
 - \frac12( \nabla \cdot (H^{p_i+p_j+1} \mbox{\boldmath$u$}) )
  \nabla\varphi_j \cdot \nabla\varphi_i \biggr) \\
&\qquad\qquad\quad
 - \frac12 \delta^{-2} \frac{p_ip_j}{p_i+p_j-1}
  (\nabla \cdot (H^{p_i+p_j-1} \mbox{\boldmath$u$}) ) \varphi_j \varphi_i
   \biggr\} {\rm d}x,
\end{align*}
where we used integration by parts. 
This yields 
$|(L^{\rm pr} \mbox{\boldmath$\varphi$},
  (\mbox{\boldmath$u$} \cdot \nabla) \mbox{\boldmath$\varphi$})_{L^2}|
 \lesssim \|\nabla \mbox{\boldmath$\varphi$}\|^2 + \delta^{-2} \|\mbox{\boldmath$\varphi$}\|^2$. 
Concerning the terms with $f_{N+1,\alpha}$, by Lemma \ref{estimate I:lemma 3} and 
$\|\mbox{\boldmath$u$}\|_m \lesssim E_m(t)$ we evaluate it as 
\[
|(f_{N+1,\alpha},(\partial_t + \mbox{\boldmath$u$} \cdot \nabla) \partial^\alpha \eta)_{L^2}| 
 \leq \|f_{N+1,\alpha}\|_1 \|(\partial_t + \mbox{\boldmath$u$} \cdot \nabla) \partial^\alpha \eta\|_{-1}
 \leq C_2 E_m(t).
\]
The term with $\mbox{\boldmath$f$}_\alpha$ and the last two terms in the right-hand side of 
\eqref{estimate I: energy estimate} can be evaluated directly by the  Cauchy--Schwarz inequality. 
Therefore, in view of Lemmas \ref{estimate I:lemma 3} and \ref{estimate I:lemma 4} we obtain 
$\frac{\rm d}{{\rm d}t} \mathscr{E}_m(t) \leq C_2 \mathscr{E}_m(t)$, 
which together with Gronwall's inequality and the equivalence \eqref{estimate I:equivalence} implies 
\[
E_m(t) \leq C_1 E_m(0) e^{C_2t} \leq C_1 M_0^2 e^{C_2t}.
\]
On the other hand, by the fundamental theorem of calculus, the Sobolev imbedding theorem, 
and Lemmas \ref{t-derivative:lemma 4} and \ref{estimate I:lemma 4} we have 
\[
|H(x,t) - H(x,0)| + |a(x,t) - a(x,0)| \leq C_2t.
\]
By taking into account these two inequalities, we define the positive constant $M_1$ and the time $T$ 
so that $M_1=2C_1 M_0^2$ and then $T = C_2^{-1}\min\{ \log 2, C_0, c_0/2 \}$. 
Then, the above arguments show that the solution in fact satisfy \eqref{estimate I:estimate} for 
$0\leq t\leq T$ uniformly in $\delta \in (0,1]$.

\section{Uniform estimate of the solution II}
\label{section:estimate II}
\setcounter{equation}{0}
\setcounter{theorem}{0}
In this section we will prove the second half of Theorem \ref{result:theorem 1}, that is, 
the uniform bound \eqref{result:uniform estimate 2} to the original variables 
$(\eta,\mbox{\boldmath$\phi$})$ by using the uniform bound \eqref{result:uniform estimate 1} 
obtained in the previous section and the necessary conditions \eqref{result:compatibility}. 
To this end we have to use the advantage of our specific choice of the indices $p_i$, that is, 
$p_i=2i$ in the case of the flat bottom and $p_i=i$ in the case with general bottom topographies.

\subsection{The case $p_i=2i$ with the flat bottom}

\begin{lemma}\label{estimate II:formula 1}
Choose $p_i = 2i$ $(i = 0,1,\ldots,N)$ and suppose that the bottom is flat. 
If $\mbox{\boldmath$\varphi$} = (\varphi_0,\ldots,\varphi_N)^{\rm T}$ satisfies 
$\mathscr{L}_i \mbox{\boldmath$\varphi$} = 0$ for $i=1,\ldots,N$, then we have 
\[
\varphi_j = \delta^2 \biggl\{ -\frac{1}{2j(2j-1)}\Delta\varphi_{j-1}
 + \beta_{j,N} H^{2(N-j)}\Delta\varphi_N \biggr\}
\]
for $j = 1,\ldots,N$, where the constant $\beta_{j,N}$ is defined by \eqref{estimate II:beta} below.
\end{lemma}

\noindent
{\bf Proof}. \ 
It follows from \eqref{estimate I:Li} with $F_i=0$ that 
\begin{align*}
\sum_{j=1}^N \frac{4ij}{2(i+j)-1} H^{2j} \varphi_j
&= \sum_{j=1}^N \frac{4ij}{2(i+j)-1} \biggl( -\frac{\delta^2}{2j(2j-1)} H^{2j} \Delta\varphi_{j-1} \biggr) \\
&\quad\;
 - \frac{2i}{(2(N+i)+1)(2N+1)} \delta^2 H^{2N} \Delta\varphi_N. 
\end{align*}
In view of this, we define constants $\beta_{j,N}$ for $j=1,2,\ldots,N$ by 
\begin{equation}\label{estimate II:beta}
\sum_{j=1}^N \frac{4ij}{2(i+j)-1} \beta_{j,N} = - \frac{2i}{(2(N+i)+1)(2N+1)}
\end{equation}
for $i=1,2,\ldots,N$. 
Since the matrix $\bigl( \frac{4ij}{2(i+j)-1} )_{1 \leq i,j \leq N}$ is nonsingular, 
the constants $\beta_{j,N}$ $(j=1,2,\ldots,N)$ are uniquely determined. 
Then, we obtain the desired identity. 
\quad$\Box$

\begin{lemma}\label{estimate II:estimate 1}
Under the same hypothesis of Lemma \ref{estimate II:formula 1}, for any integer $k$ we have 
\[
\|(\varphi_j,\ldots,\varphi_N)\|_k
\leq \delta^{2j} C(\|\eta\|_{|k| \vee |k+2(j-1)| \vee t_0}) \|\nabla \mbox{\boldmath$\varphi$}\|_{k+2j-1}
\]
for $j=1,\ldots,N$. 
\end{lemma}

\noindent
{\bf Proof}. \ 
It follows from Lemma \eqref{estimate II:formula 1} that 
$\|\varphi_j\|_k \leq \delta^2 \|\nabla \varphi_{j-1}\|_{k+1}
 + \delta^2 C(\|\eta\|_{|k| \vee t_0})\|\nabla \varphi_N\|_{k+1}$, 
so that 
\begin{align*}
\|(\varphi_j,\ldots,\varphi_N)\|_k
&\leq \delta^2 C(\|\eta\|_{|k| \vee t_0})\|\nabla (\varphi_{j-1},\ldots,\varphi_N)\|_{k+1} \\
&\leq \delta^2 C(\|\eta\|_{|k| \vee t_0})\|(\varphi_{j-1},\ldots,\varphi_N)\|_{k+2}
\end{align*}
for $j=1,\ldots,N$. 
Using this inductively, we obtain the desired estimate. 
\quad$\Box$

\bigskip
Now, we will show the second half of Theorem \ref{result:theorem 1} in the case (H1). 
Since $\mbox{\boldmath$\phi$}^\delta$ satisfies the necessary conditions \eqref{result:compatibility}, 
we can apply Lemma \ref{estimate II:estimate 1} with 
$\mbox{\boldmath$\varphi$} = \mbox{\boldmath$\phi$}^\delta$ and $k = m-2j+1$. 
Then, under our hypothesis we have $m-1 > n/2$ and $m \geq j \geq 1$, so that 
$|k| \vee |k+2(j-1)| \vee t_0 = m-1$. 
Therefore, we obtain 
\[
\|\phi_j^\delta(t)\|_{m-2j+1}
 \leq \delta^{2j} C(\|\eta(t)\|_{m-1}) \|\nabla \mbox{\boldmath$\phi$}^\delta(t)\|_m,
\]
which together with \eqref{result:uniform estimate 1} yields the desired estimate 
\eqref{result:uniform estimate 2} in the case (H1).

\subsection{The case $p_i=i$ with general bottom topographies}
For simplify the description, we introduce a differential operator $Q=Q(b)$ depending on 
the bottom topography $b$ by 
\begin{equation}\label{estimate II:Q}
Q\psi = \nabla \cdot (\psi \nabla b) + \nabla b \cdot \nabla\psi.
\end{equation}

\begin{lemma}\label{estimate II:formula 2}
Choose $p_i = i$ $(i = 0,1,\ldots,N)$. 
If $\mbox{\boldmath$\varphi$} = (\varphi_0,\ldots,\varphi_N)^{\rm T}$ satisfies 
$\mathscr{L}_i \mbox{\boldmath$\varphi$} = 0$ for $i=1,\ldots,N$, then we have 
\[
\begin{cases}
 \displaystyle
 \varphi_1 = \frac{\delta^2}{1+\delta^2|\nabla b|^2} \Bigl\{ \nabla b \cdot \nabla\varphi_0
  + \gamma_{1,N-1} H^N ( \Delta\varphi_{N-1} - N Q(b)\varphi_N ) 
  + \gamma_{1,N} H^{N+1} \Delta\varphi_N \Bigr\}, \\[2ex]
 \displaystyle
 \varphi_j = \frac{\delta^2}{1+\delta^2|\nabla b|^2} \biggl\{
  - \frac{1}{j(j-1)} \Delta\varphi_{j-2} + \frac{1}{j} Q(b)\varphi_{j-1} \\
 \displaystyle
 \makebox[6ex]{}
  + \gamma_{j,N-1} H^{N-j+1} ( \Delta\varphi_{N-1} - N Q(b)\varphi_N)
  + \gamma_{j,N} H^{N-j+2} \Delta\varphi_N \biggr\}
 \quad\mbox{for}\quad j=2,\ldots,N,
\end{cases}
\]
where the constant $\gamma_{j,k}$ is defined by \eqref{estimate II:gamma} below.
\end{lemma}

\noindent
{\bf Proof}. \ 
It follows from \eqref{estimate I:Li} with $F_i=0$ that 
\begin{align*}
&\sum_{j=1}^N \frac{ij}{i+j-1} H^j \varphi_j
= \frac{\delta^2}{1+\delta^2|\nabla b|^2} \biggl\{
 \sum_{j=2}^N \frac{ij}{i+j-1} H^j \biggl( -\frac{1}{j(j-1)} \Delta\phi_{j-2} 
  + \frac{1}{j} \nabla \cdot (\varphi_{j-1} \nabla b) \biggr) \\
&\makebox[12ex]{}
 + \sum_{j=1}^N \frac{ij}{i+j-1} H^j \biggl( \frac{1}{j} \nabla b \cdot \nabla\varphi_{j-1} \biggr) \\
&\makebox[12ex]{}
 - \frac{i}{(N+i)N} H^{N+1} (\Delta\varphi_{N-1} - N Q(b)\varphi_N)
 - \frac{i}{(N+i+1)(N+1)} H^{N+2} \Delta\phi_N^\delta \biggr\}.
\end{align*}
In view of this, we define constants $\gamma_{j,k}$ for $j=1,2,\ldots,N$ and $k \geq 0$ by 
\begin{equation}\label{estimate II:gamma}
\sum_{j=1}^N \frac{ij}{i+j-1} \gamma_{j,k} = -\frac{i}{(k+i+1)(k+1)}
\end{equation}
for $i=1,2,\ldots,N$. 
Since the matrix $\bigl( \frac{ij}{i+j-1} )_{1 \leq i,j \leq N}$ is nonsingular, 
the constants $\gamma_{j,k}$ $(j=1,2,\ldots,N, \; k\geq0)$ are uniquely determined. 
Then, we obtain the desired identity. 
\quad$\Box$

\begin{lemma}\label{estimate II:estimate 2}
Under the same hypothesis of Lemma \ref{estimate II:formula 2}, for any integer $k$ we have 
\begin{align*}
& \|(\varphi_{2j-1},\varphi_{2j},\ldots,\varphi_N)\|_k \\
& \leq \delta^{2j} C( \|\eta\|_{|k| \vee |k+2(j-1)| \vee t_0},\|b\|_{W^{|k|+1 \vee |k+2j-1|+1,\infty}} )
 ( \|\nabla \varphi_0\|_{k+2j-1} + \|\mbox{\boldmath$\varphi$}'\|_{k+2j} )
\end{align*}
for $j=1,\ldots,[(N+1)/2]$, where $\mbox{\boldmath$\varphi$}' = (\varphi_1,\ldots,\varphi_N)^{\rm T}$. 
\end{lemma}

\noindent
{\bf Proof}. \ 
We note that in view of \eqref{estimate II:Q} we have 
$\|Q\psi\|_k \lesssim \|\nabla b\|_{W^{|k| \vee |k+1|,\infty}} \|\psi\|_{k+1}$. 
It follows from Lemma \ref{estimate II:formula 2} that 
\begin{equation}\label{estimate II:recursion}
\begin{cases}
 \|\varphi_1\|_k \leq \delta^2 C( \|\eta\|_{|k| \vee t_0},\|b\|_{W^{|k|+1 \vee |k+1|+1,\infty}} )
  ( \|\nabla \varphi_0\|_{k+1} + \|\mbox{\boldmath$\varphi$}'\|_{k+2} ), \\
 \|\varphi_j\|_k \leq \delta^2 C( \|\eta\|_{|k| \vee t_0},\|b\|_{W^{|k|+1 \vee |k+1|+1,\infty}} )
  ( \|\nabla \varphi_{j-2}\|_{k+1} + \|(\varphi_{j-1},\ldots,\varphi_N)\|_{k+2} )
\end{cases}
\end{equation}
for $j=2,\ldots,N$, so that 
\[
\|(\varphi_{2j-1},\varphi_{2j},\ldots,\varphi_N)\|_k
\leq \delta^2 C( \|\eta\|_{|k| \vee t_0},\|b\|_{W^{|k|+1 \vee |k+1|+1,\infty}} )
 \|(\varphi_{2j-3},\varphi_{2j-2},\ldots,\varphi_N)\|_k
\]
for $j=2,\ldots,N$. 
Using this inductively, we obtain 
\begin{align*}
& \|(\varphi_{2j-1},\varphi_{2j},\ldots,\varphi_N)\|_k \\
&\leq \delta^{2(j-1)} C( \|\eta\|_{|k| \vee |k+2(j-2)| \vee t_0},\|b\|_{W^{|k|+1 \vee |k+2j-3|+1,\infty}}) 
 \|(\varphi_1,\varphi_2,\ldots,\varphi_N)\|_{k+2(j-1)}
\end{align*}
for $j=2,\ldots,N$. 
Applying \eqref{estimate II:recursion} with $k$ replaced by $k+2(j-1)$ to the last term in the above 
inequality, we obtain the desired estimate. 
\quad$\Box$

\bigskip
Now, we can show the second half of Theorem \ref{result:theorem 1} in the case (H2). 
We apply Lemma \ref{estimate II:estimate 2} with 
$\mbox{\boldmath$\varphi$} = \mbox{\boldmath$\phi$}^\delta$ and $k = m-2j+1$. 
Then, $|k| \vee |k+2(j-1)| \vee t_0 = m-1$ and 
$|k|+1 \vee |k+2j-1|+1 =m+1$ hold if and only if the integer $j$ satisfies $1\leq j\leq m$. 
We remind our hypothesis $m \geq [(N+1)/2]$. 
Therefore, in the case of even $N$, we obtain 
\[
\|(\phi_{2j-1}^\delta,\phi_{2j}^\delta)\|_{m-2j-1} 
 \leq \delta^{2j} C_m (\|\nabla \phi_0^\delta\|_m + \|\mbox{\boldmath$\phi$}^{\delta \prime}\|_{m+1})
 \quad\mbox{for}\quad j=1,\ldots,N/2,
\]
where $C_m = C(\|\eta\|_{m-1},\|b\|_{W^{m+1,\infty}})$. 
In the case of odd $N$, we have 
\[
\begin{cases}
 \|\phi_{2j-1}^\delta\|_{m-2j-1} 
 \leq \delta^{2j} C_m (\|\nabla \phi_0^\delta\|_m + \|\mbox{\boldmath$\phi$}^{\delta \prime}\|_{m+1})
  & \mbox{for}\quad j=1,\ldots,(N+1)/2, \\
 \|\phi_{2j}^\delta\|_{m-2j-1} 
 \leq \delta^{2j} C_m (\|\nabla \phi_0^\delta\|_m + \|\mbox{\boldmath$\phi$}^{\delta \prime}\|_{m+1})
  & \mbox{for}\quad j=1,\ldots,(N+1)/2-1.
\end{cases}
\]
These estimates together with \eqref{result:uniform estimate 1} yield the desired estimate 
\eqref{result:uniform estimate 2} in the case (H2).

The proof of Theorem \ref{result:theorem 1} is complete.

\section{Consistency of the Isobe--Kakinuma model I}
\label{section:consistency I}
\setcounter{equation}{0}
\setcounter{theorem}{0}
In this and the following two sections, we will prove Theorem \ref{result:theorem 2}. 
Suppose that $(\eta,\mbox{\boldmath$\phi$}^\delta)$ is a solution of the Isobe--Kakinuma model 
\eqref{result:IK model} and define $\phi$ by \eqref{result:phi}, which is an approximation of the 
trace of the velocity potential on the water surface. 
We will show that $(\eta,\phi)$ satisfies the water wave equations in 
Zakharov--Craig--Sulem formulation \eqref{intro:WW} with an error of order $O(\delta^{4N+2})$ 
in the case (H1) and of order $O(\delta^{4[N/2]+2})$ in the case (H2). 
Here, we remind that the water wave equations in terms of the surface elevation $\eta$ and the 
velocity potential $\Phi$ have the form 
\begin{equation}\label{cst 1:Laplace}
\Delta \Phi + \delta^{-2} \partial_z^2 \Phi = 0 
 \quad\mbox{in}\quad \Omega(t),
\end{equation}
\begin{equation}\label{cst 1:BCf}
 \begin{cases}
  \displaystyle
  \partial_t \Phi + \frac12 \Bigl( |\nabla \Phi|^2 + \delta^{-2}(\partial_z \Phi)^2 \Bigr)
   + \eta = 0
   & \mbox{on} \quad \Gamma(t), \\
  \displaystyle
  \partial_t \eta + \nabla \eta \cdot \nabla \Phi - \delta^{-2} \partial_z \Phi = 0
   & \mbox{on} \quad \Gamma(t), 
 \end{cases}
\end{equation}
\begin{equation}\label{cst 1:BCb}
 \delta^{-2} \partial_z \Phi - \nabla \eta \cdot \nabla \Phi = 0
  \quad\mbox{on}\quad \Sigma,
\end{equation}
where $\Omega(t)$, $\Gamma(t)$, and $\Sigma$ denote the water region, the water surface, and the bottom, 
respectively. 
Our strategy to show the desired consistency is to use an approximate velocity potential which 
satisfy \eqref{cst 1:Laplace}--\eqref{cst 1:BCb} approximately.

We define an approximate velocity potential $\Phi^{\mbox{\rm\tiny app}}$ in the water region 
by \eqref{intro:nd approximation}. 
Then, we see that the second equation in \eqref{result:IK model} is equivalent to 
\begin{equation}\label{cst 1:Bernoulli}
\partial_t \Phi^{\mbox{\rm\tiny app}} + \frac12 \Bigl( 
 |\nabla \Phi^{\mbox{\rm\tiny app}}|^2 
 + \delta^{-2}(\partial_z \Phi^{\mbox{\rm\tiny app}})^2 \Bigr) + \eta = 0
 \quad\mbox{on}\quad z=\eta(x,t),
\end{equation}
which is exactly the first equation in \eqref{cst 1:BCf}, that is, 
Bernoulli's law restricted on the water surface. 
However, $\Phi^{\mbox{\rm\tiny app}}$ satisfies the other equations approximately with an error of 
order $O(\delta^{2N})$ in the case (H1) and of order $O(\delta^{2[N/2]})$ in the case (H2). 
These orders of the error are not sufficient to show the desired result, so that we have to modify 
$\Phi^{\mbox{\rm\tiny app}}$ appropriately.

In the following arguments, the time $t$ is arbitrarily fixed so that we omit it in the notation. 
In \eqref{t-derivative:scrLi}--\eqref{t-derivative:scrL} we defined operators 
$\mathscr{L}\mbox{\boldmath$\varphi$} = (\mathscr{L}_0\mbox{\boldmath$\varphi$},\ldots,
\mathscr{L}_N\mbox{\boldmath$\varphi$})^{\rm T}$ ,which act on $(N+1)$ vector-valued functions 
$\mbox{\boldmath$\varphi$} = (\varphi_0,\ldots,\varphi_N)^{\rm T}$. 
We denote these operators by $\mathscr{L}^{(N)}$ and $\mathscr{L}_i^{(N)}$ for $i=0,1,\ldots,N$. 
We assume that $\eta$, $\phi$, $\mbox{\boldmath$\phi$}^\delta = (\phi_0^\delta,\ldots,\phi_N^\delta)$, 
and $b$ are given so that 
\begin{equation}\label{cst 1:assumption 1}
 \mathscr{L}^{(N)}_0 \mbox{\boldmath$\phi$}^\delta = \phi, \qquad
 \mathscr{L}^{(N)}_i \mbox{\boldmath$\phi$}^\delta = 0
  \quad\mbox{for}\quad i=1,\ldots,N,
\end{equation}
and that 
\begin{equation}\label{cst 1:assumption 2}
 \begin{cases}
  \|\eta\|_m + \|\nabla \phi\|_{m-1} + \|b\|_{W^{m+1,\infty}} \leq M, \\
  H(x) = 1 + \eta(x) -b(x) \geq c_0  \quad\mbox{for}\quad x\in\mathbf{R}^n,
 \end{cases}
\end{equation}
where $m$ is an integer satisfying $m \geq n/2+1$. 
Now, we define $\widetilde{\mbox{\boldmath$\phi$}}^\delta
 = (\widetilde{\phi}_0^\delta,\widetilde{\phi}_1^\delta,\ldots,\widetilde{\phi}_{2N+2}^\delta)^{\rm T}$ 
by 
\begin{equation}\label{cst 1:modify phi}
 \mathscr{L}_0^{(2N+2)} \widetilde{\mbox{\boldmath$\phi$}}^\delta = \phi, \qquad
 \mathscr{L}_i^{(2N+2)} \widetilde{\mbox{\boldmath$\phi$}}^\delta = 0 
  \quad\mbox{for}\quad i=1,\ldots,2N+2,
\end{equation}
and then a modified approximate velocity potential $\widetilde{\Phi}^{\mbox{\rm\tiny app}}$ by 
\begin{equation}\label{sct 1:modify potential}
\widetilde{\Phi}^{\mbox{\rm\tiny app}}(x,z,t)
 = \sum_{i=0}^{2N+2} (z + 1 - b(x))^{p_i} \widetilde{\phi}_i^\delta(x,t).
\end{equation}
We will show that $\eta$ and $\widetilde{\Phi}^{\mbox{\rm\tiny app}}$ satisfy the water wave 
equations \eqref{cst 1:Laplace}--\eqref{cst 1:BCb} with an error of desirable order.

To compare $\widetilde{\phi}_j^\delta$ with $\phi_j^\delta$ for $j=0,1,\ldots,N$, 
we introduce a new function $\mbox{\boldmath$\varphi$}^\delta$ by 
\begin{equation}\label{cst 1:varphi}
\mbox{\boldmath$\varphi$}^\delta = (\varphi_0^\delta,\varphi_1^\delta,\ldots,\varphi_N^\delta)^{\rm T},
 \qquad \varphi_j^\delta = \phi_j^\delta - \widetilde{\phi}_j^\delta
 \quad\mbox{for}\quad j=0,1,\ldots,N.
\end{equation}
Then, we see that $\mbox{\boldmath$\varphi$}^\delta$ satisfies 
\begin{equation}\label{cst 1:eq varphi pre}
\mathscr{L}^{(N)} \mbox{\boldmath$\varphi$}^\delta = \mbox{\boldmath$R$} = (R_0,R_1,\ldots,R_N)^{\rm T},
\end{equation}
where 
\begin{equation}\label{cst 1:Ri}
R_0 = \sum_{j=N+1}^{2N+2} H^{p_j} \widetilde{\phi}_j^\delta, \qquad
R_i = \sum_{j=N+1}^{2N+2} (L_{ij} - H^{p_i}L_{0j}) \widetilde{\phi}_j^\delta
 \quad\mbox{for}\quad i=1,2,\ldots,N.
\end{equation}
We decompose $R_i = R_{1,i} + \delta^{-2} R_{2,i}$, where 
\begin{align}
R_{1,i} \label{cst 1:R1i}
&= \sum_{j=N+1}^{2N+2} \biggl\{
 - \biggl( \frac{1}{p_i+p_j+1} - \frac{1}{p_j+1} \biggr) H^{p_i+p_j+1} \Delta\varphi_j \\
&\qquad\qquad
 + \biggl( \frac{p_j}{p_i+p_j} - \frac{p_j}{p_j} \biggr) H^{p_i+p_j} \nabla \cdot (\varphi_j\nabla b)
 \nonumber \\
&\qquad\qquad
 - \frac{p_i}{p_i+p_j} H^{p_i+p_j} \nabla b \cdot \nabla\varphi_j
   + \frac{p_ip_j}{p_i+p_j-1} H^{p_i+p_j-1} |\nabla b|^2 \varphi_j \biggr\}, 
   \nonumber \\
R_{2,i} \label{cst 1:R2i}
&= \sum_{j=N+1}^{2N+2} \frac{p_ip_j}{p_i+p_j-1} H^{p_i+p_j-1} \widetilde{\phi}_j^\delta,
\end{align}
for $i=1,2,\ldots,N$. 
These decompositions lead a decomposition 
$\mbox{\boldmath$R$} = \mbox{\boldmath$R$}_1 + \delta^{-2} \mbox{\boldmath$R$}_2$, where 
$\mbox{\boldmath$R$}_1 = (R_0,R_{1,1},\ldots,R_{1,N})^{\rm T}$ and 
$\mbox{\boldmath$R$}_2 = (0,R_{2,1},\ldots,R_{2,N})^{\rm T}$. 
Then, we have 
\begin{equation}\label{cst 1:eq varphi}
\mathscr{L}^{(N)} \mbox{\boldmath$\varphi$}^\delta
 = \mbox{\boldmath$R$}_1 + \delta^{-2} \mbox{\boldmath$R$}_2.
\end{equation}
By using equations \eqref{cst 1:modify phi} and \eqref{cst 1:eq varphi}, we will evaluate 
$\widetilde{\mbox{\boldmath$\phi$}}^\delta$ and $\mbox{\boldmath$\varphi$}^\delta$. 
We also note that the difference between the two approximate velocity potentials 
$\widetilde{\Phi}^{\mbox{\rm\tiny app}}$ and $\Phi^{\mbox{\rm\tiny app}}$ is represented as 

\begin{equation}\label{cst 1:difference}
\widetilde{\Phi}^{\mbox{\rm\tiny app}} - \Phi^{\mbox{\rm\tiny app}}
= \sum_{j=0}^N (z+1-b)^{p_j} \varphi_j^\delta
 + \sum_{j=N+1}^{2N+2} (z+1-b)^{p_j} \widetilde{\phi}_j^\delta.
\end{equation}

\subsection{The case $p_i=2i$ with the flat bottom}
\begin{lemma}\label{cst 1:lemma 1}
Choose $p_i = 2i$ $(i = 0,1,\ldots,N)$ and suppose that $b=0$ and that $(\eta,\phi)$ satisfy 
\eqref{cst 1:assumption 2}. 
For any $j = 1,2,\ldots,2N+2$, if an integer $k$ satisfies $|k| \leq m$ and $|k+2j-1| \leq m-1$, 
then we have 
\[
\| (\widetilde{\phi}_j^\delta,\widetilde{\phi}_{j+1}^\delta,\ldots,\widetilde{\phi}_{2N+2}^\delta) \|_k
 \leq C\delta^{2j},
\]
where $C=C(M,c_0,m,j,k,N)$ is a positive constant independent of $\delta\in(0,1]$. 
\end{lemma}

\noindent
{\bf Proof}. \ 
By Lemma \ref{t-derivative:lemma 3}, particularly, the first estimate in 
\eqref{t-derivative:elliptic estimate 2}, we have 
$\|\nabla \widetilde{\phi}_0^\delta\|_k + \|\widetilde{\phi}^{\delta \prime}\|_{k+1}
 \lesssim \|\nabla \phi\|_k \lesssim 1$ if $|k| \leq m-1$, so that 
$\|\nabla \widetilde{\mbox{\boldmath$\phi$}}^\delta\|_{k+2j-1} \lesssim 1$ if 
$|k+2j-1| \leq m-1$. 
On the other hand, it follows from Lemma \ref{estimate II:estimate 1} that 
$\| (\widetilde{\phi}_j^\delta,\widetilde{\phi}_{j+1}^\delta,\ldots,\widetilde{\phi}_N^\delta) \|_k
 \lesssim \delta^{2j} \|\nabla \widetilde{\mbox{\boldmath$\phi$}}^\delta\|_{k+2j-1}$ if 
$|k| \vee |k+2(j-1)| \leq m$. 
In view of $|k+2(j-1)| \leq |k+2j-1|+1$, these two estimates give the desired one. 
\quad$\Box$

\begin{lemma}\label{cst 1:lemma 2}
Choose $p_i = 2i$ $(i = 0,1,\ldots,N)$ and suppose that $b=0$ and that $(\eta,\phi)$ satisfy 
\eqref{cst 1:assumption 2}. 
For any $j = 0,1,\ldots,N+1$, if an integer $k$ satisfies $|k-1| \vee |k| \vee |k+2j-1| \leq m-1$, 
then we have 
\[
\|\mbox{\boldmath$\varphi$}^\delta\|_k
 + \|(\widetilde{\phi}_{N+1}^\delta,\ldots,\widetilde{\phi}_{2N+2}^\delta)\|_k \leq C \delta^{2j},
\]
where $C=C(M,c_0,m,j,k,N)$ is a positive constant independent of $\delta\in(0,1]$. 
\end{lemma}

\noindent
{\bf Proof}. \ 
It follows from Lemma \ref{t-derivative:lemma 3}, particularly, the second estimate in 
\eqref{t-derivative:elliptic estimate 2} with $k$ replaced by $k-1$ that 
$\|(\mathscr{L}^{(N)})^{-1} \mbox{\boldmath$F$}\|_k
 \lesssim \|F_0\|_k + \|\mbox{\boldmath$F$}'\|_{k-2}$ if $|k-1| \leq m-1$. 
Moreover, if $F_0=0$, then we can apply \eqref{t-derivative:elliptic estimate 3} and obtain 
$\|(\mathscr{L}^{(N)})^{-1} \mbox{\boldmath$F$}\|_k
 \lesssim \delta^2 \|\mbox{\boldmath$F$}'\|_k$ if $|k| \leq m-1$. 
Therefore, in view of \eqref{cst 1:eq varphi} we obtain 
\begin{align*}
\| \mbox{\boldmath$\varphi$}^\delta \|_k
&\leq \| (\mathscr{L}^{(N)})^{-1} \mbox{\boldmath$R$}_1\|_k + \delta^{-2} 
      \| (\mathscr{L}^{(N)})^{-1} \mbox{\boldmath$R$}_2\|_k \\
&\lesssim \|R_0\|_k + \|\mbox{\boldmath$R$}_1'\|_{k-2} + \|\mbox{\boldmath$R$}_2'\|_k
 \quad\mbox{if}\quad |k-1| \vee |k| \leq m-1.
\end{align*}
Here, by the explicit form \eqref{cst 1:Ri}--\eqref{cst 1:R2i} of $R_0$, $\mbox{\boldmath$R$}_1'$, 
and $\mbox{\boldmath$R$}_2'$, we see that 
\[
\begin{cases}
 \|R_0\|_k + \|\mbox{\boldmath$R$}_2'\|_k
  \leq C(\|\eta\|_{|k| \vee t_0})
   \|(\widetilde{\phi}_{N+1}^\delta,\ldots,\widetilde{\phi}_{2N+2}^\delta)\|_k, \\
 \|\mbox{\boldmath$R$}_1'\|_{k-2}
  \leq C(\|\eta\|_{|k-2| \vee |k-1| \vee t_0})
   \|(\widetilde{\phi}_{N+1}^\delta,\ldots,\widetilde{\phi}_{2N+2}^\delta)\|_k,
\end{cases}
\]
so that 
\[
\|R_0\|_k + \|\mbox{\boldmath$R$}_1'\|_{k-2} + \|\mbox{\boldmath$R$}_2'\|_k
 \lesssim \|(\widetilde{\phi}_{N+1}^\delta,\ldots,\widetilde{\phi}_{2N+2}^\delta)\|_k
 \quad\mbox{if}\quad |k-2| \vee |k| \leq m.
\]
On the other hand, if $0\leq j\leq N+1$, then by Lemma \ref{cst 1:lemma 1} we have 
\[
\|(\widetilde{\phi}_{N+1}^\delta,\ldots,\widetilde{\phi}_{2N+2}^\delta)\|_k
 \leq \|(\widetilde{\phi}_j^\delta,\ldots,\widetilde{\phi}_{2N+2}^\delta)\|_k
 \lesssim \delta^{2j}
\]
if $|k| \leq m$ and $|k+2j-1| \leq m-1$. 
These three estimates yield $\|\mbox{\boldmath$\varphi$}^\delta\|_k \leq C \delta^{2j}$ 
if $|k-1| \vee |k| \leq m-1$, $|k-2| \vee |k| \leq m$, $|k| \leq m$, and $|k+2j-1| \leq m-1$. 
Since these last conditions on $k$ are equivalent to $|k-1| \vee |k| \vee |k+2j-1| \leq m-1$, 
we obtain the desired result. 
\quad$\Box$

\begin{remark}\label{cst 1:remark 1}
{\rm 
Lemmas \ref{cst 1:lemma 1} and \ref{cst 1:lemma 2} imply that 
$(\widetilde{\phi}_{N+1}^\delta,\ldots,\widetilde{\phi}_{2N+2}^\delta)$ and 
$(\varphi_0^\delta,\ldots,\varphi_N^\delta)$ are both of order $O(\delta^{2N+2})$ 
if $m$ is sufficiently large. 
In view of \eqref{cst 1:difference}, the difference between the two approximate velocity potentials 
$\widetilde{\Phi}^{\mbox{\rm\tiny app}}$ and $\Phi^{\mbox{\rm\tiny app}}$ is of order $O(\delta^{2N+2})$. 
}
\end{remark}

We remind that $\widetilde{\Phi}^{\mbox{\rm\tiny app}}$ was defined by \eqref{sct 1:modify potential}. 
In the case (H1), by direct calculation and Lemma \ref{estimate II:formula 1}, we see that 
\begin{equation}\label{cst 1:BVP1}
\begin{cases}
 \Delta \widetilde{\Phi}^{\mbox{\rm\tiny app}}
  + \delta^{-2} \partial_z^2 \widetilde{\Phi}^{\mbox{\rm\tiny app}} = R
   & \mbox{in}\quad \Omega, \\
 \widetilde{\Phi}^{\mbox{\rm\tiny app}} = \phi
   & \mbox{on}\quad \Gamma, \\
 \delta^{-2} \partial_z \widetilde{\Phi}^{\mbox{\rm\tiny app}} = 0
   & \mbox{on}\quad \Sigma,
\end{cases}
\end{equation}
where 
\begin{equation}\label{cst 1:R1}
R(x,z) = \sum_{j=0}^{2N+2} (z+1)^{2j} r_j(x)
\end{equation}
and
\begin{align*}
r_j(x) = 
\begin{cases}
 (2j+2)(2j+1)\beta_{j+1,2N+2} H^{4N+2-2j} \Delta \widetilde{\phi}_{2N+2}^\delta
  & \mbox{for}\quad j=0,1,\ldots,2N+1, \\
 \Delta \widetilde{\phi}_{2N+2}^\delta
  & \mbox{for}\quad j=2N+2.
\end{cases}
\end{align*}
Concerning the remainder term $R$, we have the following lemma.

\begin{lemma}\label{cst 1:lemma 3}
Choose $p_i = 2i$ $(i = 0,1,\ldots,N)$ and suppose that $b=0$ and that $(\eta,\phi)$ satisfy 
\eqref{cst 1:assumption 2}. 
For any $j = 0,1,\ldots,2N+2$, if an integer $k$ satisfies $|k| \vee |k+2| \leq m$ and $|k+2j+1| \leq m-1$, 
then we have 
\[
\|(r_0,r_1,\ldots,r_{2N+2})\|_k \leq C \delta^{2j},
\]
where $C=C(M,c_0,m,j,k,N)$ is a positive constant independent of $\delta\in(0,1]$. 
\end{lemma}

\noindent
{\bf Proof}. \ 
It is easy to see that 
$\|(r_0,r_1,\ldots,r_{2N+2})\|_k \lesssim \|\widetilde{\phi}_{2N+2}\|_{k+2}$ if $|k| \leq m$. 
By Lemma \ref{cst 1:lemma 1} with $k$ replaced by $k+2$, we have 
$\|\widetilde{\phi}_{2N+2}\|_{k+2} \lesssim \delta^{2j}$ if 
$|k+2| \vee |k+2j+1|+1 \leq m$ for $j=0,1,\ldots,2N+2$. 
Combining these estimates we obtain the desired one. 
\quad$\Box$

\begin{remark}\label{cst 1:remark 2}
{\rm 
Lemma \ref{cst 1:lemma 3} implies that the remainder term $R$ is of order $O(\delta^{4N+2})$ 
if $m$ is sufficiently large, so that the approximate velocity potential 
$\widetilde{\Phi}^{\mbox{\rm\tiny app}}$ satisfies the continuity equation 
\eqref{cst 1:Laplace} with an error of order $O(\delta^{4N+2})$ while it satisfies the 
boundary condition \eqref{cst 1:BCb} on the bottom precisely in the case of the flat bottom. 
}
\end{remark}

\subsection{The case $p_i=i$ with general bottom topographies}
\begin{lemma}\label{cst 1:lemma 4}
Choose $p_i = i$ $(i = 0,1,\ldots,N)$ and suppose that $(\eta,\phi)$ and $b$ satisfy 
\eqref{cst 1:assumption 2}. 
For any $j = 1,2,\ldots,N+1$, if an integer $k$ satisfies $|k| \leq m$ and $|k+2j-1| \leq m-1$, 
then we have 
\[
\| (\widetilde{\phi}_{2j-1}^\delta,\widetilde{\phi}_{2j}^\delta,\ldots,
  \widetilde{\phi}_{2N+1}^\delta,\widetilde{\phi}_{2N+2}^\delta) \|_k
 \leq C\delta^{2j},
\]
where $C=C(M,c_0,m,j,k,N)$ is a positive constant independent of $\delta\in(0,1]$. 
\end{lemma}

\noindent
{\bf Proof}. \ 
As in the proof of Lemma \ref{cst 1:lemma 1}, by Lemma \ref{t-derivative:lemma 3} we have 
$\|\nabla \widetilde{\phi}_0^\delta\|_{k+2j-1} + \|\widetilde{\phi}^{\delta \prime}\|_{k+2j} \lesssim 1$ 
if $|k+2j-1| \leq m-1$. 
On the other hand, it follows from Lemma \ref{estimate II:estimate 2} that 
$\| (\widetilde{\phi}_{2j-1}^\delta,\widetilde{\phi}_{2j}^\delta,\ldots,\widetilde{\phi}_{2N+2}^\delta) \|_k
 \lesssim \delta^{2j} (\|\nabla \widetilde{\phi}_0^\delta\|_{k+2j-1}
  + \|\widetilde{\phi}^{\delta \prime}\|_{k+2j})$ if 
$|k| \vee |k+2j-1| \leq m$. 
These two estimates give the desired one. 
\quad$\Box$

\begin{lemma}\label{cst 1:lemma 5}
Choose $p_i = i$ $(i = 0,1,\ldots,N)$ and suppose that $(\eta,\phi)$ and $b$ satisfy 
\eqref{cst 1:assumption 2}. 
For any $j = 0,1,\ldots,[N/2]+1$, if an integer $k$ satisfies 
$|k-1| \vee |k| \vee |k+2j-1| \leq m-1$, then we have 
\[
\|\mbox{\boldmath$\varphi$}^\delta\|_k 
 +  \|(\widetilde{\phi}_{N+1}^\delta,\ldots,\widetilde{\phi}_{2N+2}^\delta)\|_k \leq C \delta^{2j},
\]
where $C=C(M,c_0,m,j,k,N)$ is a positive constant independent of $\delta\in(0,1]$. 
\end{lemma}

\noindent
{\bf Proof}. \ 
As in the proof of Lemma \ref{cst 1:lemma 2}, we have 
\[
\| \mbox{\boldmath$\varphi$}^\delta \|_k
\lesssim \|R_0\|_k + \|\mbox{\boldmath$R$}_1'\|_{k-2} + \|\mbox{\boldmath$R$}_2'\|_k
 \quad\mbox{if}\quad |k-1| \vee |k| \leq m-1.
\]
Here, by the explicit form \eqref{cst 1:Ri}--\eqref{cst 1:R2i} of $R_0$, $\mbox{\boldmath$R$}_1'$, 
and $\mbox{\boldmath$R$}_2'$, we see that 
\[
\begin{cases}
 \|R_0\|_k + \|\mbox{\boldmath$R$}_2'\|_k
  \leq C(\|\eta\|_{|k| \vee t_0},\|b\|_{W^{|k|,\infty}})
   \|(\widetilde{\phi}_{N+1}^\delta,\ldots,\widetilde{\phi}_{2N+2}^\delta)\|_k, \\
 \|\mbox{\boldmath$R$}_1'\|_{k-2}
  \leq C(\|\eta\|_{|k-2| \vee |k-1| \vee t_0},\|b\|_{W^{|k-2|+1 \vee |k-1|+1,\infty}})
   \|(\widetilde{\phi}_{N+1}^\delta,\ldots,\widetilde{\phi}_{2N+2}^\delta)\|_k,
\end{cases}
\]
so that 
\[
\|R_0\|_k + \|\mbox{\boldmath$R$}_1'\|_{k-2} + \|\mbox{\boldmath$R$}_2'\|_k
 \lesssim \|(\widetilde{\phi}_{N+1}^\delta,\ldots,\widetilde{\phi}_{2N+2}^\delta)\|_k
 \quad\mbox{if}\quad |k-2| \vee |k| \leq m.
\]

{\it (i) The case of even $N=2N_1$: }
if $0 \leq j \leq N_1+1 = [N/2]+1$, then by Lemma \ref{cst 1:lemma 1} we have 
\begin{align*}
\|(\widetilde{\phi}_{N+1}^\delta,\ldots,\widetilde{\phi}_{2N+2}^\delta)\|_k
&= \|(\widetilde{\phi}_{2(N_1+1)-1}^\delta,\ldots,\widetilde{\phi}_{2N+2}^\delta)\|_k
 \leq \|(\widetilde{\phi}_{2j-1}^\delta,\ldots,\widetilde{\phi}_{2N+2}^\delta)\|_k \\
&\lesssim \delta^{2j}
 \quad\mbox{if}\quad |k| \vee |k+2j-1|+1 \leq m.
\end{align*}

{\it (ii) The case of odd $N=2N_1-1$: }
if $0 \leq j \leq N_1 = [N/2]+1$, then by Lemma \ref{cst 1:lemma 1} we have 
\begin{align*}
\|(\widetilde{\phi}_{N+1}^\delta,\ldots,\widetilde{\phi}_{2N+2}^\delta)\|_k
&= \|(\widetilde{\phi}_{2N_1}^\delta,\ldots,\widetilde{\phi}_{2N+2}^\delta)\|_k
 \leq \|(\widetilde{\phi}_{2j-1}^\delta,\ldots,\widetilde{\phi}_{2N+2}^\delta)\|_k \\
&\lesssim \delta^{2j}
 \quad\mbox{if}\quad |k| \vee |k+2j-1|+1 \leq m.
\end{align*}
Combining the above estimates, we obtain the desired result. 
\quad$\Box$

\begin{remark}\label{cst 1:remark 3}
{\rm 
Lemma \ref{cst 1:lemma 5} imply that 
$(\widetilde{\phi}_{N+1}^\delta,\ldots,\widetilde{\phi}_{2N+2}^\delta)$ and 
$(\varphi_0^\delta,\ldots,\varphi_N^\delta)$ are both of order 
$O(\delta^{2[N/2]+2})$ if $m$ is sufficiently large. 
In view of \eqref{cst 1:difference}, the difference between the two approximate velocity potentials 
$\widetilde{\Phi}^{\mbox{\rm\tiny app}}$ and $\Phi^{\mbox{\rm\tiny app}}$ is of order 
$O(\delta^{2[N/2]+2})$. 
}
\end{remark}

In the case (H2), by Lemma \ref{estimate II:formula 2} 
we see that the approximate velocity potential $\widetilde{\Phi}^{\mbox{\rm\tiny app}}$ defined by 
\eqref{sct 1:modify potential} satisfies 
\begin{equation}\label{cst 1:BVP2}
\begin{cases}
 \Delta \widetilde{\Phi}^{\mbox{\rm\tiny app}}
  + \delta^{-2} \partial_z^2 \widetilde{\Phi}^{\mbox{\rm\tiny app}} = R
   & \mbox{in}\quad \Omega, \\
 \widetilde{\Phi}^{\mbox{\rm\tiny app}} = \phi
   & \mbox{on}\quad \Gamma, \\
 \delta^{-2} \partial_z \widetilde{\Phi}^{\mbox{\rm\tiny app}}
  - \nabla b \cdot \nabla \widetilde{\Phi}^{\mbox{\rm\tiny app}} = r_B
   & \mbox{on}\quad \Sigma,
\end{cases}
\end{equation}
where 
\begin{equation}\label{cst 1:R2}
R(x,z) = \sum_{j=0}^{2N+2} (z+1-b(x))^j r_j(x)
\end{equation}
and
\[
r_j(x) = 
\begin{cases}
 (j+2)(j+1)\{ \gamma_{j+2,2N+1} H^{2N+1-j} ( \Delta \widetilde{\phi}_{2N+1}^\delta
   -(2N+2)Q(b) \widetilde{\phi}_{2N+2}^\delta ) \\
 \phantom{ (j+2)(j+1)\{ }
  + \gamma_{j+2,2N+2} H^{2N+2-j} \Delta \widetilde{\phi}_{2N+2}^\delta \}
  \qquad\mbox{for}\quad j=0,1,\ldots,2N, \\
 \Delta \widetilde{\phi}_{2N+1}^\delta - (2N+2) Q(b) \widetilde{\phi}_{2N+2}^\delta
  \makebox[20ex]{}\mbox{for}\quad j=2N+1, \\
 \Delta \widetilde{\phi}_{2N+2}^\delta
  \makebox[43ex]{}\mbox{for}\quad j=2N+2,
\end{cases}
\]
\begin{equation}\label{cst 1:rB}
r_B(x) = 
 \gamma_{1,2N+1} H^{2N+2} (\Delta \widetilde{\phi}_{2N+1}^\delta
   -(2N+2)Q(b) \widetilde{\phi}_{2N+2}^\delta )
  + \gamma_{1,2N+2} H^{2N+3} \Delta \widetilde{\phi}_{2N+2}^\delta.
\end{equation}
Concerning the remainder term $R$ and $r_B$, we have the following lemma.

\begin{lemma}\label{cst 1:lemma 6}
Choose $p_i = i$ $(i = 0,1,\ldots,N)$ and suppose that $(\eta,\phi)$ and $b$ satisfy 
\eqref{cst 1:assumption 2}. 
For any $j = 0,1,\ldots,N+1$, if an integer $k$ satisfies $|k| \vee |k+2| \leq m$ and $|k+2j+1| \leq m-1$, 
then we have 
\[
\|(r_0,r_1,\ldots,r_{2N+2})\|_k +\|r_B\|_k \leq C \delta^{2j},
\]
where $C=C(M,c_0,m,j,k,N)$ is a positive constant independent of $\delta\in(0,1]$. 
\end{lemma}

\noindent
{\bf Proof}. \ 
In view of $\|Q(b) \psi\|_k \leq C(\|\nabla b\|_{W^{|k| \vee |k+1|,\infty}}) \|\psi\|_k$, we see that 
$\|(r_0,r_1,\ldots,r_{2N+2})\|_k +\|r_B\|_k
 \lesssim \|(\widetilde{\phi}_{2N+1}^\delta,\widetilde{\phi}_{2N+2}^\delta)\|_{k+2}$ if 
$|k| \vee |k+1| \leq m$. 
By Lemma \ref{cst 1:lemma 4} with $k$ replaced by $k+2$, we have 
$\|(\widetilde{\phi}_{2N+1}^\delta,\widetilde{\phi}_{2N+2}^\delta)\|_{k+2} \lesssim \delta^{2j}$ if 
$|k+2| \vee |k+2j+1|+1 \leq m$ for $j=0,1,\ldots,N+1$. 
Combining these estimates we obtain the desired one. 
\quad$\Box$

\begin{remark}\label{cst 1:remark 4}
{\rm 
Lemma \ref{cst 1:lemma 6} implies that the remainder term $R$ and $r_B$ are of order $O(\delta^{2N+2})$ 
if $m$ is sufficiently large, so that the approximate velocity potential 
$\widetilde{\Phi}^{\mbox{\rm\tiny app}}$ satisfies the continuity equation \eqref{cst 1:Laplace} 
and the boundary condition \eqref{cst 1:BCb} on the bottom with an error of order $O(\delta^{2N+2})$. 
We also note that $4[N/2]+2 \leq 2N+2$. 
}
\end{remark}

\section{Consistency of the Isobe--Kakinuma model II}
\label{section:consistency II}
\setcounter{equation}{0}
\setcounter{theorem}{0}
We proceed to show that $\eta$ and $\widetilde{\Phi}^{\mbox{\rm\tiny app}}$ satisfy approximately 
the second condition in \eqref{cst 1:BCf}, that is, the kinematic boundary condition on the water 
surface. 
Since the solution $(\eta,\mbox{\boldmath$\phi$}^\delta)$ to the Isobe--Kakinuma model satisfies 
\begin{equation}\label{cst 2:IK0}
\partial_t \eta - \sum_{j=0}^N L_{0j}\phi_j^\delta = 0,
\end{equation}
we need to compare 
$(\delta^{-2} \partial_z \widetilde{\Phi}^{\mbox{\rm\tiny app}}
 - \nabla\eta \cdot \nabla\widetilde{\Phi}^{\mbox{\rm\tiny app}})|_{z=\eta}$ 
with $\sum_{j=0}^N L_{0j}\phi_j^\delta$. 
To estimate the difference between them in Sobolev spaces, we will utilize fully the duality 
$(H^k)^* = H^{-k}$, so that we will evaluate the quantity 
\begin{equation}\label{cst 2:I}
I = ((\delta^{-2} \partial_z \widetilde{\Phi}^{\mbox{\rm\tiny app}}
  - \nabla\eta \cdot \nabla\widetilde{\Phi}^{\mbox{\rm\tiny app}})|_{z=\eta}
 - \sum_{j=0}^N L_{0j}\phi_j^\delta,\psi)_{L^2}
\end{equation}
for arbitrarily fixed $\psi$. 
Regarding $\psi$ as a function on the water surface, we extend it into the water region by 
\begin{equation}\label{cst 2:Psi}
\Psi(x,z) = \sum_{j=0}^{2N+2} (z+1-b(x))^{p_j}\psi_j^\delta(x),
\end{equation}
where $\mbox{\boldmath$\psi$}^\delta = (\psi_0^\delta,\ldots,\psi_{2N+2}^\delta)^{\rm T}$ is 
defined by 
\begin{equation}\label{cst 2:psi}
 \mathscr{L}_0^{(2N+2)} \widetilde{\mbox{\boldmath$\psi$}}^\delta = \psi, \qquad
 \mathscr{L}_i^{(2N+2)} \widetilde{\mbox{\boldmath$\psi$}}^\delta = 0 
  \quad\mbox{for}\quad i=1,\ldots,2N+2.
\end{equation}
This construction of $\Psi$ from $\psi$ is the same as that of 
$\widetilde{\Phi}^{\mbox{\rm\tiny app}}$ from $\phi$. 
See \eqref{cst 1:modify phi}--\eqref{sct 1:modify potential}. 
By Green's formula, we have 
\begin{align*}
& ((\delta^{-2} \partial_z \widetilde{\Phi}^{\mbox{\rm\tiny app}}
 - \nabla\eta \cdot \nabla\widetilde{\Phi}^{\mbox{\rm\tiny app}})|_{z=\eta},\psi)_{L^2}
- ((\delta^{-2} \partial_z \widetilde{\Phi}^{\mbox{\rm\tiny app}}
 - \nabla\eta \cdot \nabla\widetilde{\Phi}^{\mbox{\rm\tiny app}})|_{z=-1+b},\psi_0^\delta)_{L^2} \\
&= \int_{\Omega} \{ \nabla \cdot (\Psi \nabla \widetilde{\Phi}^{\mbox{\rm\tiny app}})
 - \delta^{-2} \partial_z (\Psi \partial_z \widetilde{\Phi}^{\mbox{\rm\tiny app}}) \}\,{\rm d}X \\
&= \int_{\Omega} (\Delta \widetilde{\Phi}^{\mbox{\rm\tiny app}} 
  - \delta^{-2} \partial_z^2 \widetilde{\Phi}^{\mbox{\rm\tiny app}}) \Psi \,{\rm d}X
 + \int_{\Omega} \{ \nabla\widetilde{\Phi}^{\mbox{\rm\tiny app}} \cdot \nabla\Psi
  + \delta^{-2} (\partial_z \widetilde{\Phi}^{\mbox{\rm\tiny app}})(\partial_z \Psi) \} \,{\rm d}X.
\end{align*}
Since $\widetilde{\Phi}^{\mbox{\rm\tiny app}}$ satisfies \eqref{cst 1:BVP1} in the case (H1) 
and \eqref{cst 1:BVP2} in the case (H2), we obtain 
\begin{align}\label{cst 2:I123}
& ((\delta^{-2} \partial_z \widetilde{\Phi}^{\mbox{\rm\tiny app}}
 - \nabla\eta \cdot \nabla\widetilde{\Phi}^{\mbox{\rm\tiny app}})|_{z=\eta},\psi)_{L^2} \\
&= (r_B,\psi_0^\delta)_{L^2} + \int_{\Omega} R\Psi\,{\rm d}X
 + \int_{\Omega} \{ \nabla\widetilde{\Phi}^{\mbox{\rm\tiny app}} \cdot \nabla\Psi
  + \delta^{-2} (\partial_z \widetilde{\Phi}^{\mbox{\rm\tiny app}})(\partial_z \Psi) \} \,{\rm d}X 
  \nonumber \\
&=: I_1 + I_2 + I_3, \nonumber
\end{align}
where $r_B=0$ in the case (H1). 
In view of \eqref{cst 2:Psi} and the definition of $R$, that is, \eqref{cst 1:R1} in the case (H1) 
and \eqref{cst 1:R2} in the case (H2), we have 
\begin{equation}\label{cst 2:I2}
I_2 = \sum_{i,j=0}^{2N+2} \frac{1}{p_i+p_j+1} (H^{p_i+p_j+1}r_i,\psi_j^\delta)_{L^2}.
\end{equation}
In view of \eqref{sct 1:modify potential} and \eqref{cst 2:Psi}, 
by direct calculation we see that 
\[
I_3 = \sum_{i,j=0}^{2N+2} (L_{ij} \widetilde{\phi}_j^\delta,\psi_i^\delta)_{L^2}.
\]
Here, we remind that $\widetilde{\mbox{\boldmath$\phi$}}^\delta$ and 
$\mbox{\boldmath$\psi$}^\delta$ were defined by \eqref{cst 1:modify phi} and \eqref{cst 2:psi}, 
respectively, so that we have 
\[
\sum_{j=0}^{2N+2} H^{p_j} \widetilde{\phi}_j^\delta = \phi, \qquad
\sum_{j=0}^{2N+2} L_{ij} \widetilde{\phi}_j^\delta
 = \sum_{j=0}^{2N+2} H^{p_i}L_{0j} \widetilde{\phi}_j^\delta
 \quad\mbox{for}\quad 1\leq i\leq 2N+2,
\]
and similar relations hold for $\mbox{\boldmath$\psi$}^\delta$. 
Moreover, $\mbox{\boldmath$\phi$}^\delta$ satisfies \eqref{cst 1:assumption 1}, 
so that we have also 
\[
\sum_{j=0}^N H^{p_j} \phi_j^\delta = \phi, \qquad
\sum_{j=0}^N L_{ij} \phi_j^\delta
 = \sum_{j=0}^N H^{p_i}L_{0j} \phi_j^\delta
 \quad\mbox{for}\quad 1\leq i\leq N.
\]
Using these relations and $L_{ij}^*=L_{ji}$, we can rewrite $I_3$ as 
\begin{align*}
I_3 &= \sum_{i,j=0}^{2N+2} (\widetilde{\phi}_j^\delta, L_{ji} \psi_i^\delta)_{L^2} 
 = \sum_{i,j=0}^{2N+2} (\widetilde{\phi}_j^\delta, H^{p_j} L_{0i} \psi_i^\delta)_{L^2} 
 = \sum_{i,j=0}^{2N+2} (H^{p_j} \widetilde{\phi}_j^\delta, L_{0i} \psi_i^\delta)_{L^2} \\
&= \sum_{i=0}^{2N+2} (\phi, L_{0i} \psi_i^\delta)_{L^2}
 = \sum_{j=0}^N \sum_{i=0}^{2N+2} (H^{p_j} \phi_j^\delta, L_{0i} \psi_i^\delta)_{L^2}
 = \sum_{j=0}^N \sum_{i=0}^{2N+2} (\phi_j^\delta, L_{ji} \psi_i^\delta)_{L^2} \\
&= \sum_{i,j=0}^N (L_{ij} \phi_j^\delta, \psi_i^\delta)_{L^2}
 + \sum_{j=0}^N \sum_{i=N+1}^{2N+2} (L_{ij} \phi_j^\delta, \psi_i^\delta)_{L^2}. 
\end{align*}
Here, for the first term of the right-hand side we see that 
\begin{align*}
\sum_{i,j=0}^N (L_{ij} \phi_j^\delta, \psi_i^\delta)_{L^2}
&= \sum_{i,j=0}^N (H^{p_i}L_{0j} \phi_j^\delta, \psi_i^\delta)_{L^2}
 = \sum_{i,j=0}^N (L_{0j} \phi_j^\delta, H^{p_i}\psi_i^\delta)_{L^2} \\
&= \sum_{j=0}^N (L_{0j} \phi_j^\delta, \psi)_{L^2}
 - \sum_{j=0}^N \sum_{i=N+1}^{2N+2} (L_{0j} \phi_j^\delta, H^{p_i}\psi_i^\delta)_{L^2},
\end{align*}
so that 
\begin{align}\label{cst 2:I312}
&I_3 - \sum_{j=0}^N (L_{0j} \phi_j^\delta, \psi)_{L^2} 
 = \sum_{j=0}^N \sum_{i=N+1}^{2N+2} ((L_{ij}-H^{p_i}L_{0j}) \phi_j^\delta, \psi_i^\delta)_{L^2} \\
&= \sum_{j=0}^N \sum_{i=N+1}^{2N+2}
 ((L_{ij}-H^{p_i}L_{0j}) (\phi_j^\delta - \widetilde{\phi}_j^\delta), \psi_i^\delta)_{L^2} 
 + \sum_{j=0}^N \sum_{i=N+1}^{2N+2} ((L_{ij}-H^{p_i}L_{0j}) \widetilde{\phi}_j^\delta, \psi_i^\delta)_{L^2} 
 \nonumber \\
&= \sum_{j=0}^N \sum_{i=N+1}^{2N+2} ((L_{ij}-H^{p_i}L_{0j}) \varphi_j^\delta, \psi_i^\delta)_{L^2} 
 - \sum_{j=N+1}^{2N+2} \sum_{i=N+1}^{2N+2} ((L_{ij}-H^{p_i}L_{0j}) \widetilde{\phi}_j^\delta, \psi_i^\delta)_{L^2} 
 \nonumber \\
&=: I_{3,1} + I_{3,2}, \nonumber
\end{align}
where $\mbox{\boldmath$\varphi$}^\delta = (\varphi_0^\delta,\ldots,\varphi_N^\delta)^{\rm T}$ 
was defined by \eqref{cst 1:varphi}. 
Summarizing the above calculations, the quantity $I$ defined by \eqref{cst 2:I} is decomposed as 
\[
I = I_1 + I_2 + I_{3,1} + I_{3,2}.
\]
By using this expression, we will evaluate the quantity $I$ in the following.

\subsection{The case $p_i=2i$ with the flat bottom}
\begin{lemma}\label{cst 2:lemma 1}
Choose $p_i = 2i$ $(i = 0,1,\ldots,N)$ and suppose that $b=0$ and that $(\eta,\phi)$ satisfy 
\eqref{cst 1:assumption 2}. 
For any $j = 0,1,\ldots,2N+2$, if an integer $k$ satisfies $|k+2j| \leq m$ and $|k+1| \leq m-1$, 
then we have 
\[
\| (\psi_j^\delta,\psi_{j+1}^\delta,\ldots,\psi_{2N+2}^\delta) \|_{-(k+2j)}
 \leq C\delta^{2j} \|\psi\|_{-k},
\]
where $C=C(M,c_0,m,j,k,N)$ is a positive constant independent of $\delta\in(0,1]$. 
\end{lemma}

\noindent
{\bf Proof}. \ 
By Lemma \ref{t-derivative:lemma 3}, particularly, the second estimate in 
\eqref{t-derivative:elliptic estimate 2} with $k$ replaced by $k-1$, we have 
$\|\mbox{\boldmath$\psi$}^\delta\|_{k} \lesssim \|\psi\|_k$ if $|k-1| \leq m-1$. 
On the other hand, it follows from Lemma \ref{estimate II:estimate 1} that 
$\| (\psi_j^\delta,\psi_{j+1}^\delta,\ldots,\psi_N^\delta) \|_k
 \lesssim \delta^{2j} \|\mbox{\boldmath$\psi$}^\delta\|_{k+2j}$ if 
$|k| \vee |k+2(j-1)| \leq m$. 
These two estimates give 
$\| (\psi_j^\delta,\psi_{j+1}^\delta,\ldots,\psi_N^\delta) \|_k
 \lesssim \delta^{2j} \|\psi\|_{k+2j}$ if $|k| \vee |k+2j-1|+1 \leq m$. 
Replacing $k$ with $-(k+2j)$, we obtain the desired result. 
\quad$\Box$

\begin{lemma}\label{cst 2:lemma 2}
Choose $p_i = 2i$ $(i = 0,1,\ldots,N)$ and suppose that $b=0$ and that $(\eta,\phi)$ and 
$\mbox{\boldmath$\phi$}^\delta$ satisfy \eqref{cst 1:assumption 1}--\eqref{cst 1:assumption 2}. 
For any $l = 0,1,\ldots,2N+1$, if $m \geq l+1+\delta_{l1}$, then we have 
\[
\| (\delta^{-2} \partial_z \widetilde{\Phi}^{\mbox{\rm\tiny app}}
  - \nabla\eta \cdot \nabla\widetilde{\Phi}^{\mbox{\rm\tiny app}})|_{z=\eta}
 - \sum_{j=0}^N L_{0j}\phi_j^\delta \|_{m-2(l+1)} \leq C\delta^{2l},
\]
where $C=C(M,c_0,m,l,N)$ is a positive constant independent of $\delta\in(0,1]$ and 
$\delta_{l1}$ is the Kronecker delta. 
\end{lemma}

\noindent
{\bf Proof}. \ 
In the case $l=0$, we do not need to use the duality argument, and 
by direct evaluation and Lemma \ref{cst 1:lemma 1} we obtain the estimate of the lemma. 
Therefore, we will consider the case $1 \leq l \leq 2N+1$. 
By assumption we have $I_1=0$, so that it is sufficient to evaluate $I_2$, $I_{3,1}$, and $I_{3,2}$. 
It follows from \eqref{cst 2:I2} that 
\[
|I_2| \lesssim \|(r_0,r_1,\ldots,r_{2N+2})\|_{k} \|\mbox{\boldmath$\psi$}^\delta\|_{-k}
 \quad\mbox{if}\quad |k| \leq m.
\]
Here, by Lemmas \ref{cst 1:lemma 3} and \ref{cst 2:lemma 1} we have 
\[
\begin{cases}
 \|(r_0,r_1,\ldots,r_{2N+2})\|_{k} \lesssim \delta^{2l}
  & \mbox{if}\quad |k| \vee |k+2| \leq m, \; |k+2l+1| \leq m-1, \\
 \|\mbox{\boldmath$\psi$}^\delta\|_{-k} \lesssim \|\psi\|_{-k}
  & \mbox{if}\quad |k| \leq m, \; |k+1| \leq m-1,
\end{cases}
\]
so that $|I_2| \lesssim \delta^{2l}\|\psi\|_{-k}$ if 
$ |k| \vee |k+2| \leq m$ and $|k+1| \vee |k+2l+1| \leq m-1$. 
Since these conditions on $k$ are equivalent to $-m \leq k \leq m-2(l+1)$, 
if $m \geq l+1$, then we can take $k=m-2(l+1)$ and obtain 
\begin{equation}\label{cst 2:estimate 1}
|I_2| \lesssim \delta^{2l} \|\psi\|_{-(m-2(l+1))}.
\end{equation}
We proceed to evaluate $I_{3,1}$ and $I_{3,2}$. 
To this end, we note that in the case of the flat bottom we have 
\[
\|(L_{ij}-H^{p_i}L_{0j})\varphi\|_k
 \leq C(\|\eta\|_{|k| \vee t_0}) ( \|\varphi\|_{k+2} + \delta^{-2} \|\varphi\|_k). 
\]
We decompose $l$ into a sum of two integers $l_1$ and $l_2$ satisfying $0 \leq l_1 \leq N+1$ and 
$1\leq l_2 \leq N$. 
It follows from \eqref{cst 2:I312} that 
\[
|I_{3,1}| \lesssim ( \|\mbox{\boldmath$\varphi$}^\delta\|_{k+2l_1+2}
 + \delta^{-2} \|\mbox{\boldmath$\varphi$}^\delta\|_{k+2l_1} )
 \|(\psi_{N+1}^\delta,\ldots,\psi_{2N+2}^\delta)\|_{-(k+2l_1)}
 \quad\mbox{if}\quad |k+2l_1| \leq m.
\]
Here, by Lemma \ref{cst 1:lemma 2} with $(k,j)$ replaced by $(k+2l_1+2,l_2)$ and $(k+2l_2,l_2+1)$ 
\[
\begin{cases}
 \|\mbox{\boldmath$\varphi$}^\delta\|_{k+2l_1+2} \lesssim \delta^{2l_2}
  \quad\mbox{if}\quad |k+2l_1+1| \vee |k+2l_1+2| \vee |k+2(l_1+l_2)+1| \leq m-1, \\
 \|\mbox{\boldmath$\varphi$}^\delta\|_{k+2l_1} \lesssim \delta^{2l_2+2}
  \quad\mbox{if}\quad |k+2l_1-1| \vee |k+2l_1| \vee |k+2(l_1+l_2)+1| \leq m-1, 
\end{cases}
\]
and by Lemma \ref{cst 2:lemma 1} with $l$ replaced by $l_1$
\[
\\
 \|(\psi_{N+1}^\delta,\ldots,\psi_{2N+2}^\delta)\|_{-(k+2l_1)} \lesssim \delta^{2l_1} \|\psi\|_{-k}
  \quad\mbox{if}\quad |k+2l_1| \leq m, \; |k+1| \leq m-1,
\]
so that $|I_{3,1}| \lesssim \delta^{2l} \|\psi\|_{-k}$ if 
$|k+1| \vee |k+2l_1-1| \vee |k+2(l_1+l_2)+1| \leq m-1$. 
In the case $l \geq 2$, we can take $l_1 \geq 1$ so that these conditions on $k$ is 
equivalent to $-m \leq k \leq m-2l-2$. 
Therefore, if $m \geq l+1$, then we can take $k=m-2(l+1)$. 
In the case $l=1$, we have $l_1=0$ and $l_2=1$ so that these conditions on $k$ is equivalent to 
$-m +2 \leq k \leq m - 4$. 
Therefore, if $m \geq 3$, then we can take $k=m-4$. 
In any case, if $m \geq l+1+\delta_{l1}$, then we obtain 
\begin{equation}\label{cst 2:estimate 2}
|I_{3,1}| \lesssim \delta^{2l} \|\psi\|_{-(m-2(l+1))}.
\end{equation}
Similarly, it follows from \eqref{cst 2:I312} that 
\begin{align*}
|I_{3,2}| 
&\lesssim ( \|(\widetilde{\phi}_{N+1}^\delta,\ldots,\widetilde{\phi}_{2N+2}^\delta)\|_{k+2l_1+2}
 + \delta^{-2} \|(\widetilde{\phi}_{N+1}^\delta,\ldots,\widetilde{\phi}_{2N+2}^\delta)\|_{k+2l_1} ) \\
&\qquad
 \times \|(\psi_{N+1}^\delta,\ldots,\psi_{2N+2}^\delta)\|_{-(k+2l_1)}
 \qquad\mbox{if}\quad |k+2l_1| \leq m.
\end{align*}
Here, by Lemma \ref{cst 1:lemma 1} with $(k,j)$ replaced by $(k+2l_1+2,l_2)$ and $(k+2l_2,l_2+1)$ 
\[
\begin{cases}
\|(\widetilde{\phi}_{N+1}^\delta,\ldots,\widetilde{\phi}_{2N+2}^\delta)\|_{k+2l_1+2}
 \lesssim \delta^{2l_2} 
 & \mbox{if}\quad |k+2l_1+2| \leq m, \; |k+2(l_1+l_2)+1| \leq m-1, \\
\|(\widetilde{\phi}_{N+1}^\delta,\ldots,\widetilde{\phi}_{2N+2}^\delta)\|_{k+2l_1}
 \lesssim \delta^{2l_2+2} 
 & \mbox{if}\quad |k+2l_1| \leq m, \; |k+2(l_1+l_2)+1| \leq m-1,
\end{cases}
\]
so that $|I_{3,2}| \lesssim \delta^{2l} \|\psi\|_{-k}$ if 
$|k+2l_1| \vee |k+2l_1+2| \leq m$ and $|k+1| \vee |k+2(l_1+l_2)+1| \leq m-1$. 
Since these conditions on $k$ are equivalent to $-m \leq k \leq m-2(l+1)$, 
if $m \geq l+1$, then we can take $k=m-2(l+1)$ and obtain 
$|I_{3,2}| \lesssim \delta^{2l} \|\psi\|_{-(m-2(l+1))}$. 
This together with \eqref{cst 2:estimate 1} and \eqref{cst 2:estimate 2} yields 
$|I| \lesssim \delta^{2l} \|\psi\|_{-(m-2(l+1))}$, that is, 
\[
|((\delta^{-2} \partial_z \widetilde{\Phi}^{\mbox{\rm\tiny app}}
  - \nabla\eta \cdot \nabla\widetilde{\Phi}^{\mbox{\rm\tiny app}})|_{z=\eta}
  - \sum_{j=0}^N L_{0j}\phi_j^\delta, \psi)_{L^2}|
\lesssim \delta^{2l} \|\psi\|_{-(m-2(l+1))}
\]
for any $\psi$. 
Therefore, by the duality $(H^{m-2(l+1)})^* = H^{-(m-2(l+1))}$ we obtain the desired estimate. 
\quad$\Box$

\begin{remark}
{\rm 
Lemma \ref{cst 2:lemma 2} implies that for the solution $(\eta,\mbox{\boldmath$\phi$}^\delta)$ of 
the Isobe--Kakinuma model, $(\eta,\widetilde{\Phi}^{\mbox{\rm\tiny app}})$ satisfies the second 
condition in \eqref{cst 1:BCf} with an error of order $O(\delta^{4N+2})$ if $m$ is sufficiently large. 
}
\end{remark}

\subsection{The case $p_i=i$ with general bottom topographies}
\begin{lemma}\label{cst 2:lemma 3}
Choose $p_i = i$ $(i = 0,1,\ldots,N)$ and suppose that $(\eta,\phi)$ and $b$ satisfy 
\eqref{cst 1:assumption 2}. 
For any $j = 0,1,\ldots,N+1$, if an integer $k$ satisfies $|k+2j| \leq m$ and $|k+1| \leq m-1$, 
then we have 
\[
\| (\psi_{2j-1}^\delta,\psi_{2j}^\delta,\ldots,
  \psi_{2N+1}^\delta,\psi_{2N+2}^\delta) \|_{-(k+2j)}
 \leq C\delta^{2j} \|\psi\|_{-k},
\]
where $C=C(M,c_0,m,j,k,N)$ is a positive constant independent of $\delta\in(0,1]$ and we used a 
notational convention $\psi_{-1}^\delta=0$. 
Particularly, for any $j = 0,1,\ldots,[N/2]+1$, if an integer $k$ satisfies 
$|k+2j| \leq m$ and $|k+1| \leq m-1$, then we have 
\[
\| (\psi_{N+1}^\delta,\psi_{N+2}^\delta,\ldots,\psi_{2N+2}^\delta) \|_{-(k+2j)}
 \leq C\delta^{2j} \|\psi\|_{-k}.
\]
\end{lemma}

\noindent
{\bf Proof}. \ 
As in the proof of Lemma \ref{cst 2:lemma 1}, we have 
$\|\mbox{\boldmath$\psi$}^\delta\|_k \lesssim \|\psi\|_k$ if $|k-1| \leq m-1$. 
It follows from Lemma \ref{estimate II:estimate 2} that 
$\|(\psi_{2j-1}^\delta,\psi_{2j}^\delta,\ldots,\psi_{2N+2}^\delta) \|_k
 \lesssim \delta^{2j} \|\mbox{\boldmath$\psi$}^\delta\|_{k+2j}$ if 
$|k| \vee |k+2(j-1)| \vee |k+2j-1| \leq m$. 
These two estimates give 
$\|(\psi_{2j-1}^\delta,\psi_{2j}^\delta,\ldots,\psi_{2N+2}^\delta) \|_k
 \lesssim \delta^{2j} \|\psi\|_{k+2j}$ if $|k| \leq m$ and $|k+2j-1| \leq m-1$. 
Replacing $k$ with $-(k+2j)$ we obtain the desired result. 
The later part of the lemma comes from the former one as in the proof of 
Lemma \ref{cst 1:lemma 5}. 
\quad$\Box$

\begin{lemma}\label{cst 2:lemma 4}
Choose $p_i = i$ $(i = 0,1,\ldots,N)$ and suppose that $(\eta,\phi)$ and $b$ satisfy 
\eqref{cst 1:assumption 1}--\eqref{cst 1:assumption 2}. 
For any $l = 0,1,\ldots,2[N/2]+1$, if $m \geq l+1+\delta_{l1}$, then we have 
\[
\| (\delta^{-2} \partial_z \widetilde{\Phi}^{\mbox{\rm\tiny app}}
  - \nabla\eta \cdot \nabla\widetilde{\Phi}^{\mbox{\rm\tiny app}})|_{z=\eta}
 - \sum_{j=0}^N L_{0j}\phi_j^\delta \|_{m-2(l+1)-\delta_{l1}} \leq C\delta^{2l},
\]
where $C=C(M,c_0,m,l,N)$ is a positive constant independent of $\delta\in(0,1]$ and 
$\delta_{l1}$ is the Kronecker delta. 
\end{lemma}

\noindent
{\bf Proof}. \ 
In the case $l=0$, we do not need to use the duality argument, and 
by direct evaluation and Lemma \ref{cst 1:lemma 4} we obtain the estimate of the lemma. 
Therefore, we will consider the case $1 \leq l \leq 2[N/2]+1$. 
It follows from \eqref{cst 2:I123}--\eqref{cst 2:I2} that 
\[
|I_1|+|I_2| \lesssim (\|r_B\|_k + \|(r_0,r_1,\ldots,r_{2N+2})\|_{k}) \|\mbox{\boldmath$\psi$}^\delta\|_{-k}
 \quad\mbox{if}\quad |k| \leq m.
\]
Here, in view of $2[N/2]+1 \leq N+1$ by Lemmas \ref{cst 1:lemma 6} and \ref{cst 2:lemma 3} we have 
\[
\begin{cases}
 \|r_B\|_k + \|(r_0,r_1,\ldots,r_{2N+2})\|_{k} \lesssim \delta^{2l}
  & \mbox{if}\quad |k| \vee |k+2| \leq m, \; |k+2l+1| \leq m-1, \\
 \|\mbox{\boldmath$\psi$}^\delta\|_{-k} \lesssim \|\psi\|_{-k}
  & \mbox{if}\quad |k| \leq m, \; |k+1| \leq m-1.
\end{cases}
\]
Therefore, as in the proof of Lemma \ref{cst 2:lemma 2}, we have 
$|I_1| + |I_2| \lesssim \delta^{2l} \|\psi\|_{-(m-2(l+1))}$.

We proceed to evaluate $I_{3,1}$ and $I_{3,2}$. 
To this end, we note that 
\[
\|(L_{ij}-H^{p_i}L_{0j})\varphi\|_k
 \leq C(\|\eta\|_{|k| \vee t_0},\|b\|_{W^{|k|+1 \vee |k+1|+1,\infty}})
  ( \|\varphi\|_{k+2} + \delta^{-2} \|\varphi\|_k). 
\]
As before, we decompose $l$ into a sum of two integers $l_1$ and $l_2$ satisfying 
$1 \leq l_1 \leq [N/2]+1$ and $0\leq l_2 \leq  [N/2]$. 
It follows from \eqref{cst 2:I312} that 
\begin{align*}
|I_3| &\lesssim \{ 
  \|(\widetilde{\phi}_{N+1}^\delta,\ldots,\widetilde{\phi}_{2N+2}^\delta)\|_{k+2l_1+2}
 + \delta^{-2} \|(\widetilde{\phi}_{N+1}^\delta,\ldots,\widetilde{\phi}_{2N+2}^\delta)\|_{k+2l_1} \\
&\quad\;
 + \|\mbox{\boldmath$\varphi$}^\delta\|_{k+2l_1+2}
 + \delta^{-2} \|\mbox{\boldmath$\varphi$}^\delta\|_{k+2l_1} \}
 \|(\psi_{N+1}^\delta,\ldots,\psi_{2N+2}^\delta)\|_{-(k+2l_1)}
\end{align*}
if $|k+2l_1| \vee |k+2l_1+1| \leq m$. 
Here, by Lemma \ref{cst 1:lemma 5} with $(k,j)$ replaced by $(k+2l_1+2,l_2)$ and $(k+2l_2,l_2+1)$ 
we have 
\[
\begin{cases}
\|(\widetilde{\phi}_{N+1}^\delta,\ldots,\widetilde{\phi}_{2N+2}^\delta)\|_{k+2l_1+2}
 + \|\mbox{\boldmath$\varphi$}^\delta\|_{k+2l_1+2} \lesssim \delta^{2l_2}, \\
\|(\widetilde{\phi}_{N+1}^\delta,\ldots,\widetilde{\phi}_{2N+2}^\delta)\|_{k+2l_1}
 + \|\mbox{\boldmath$\varphi$}^\delta\|_{k+2l_1} \lesssim \delta^{2l_2+2}
\end{cases}
\]
if $|k+2l_1-1| \vee |k+2l_1+2| \vee |k+2(l_1+l_2)+1| \leq m-1$. 
By Lemma \ref{cst 2:lemma 3} we have also 
\[
\|(\psi_{N+1}^\delta,\ldots,\psi_{2N+2}^\delta)\|_{-(k+2l_1)}
 \lesssim \delta^{2l_1} \|\psi\|_{-k}
 \quad\mbox{if}\quad |k+2l_1| \leq m, \; |k+1| \leq m-1.
\]
Therefore, we obtain $|I_3| \lesssim \delta^{2l} \|\psi\|_{-k}$ if 
$|k+1| \vee |k+2l_1-1| \vee |k+2l_1+2| \vee |k+2(l_1+l_2)+1| \leq m-1$. 
In the case $l \geq 2$, we can take $l_2 \geq 1$ so that these conditions on $k$ is 
equivalent to $-m \leq k \leq m-2l-2$. 
Therefore, if $m \geq l+1$, then we can take $k=m-2(l+1)$. 
In the case $l=1$, we have $l_1=1$ and $l_2=0$ so that these conditions on $k$ is equivalent to 
$-m \leq k \leq m - 5$. 
Therefore, if $m \geq 3$, then we can take $k=m-5$. 
In any case, we obtain $|I_3| \lesssim \delta^{2l} \|\psi\|_{-(m-2(l+1)-\delta_{l1})}$ 
if $m \geq l+1+\delta_{l1}$. 

Summarizing the above estimate, we obtain $|I| \leq \delta^{2l} \|\psi\|_{-(m-2(l+1)-\delta_{l1})}$ 
if $m \geq l+1+\delta_{l1}$. 
This implies the desired estimate. 
\quad$\Box$

\begin{remark}
{\rm 
Lemma \ref{cst 2:lemma 4} implies that for the solution $(\eta,\mbox{\boldmath$\phi$}^\delta)$ of 
the Isobe--Kakinuma model, $(\eta,\widetilde{\Phi}^{\mbox{\rm\tiny app}})$ satisfies the second 
condition in \eqref{cst 1:BCf} with an error of order $O(\delta^{4[N/2]+2})$ if $m$ is 
sufficiently large. 
}
\end{remark}

\section{Consistency of the Isobe--Kakinuma model III}
\label{section:consistency III}
\setcounter{equation}{0}
\setcounter{theorem}{0}
In this section we will finish to prove Theorem \ref{result:theorem 2}, that is, 
a consistency of the Isobe--Kakinuma model \eqref{intro:ndIK model} 
with the water wave equations \eqref{intro:WW} in Zakharov--Craig--Sulem formulation. 
To this end, in view of \eqref{cst 2:IK0} we need to correlate $\sum_{j=0}^N L_{0j}\phi_j^\delta$ 
with $\Lambda(\eta,b,\delta)\phi$, where $\Lambda(\eta,b,\delta)$ is the Dirichlet-to-Neumann map 
for Laplace's equation defined by \eqref{intro:DN}--\eqref{intro:BBP}. 
We remind that the modified approximate velocity potential $\widetilde{\Phi}^{\mbox{\rm\tiny app}}$ 
satisfies the boundary value problem \eqref{cst 1:BVP1} or \eqref{cst 1:BVP2} 
and that $\phi$ was defined by \eqref{result:phi} from the solution 
$(\eta,\mbox{\boldmath$\phi$}^\delta)$ to the Isobe--Kakinuma model.

Let $\Phi$ be the unique solution to the boundary value problem \eqref{intro:BBP} and put 
\begin{equation}\label{cst 3:Phi^res}
\Phi^{\mbox{\rm\tiny res}} = \Phi - \widetilde{\Phi}^{\mbox{\rm\tiny app}}.
\end{equation}
Then, $\Phi^{\mbox{\rm\tiny res}}$ satisfies the boundary value problem 
\begin{equation}\label{cst 3:BBP}
\left\{
 \begin{array}{lll}
  \Delta\Phi^{\mbox{\rm\tiny res}} + \delta^{-2} \partial_z^2\Phi^{\mbox{\rm\tiny res}} = -R
   & \mbox{in} & \Omega, \\
  \Phi^{\mbox{\rm\tiny res}} = 0 & \mbox{on} & \Gamma, \\
  \delta^{-2} \partial_z \Phi^{\mbox{\rm\tiny res}}
   - \nabla b \cdot \nabla\Phi^{\mbox{\rm\tiny res}} = -r_B & \mbox{on} & \Sigma,
 \end{array}
\right.
\end{equation}
where $R$ and $r_B$ were defined by \eqref{cst 1:R1} and $r_B=0$ in the case (H1) 
and by \eqref{cst 1:R2}--\eqref{cst 1:rB} in the case (H2). 
Applying the identity 
\begin{align*}
\nabla \cdot \int_{-1+b}^{\eta} \nabla \Psi \,{\rm d}z
&= \int_{-1+b}^{\eta} (\Delta \Psi+ \delta^{-2} \partial_z^2 \Psi) \,{\rm d}z \\
&\quad\;
 - (\delta^{-2} \partial_z \Psi - \nabla\eta \cdot \nabla\Psi)|_{z=\eta}
 + (\delta^{-2} \partial_z \Psi - \nabla\eta \cdot \nabla\Psi)|_{z=-1+b}
\end{align*}
to $\Psi = \Phi^{\mbox{\rm\tiny res}}$ and noting \eqref{intro:DN}, we obtain 
\begin{align}\label{cst 3:error}
(\delta^{-2} \partial_z \widetilde{\Phi}^{\mbox{\rm\tiny app}}
  - \nabla\eta \cdot \nabla\widetilde{\Phi}^{\mbox{\rm\tiny app}})|_{z=\eta}
 - \Lambda(\eta,b,\delta)\phi
&= \nabla \cdot \int_{-1+b}^{\eta} \nabla \Phi^{\mbox{\rm\tiny res}} \,{\rm d}z
 + \int_{-1+b}^{\eta}R\,{\rm d}z  + r_B \\
&=: I_1 + I_2 + I_3. \nonumber
\end{align}
In view of \eqref{cst 1:R1} and \eqref{cst 1:R2}, we have 
\begin{equation}\label{cst 3:I2}
I_2 = \sum_{j=0}^{2N+2} \frac{1}{p_j+1} H^{p_j+1} r_j.
\end{equation}
Therefore, we can evaluate $I_2$ and $I_3$ directly by using the estimates obtained in 
Section \ref{section:consistency I}. 
To evaluate $I_1$ we will use an estimate for the boundary value problem \eqref{cst 3:BBP} of elliptic type. 
To this end, it is convenient to transform the problem \eqref{cst 3:BBP} in the water region $\Omega$ 
into a problem in a simple domain $\Omega_0 = \mathbf{R}^n \times (0,1)$ by using a diffeomorphism 
$\Theta(x,z) = (x,\theta(x,z)) : \Omega_0 \to \Omega$, which simply stretches the vertical direction, 
where $\theta(x,z) = \eta(x)(z+1) + (1-b(x))z$. 
We put $\widetilde{\Phi}^{\mbox{\rm\tiny res}} = \Phi^{\mbox{\rm\tiny res}} \circ \Theta$. 
Then, we have 
\begin{equation}\label{cst 3:I1}
I_1 = \nabla \cdot \int_{-1}^0 ( H \nabla\widetilde{\Phi}^{\mbox{\rm\tiny res}}
 - (\partial_z \widetilde{\Phi}^{\mbox{\rm\tiny res}}) \nabla\theta ) \,{\rm d}z,
\end{equation}
and the boundary value problem \eqref{cst 3:BBP} is transformed into 
\begin{equation}\label{cst 3:BBP2}
\left\{
 \begin{array}{lll}
 \nabla_X \cdot \mathcal{P} \nabla_X \widetilde{\Phi}^{\mbox{\rm\tiny res}}
  = - \widetilde{R}   & \mbox{in} & -1<z<0, \\
  \widetilde{\Phi}^{\mbox{\rm\tiny res}} = 0 & \mbox{on} & z=0, \\
  \mbox{\boldmath$e$}_z \cdot \mathcal{P} \nabla_X
   \widetilde{\Phi}^{\mbox{\rm\tiny res}} = -r_B & \mbox{on} & z=-1,
 \end{array}
\right.
\end{equation}
where the coefficient matrix $\mathcal{P}$ is defined by 
\[
\mathcal{P} = \det\biggl(\frac{\partial \Theta}{\partial X}\biggr)
 \biggl(\frac{\partial \Theta}{\partial X}\biggr)^{-1}
  \begin{pmatrix}
  E_n & \mbox{\boldmath$0$} \\
  \mbox{\boldmath$0$}^{\rm T} & \delta^{-2}
 \end{pmatrix}
 \biggl( \biggl(\frac{\partial \Theta}{\partial X}\biggr)^{-1} \biggr)^{\rm T},
\]
$\mbox{\boldmath$e$}_z=(0,\ldots,0,1)^{\rm T}$, and 
\begin{equation}\label{cst 3:R}
\widetilde{R} = R \circ \Theta = \sum_{j=0}^{2N+2} (z+1)^{p_j}H^{p_j}r_j.
\end{equation}
By applying the standard theory of elliptic partial differential equations to \eqref{cst 3:BBP2}, 
we obtain the following lemma. 
For details, we refer to T. Iguchi \cite{Iguchi2009, Iguchi2011} and D. Lannes \cite{Lannes2005}.

\begin{lemma}\label{cst 3:lemma 1}
Let $c_0$, $M$ be positive constant and $m$ an integer such that $m>n/2$. 
There exists a positive constant $C=C(c_0,M,m)$ such that if $\eta$ and $b$ satisfy 
\[
\begin{cases}
 \|\eta\|_m + \|b\|_{W^{m,\infty}} \leq M, \\
 c_0 < H(x) = 1 + \eta(x) - b(x)  \quad\mbox{for}\quad x\in\mathbf{R}^n,
\end{cases}
\]
and $\widetilde{\Phi}^{\mbox{\rm\tiny res}}$ is a solution to \eqref{cst 3:BBP2}, 
then for $k=0,1,\ldots,m-1$ and $\delta \in (0,1]$ we have 
\[
\|J^k \nabla \widetilde{\Phi}^{\mbox{\rm\tiny res}}\|_{L^2(\Omega_0)}
 + \delta^{-1} \|J^k \partial_z \widetilde{\Phi}^{\mbox{\rm\tiny res}}\|_{L^2(\Omega_0)}
\leq C \delta ( \|J^k \widetilde{R}\|_{L^2(\Omega_0)} + \|r_B\|_k ).
\]
\end{lemma}

We remind that we have assumed \eqref{cst 1:assumption 2}. 
Thanks of this lemma and \eqref{cst 3:R}, we see that 
\begin{equation}\label{cst 3:estimate I1}
\|I_1\|_k \lesssim  \|J^{k+2}\partial_z  \widetilde{\Phi}^{\mbox{\rm\tiny res}}\|_{L^2(\Omega_0)}
 \lesssim \delta^2 ( \|(r_0,r_1,\ldots,r_{2N+2})\|_{k+2} + \|r_B\|_{k+2})
\end{equation}
if $0\leq k+2 \leq m-1$. 
In the above calculation, we used the Poincar\'e inequality.

Now, we suppose that $(\eta,\mbox{\boldmath$\phi$}^\delta)$ is a solution to the Isobe--Kakinuma model
\eqref{result:IK model} obtained in Theorem \ref{result:theorem 1} and define $\phi$ by 
\eqref{result:phi}. 
Then, we put 
\begin{equation}\label{cst 3:r12}
 \begin{cases}
  \mathfrak{r}_1 = \partial_t\eta-\Lambda(\eta,b,\delta)\phi, \\
  \displaystyle
  \mathfrak{r}_2 = \partial_t\phi+\eta+\frac12|\nabla\phi|^2
   -\delta^2\frac{(\Lambda(\eta,b,\delta)\phi+\nabla\eta\cdot\nabla\phi)^2}{2(1+\delta^2|\nabla\eta|^2)}.
 \end{cases}
\end{equation}
We will evaluate these remainder terms $\mathfrak{r}_1$ and $\mathfrak{r}_2$ in the following. 
To this end, we put also 
\begin{equation}\label{cst 3:r345}
 \begin{cases}
  \mathfrak{r}_3 = (\delta^{-2} \partial_z \widetilde{\Phi}^{\mbox{\rm\tiny app}}
   - \nabla\eta \cdot \nabla\widetilde{\Phi}^{\mbox{\rm\tiny app}})|_{z=\eta}
   - \Lambda(\eta,b,\delta)\phi, \\
  \displaystyle
  \mathfrak{r}_4 = (\delta^{-2} \partial_z \widetilde{\Phi}^{\mbox{\rm\tiny app}}
   - \nabla\eta \cdot \nabla\widetilde{\Phi}^{\mbox{\rm\tiny app}})|_{z=\eta}
   - \sum_{j=0}^N L_{0j} \phi_j^\delta, \\
  \mathfrak{r}_5 = ((\delta^{-2} \partial_z - \nabla\eta \cdot \nabla)
   (\Phi^{\mbox{\rm\tiny app}} - \widetilde{\Phi}^{\mbox{\rm\tiny app}}))|_{z=\eta}.
 \end{cases}
\end{equation}
It follows from \eqref{cst 3:error}, \eqref{cst 3:I2}, and \eqref{cst 3:estimate I1} that 
\begin{equation}\label{cst 3:r3}
\|\mathfrak{r}_3\|_k \lesssim \delta^2 ( \|(r_0,r_1,\ldots,r_{2N+2})\|_{k+2} + \|r_B\|_{k+2})
 \quad\mbox{if}\quad 0\leq k+2 \leq m-1.
\end{equation}
We have evaluated $\mathfrak{r}_4$ in Lemmas \ref{cst 2:lemma 2} and \ref{cst 2:lemma 4}. 
By direct calculation, we see that 
\begin{align*}
\mathfrak{r}_5
&= \sum_{j=0}^N \{ p_j H^{p_j-1} (\delta^{-2} + \nabla\eta \cdot \nabla b) \varphi_j^\delta
 - H^{p_j} \nabla\eta \cdot \nabla\varphi_j^\delta \} \\
&\quad\;
 - \sum_{j=N+1}^{2N+2} \{ p_j H^{p_j-1} (\delta^{-2} + \nabla\eta \cdot \nabla b) \widetilde{\phi}_j^\delta
 - H^{p_j} \nabla\eta \cdot \nabla\widetilde{\phi}_j^\delta \}, 
\end{align*}
so that 
\begin{align}\label{cst 3:r5}
\| \mathfrak{r}_5 \|_k
&\lesssim \delta^{-2} ( \|\mbox{\boldmath$\varphi$}^\delta\|_k
 + \|(\widetilde{\phi}_{N+1}^\delta,\ldots,\widetilde{\phi}_{2N+2}^\delta)\|_k ) \\
&\quad\;
 + \|\mbox{\boldmath$\varphi$}^\delta\|_{k+1}
 + \|(\widetilde{\phi}_{N+1}^\delta,\ldots,\widetilde{\phi}_{2N+2}^\delta)\|_{k+1}
 \quad\mbox{if}\quad |k|+1 \leq m. \nonumber
\end{align}
Therefore, we can evaluate $\mathfrak{r}_3$, $\mathfrak{r}_4$, and $\mathfrak{r}_5$ by the results 
obtained in Sections \ref{section:consistency I}--\ref{section:consistency II}, 
so that it is sufficient to express $\mathfrak{r}_1$ and $\mathfrak{r}_2$ in terms of these quantities. 
It is clear that $\mathfrak{r}_1 = \mathfrak{r}_3 - \mathfrak{r}_4$. 
Differentiating the identity $\phi = \Phi^{\mbox{\rm\tiny app}}|_{z=\eta}$ with respect to $t$ and $x$, 
we have 
\begin{equation}\label{cst 3:diff}
 \begin{cases}
  \partial_t \phi = (\partial_t \Phi^{\mbox{\rm\tiny app}}
   + (\partial_z \Phi^{\mbox{\rm\tiny app}})\partial_t \eta)|_{z=\eta}, \\
  \nabla \phi = (\nabla \Phi^{\mbox{\rm\tiny app}}
   + (\partial_z \Phi^{\mbox{\rm\tiny app}})\nabla \eta)|_{z=\eta}.
 \end{cases}
\end{equation}
Plugging these into \eqref{cst 1:Bernoulli} to eliminate 
$(\partial_t \Phi^{\mbox{\rm\tiny app}})|_{z=\eta}$ and 
$(\nabla \Phi^{\mbox{\rm\tiny app}})|_{z=\eta}$ and using the first equation in \eqref{cst 3:r12} 
to eliminate $\partial_t \eta$, we obtain 
\begin{align*}
&\partial_t \phi + \eta + \frac12 |\nabla \phi|^2 \\
&= (\partial_z \Phi^{\mbox{\rm\tiny app}})|_{z=\eta}
 (\Lambda\phi + \nabla\eta \cdot \nabla\phi + \mathfrak{r}_1) 
 - \frac12 \delta^{-2} (1 + \delta^2 |\nabla\eta|^2)(\partial_z \Phi^{\mbox{\rm\tiny app}})^2|_{z=\eta} \\
&= (\partial_z \Phi^{\mbox{\rm\tiny app}})|_{z=\eta} \mathfrak{r}_1
 + \frac12 (\partial_z \Phi^{\mbox{\rm\tiny app}})|_{z=\eta}
 \{ 2(\Lambda\phi + \nabla\eta \cdot \nabla\phi)
  - \delta^{-2} (1 + \delta^2 |\nabla\eta|^2)(\partial_z \Phi^{\mbox{\rm\tiny app}})|_{z=\eta} \},
\end{align*}
where $\Lambda = \Lambda(\eta,b,\delta)$. 
On the other hand, it follows from the definition of $\mathfrak{r}_3$ and $\mathfrak{r}_5$ that 
\[
( \delta^{-2} \partial_z \Phi^{\mbox{\rm\tiny app}}
 - \nabla\eta \cdot \nabla\Phi^{\mbox{\rm\tiny app}} )|_{z=\eta}
= \Lambda\phi + \mathfrak{r}_3 + \mathfrak{r}_5, 
\]
which together with the second equation in \eqref{cst 3:diff} yields 
\[
(1 + \delta^2 |\nabla\eta|^2)(\partial_z \Phi^{\mbox{\rm\tiny app}})|_{z=\eta}
= \delta^2 (\Lambda\phi + \nabla\eta \cdot \nabla\phi + \mathfrak{r}_3 + \mathfrak{r}_5).
\]
Therefore,
\begin{align*}
&\partial_t \phi + \eta + \frac12 |\nabla \phi|^2 \\
&= (\partial_z \Phi^{\mbox{\rm\tiny app}})|_{z=\eta} \mathfrak{r}_1
 + \delta^2 \frac{(\Lambda\phi + \nabla\eta \cdot \nabla\phi + (\mathfrak{r}_3 + \mathfrak{r}_5))
  (\Lambda\phi + \nabla\eta \cdot \nabla\phi - (\mathfrak{r}_3 + \mathfrak{r}_5) )}{
 2(1 + \delta^2 |\nabla\eta|^2)} \\
&= (\partial_z \Phi^{\mbox{\rm\tiny app}})|_{z=\eta} \mathfrak{r}_1
 + \delta^2 \frac{(\Lambda\phi + \nabla\eta \cdot \nabla\phi)^2
  - (\mathfrak{r}_3 + \mathfrak{r}_5)^2}{2(1 + \delta^2 |\nabla\eta|^2)},
\end{align*}
so that 
\begin{equation}\label{cst 3:r1}
 \begin{cases}
  \mathfrak{r}_1 = \mathfrak{r}_3 - \mathfrak{r}_4, \\
  \displaystyle
  \mathfrak{r}_2 = (\partial_z \Phi^{\mbox{\rm\tiny app}})|_{z=\eta} \mathfrak{r}_1
   - \delta^2 \frac{(\mathfrak{r}_3 + \mathfrak{r}_5)^2}{2(1 + \delta^2 |\nabla\eta|^2)}.
 \end{cases}
\end{equation}
Here, in view of 
$(\partial_z \Phi^{\mbox{\rm\tiny app}})|_{z=\eta} = \sum_{i=1}^N p_i H^{p_i-1} \phi_i^\delta
 = \delta^2 w$, we have $\|(\partial_z \Phi^{\mbox{\rm\tiny app}})|_{z=\eta}\|_m \lesssim \delta$.

\subsection{The case $p_i=2i$ with the flat bottom}
It follows directly from Lemma \ref{cst 2:lemma 2} that 
$\|\mathfrak{r}_4\|_{m-4(N+1)} \lesssim \delta^{4N+2}$ if $m \geq 2(N+1)$. 
By \eqref{cst 3:r3} and Lemma \ref{cst 1:lemma 3} with $(k,j)$ replaced by $(k+2,l-1)$, 
for $l=1,2,\ldots,2N+3$ we have $\|\mathfrak{r}_3\|_k \lesssim \delta^{2l}$ if 
$0 \leq k+2 \leq m-1$, $|k+2| \vee |k+4| \leq m$, and $|k+2l+1| \leq m-1$. 
These conditions on $k$ are equivalent to $-2 \leq k\leq m-2(l+1)$. 
Particularly, we have 
\[
\begin{cases}
 \|\mathfrak{r}_3\|_{m-4(N+1)} \lesssim \delta^{4N+2} & \mbox{if}\quad m \geq 4N+2, \\
 \|\mathfrak{r}_3\|_{m-2(N+1)} \lesssim \delta^{2N} & \mbox{if}\quad m \geq 2N,
\end{cases}
\]
so that 
\[
\|\mathfrak{r}_1\|_{m-4(N+1)} \lesssim \delta^{4N+2} \quad\mbox{if}\quad m \geq 4N+2.
\]
By \eqref{cst 3:r5} and Lemma \ref{cst 1:lemma 2} with $(k,j)$ replaced by 
$(k,N+1)$ and $(k+1,N)$ we have 
$\|\mathfrak{r}_5\|_k \lesssim \delta^{2N}$ if 
$|k| \leq m-1$, $|k-1| \vee |k| \vee |k+2N+1| \leq m-1$, and 
$|k| \vee |k+1| \vee |k+2N| \leq m-1$. 
These conditions on $k$ are equivalent to $-m+2 \leq k \leq m-2(N+1)$. 
Particularly, we have 
\[
\|\mathfrak{r}_5\|_{m-2(N+1)} \lesssim \delta^{2N} \quad\mbox{if}\quad m\geq N+2.
\]
Therefore, if $m \geq 4N+2$ and $m-2(N+1) > n/2$, then 
\[
\|\mathfrak{r}_2\|_{m-4(N+1)} 
\lesssim \|\mathfrak{r}_1\|_{m-4(N+1)}
 + \delta^2 (\|\mathfrak{r}_3\|_{m-2(N+1)} + \|\mathfrak{r}_5\|_{m-2(N+1)})^2 
\lesssim \delta^{4N+2},
\]
so that we obtain the desired estimate in Theorem \ref{result:theorem 2} in the case (H1).

\subsection{The case $p_i=i$ with general bottom topographies}
It follows directly from Lemma \ref{cst 2:lemma 4} that 
$\|\mathfrak{r}_4\|_{m-4([N/2]+1)} \lesssim \delta^{4[N/2]+2}$ if $m \geq 2([N/2]+1)+\delta_{N1}$. 
By \eqref{cst 3:r3} and Lemma \ref{cst 1:lemma 6} with $(k,j)$ replaced by $(k+2,l-1)$, 
for $l=1,2,\ldots,N+2$ we have $\|\mathfrak{r}_3\|_k \lesssim \delta^{2l}$ if 
$0 \leq k+2 \leq m-1$, $|k+2| \vee |k+4| \leq m$, and $|k+2l+1| \leq m-1$. 
These conditions on $k$ are equivalent to $-2 \leq k\leq m-2(l+1)$. 
Particularly, we have 
\[
\begin{cases}
 \|\mathfrak{r}_3\|_{m-4([N/2]+1)} \lesssim \delta^{4[N/2]+2} & \mbox{if}\quad m \geq 4[N/2]+2, \\
 \|\mathfrak{r}_3\|_{m-2([N/2]+1)} \lesssim \delta^{2[N/2]} & \mbox{if}\quad m \geq 2[N/2].
\end{cases}
\]
Here, we note that the later estimate in the case $N=1$ comes from a direct evaluation. 
Therefore, 
\[
\|\mathfrak{r}_1\|_{m-4([N/2]+1)} \lesssim \delta^{4[N/2]+2}
 \quad\mbox{if}\quad m \geq 4[N/2]+2+\delta_{N1}.
\]
By \eqref{cst 3:r5} and Lemma \ref{cst 1:lemma 5} with $(k,j)$ replaced by 
$(k,[N/2]+1)$ and $(k+1,[N/2])$ we have 
$\|\mathfrak{r}_5\|_k \lesssim \delta^{2[N/2]}$ if 
$|k| \leq m-1$, $|k-1| \vee |k| \vee |k+2[N/2]+1| \leq m-1$, and 
$|k| \vee |k+1| \vee |k+2[N/2]| \leq m-1$. 
These conditions on $k$ are equivalent to $-m+2 \leq k \leq m-2([N/2]+1)$. 
Particularly, we have 
\[
\|\mathfrak{r}_5\|_{m-2([N/2]+1)} \lesssim \delta^{2[N/2]} \quad\mbox{if}\quad m\geq [N/2]+2.
\]
Therefore, if $m \geq 4[N/2]+2+\delta_{N1}$ and $m-2([N/2]+1) > n/2$, then 
\[
\|\mathfrak{r}_2\|_{m-4([N/2]+1)} 
\lesssim \|\mathfrak{r}_1\|_{m-4([N/2]+1)}
 + \delta^2 (\|\mathfrak{r}_3\|_{m-2([N/2]+1)} + \|\mathfrak{r}_5\|_{m-2([N/2]+1)})^2 
\lesssim \delta^{4[N/2]+2},
\]
so that we obtain the desired estimate in Theorem \ref{result:theorem 2} in the case (H2).

The proof of Theorem \ref{result:theorem 2} is complete.

\section{Rigorous justification of the Isobe--Kakinuma model}
\label{section:justification}
\setcounter{equation}{0}
\setcounter{theorem}{0}
In this section we will prove Theorem \ref{result:theorem 4}. 
To this end we take advantage of the stability of the water wave equations \eqref{intro:WW}, 
which is given by the following theorem. 
Although the statement is not explicitly given in T. Iguchi \cite{Iguchi2009}, 
we can prove it in exactly the same way as the proof of Theorem \ref{result:theorem 3}, 
so that we omit the proof. See also D. Lannes \cite{Lannes2013-2}.

\begin{theorem}\label{result:stability}
In addition to hypothesis of Theorem \ref{result:theorem 3} we assume that 
$(\eta^{\mbox{\rm\tiny app}},\phi^{\mbox{\rm\tiny app}})$ satisfy the equations 
\[
\left\{
 \begin{array}{l}
  \partial_t\eta^{\mbox{\rm\tiny app}}
   -\Lambda(\eta^{\mbox{\rm\tiny app}},b,\delta)\phi^{\mbox{\rm\tiny app}}
    = f_1^{\mbox{\rm\tiny err}}, \\
  \displaystyle
  \partial_t\phi^{\mbox{\rm\tiny app}}+\eta^{\mbox{\rm\tiny app}}
   +\frac12|\nabla\phi^{\mbox{\rm\tiny app}}|^2
   -\delta^2\frac{(\Lambda(\eta^{\mbox{\rm\tiny app}},b,\delta)\phi^{\mbox{\rm\tiny app}}
    +\nabla\eta^{\mbox{\rm\tiny app}}\cdot\nabla\phi^{\mbox{\rm\tiny app}})^2}{2(1
     +\delta^2|\nabla\eta^{\mbox{\rm\tiny app}}|^2)}
   = f_2^{\mbox{\rm\tiny err}},
 \end{array}
\right.
\]
on a time interval $[0,T_1]$, the initial condition \eqref{intro:IC of WW}, and the uniform bound: 
\[
\left\{
 \begin{array}{l}
  \|\eta^{\mbox{\rm\tiny app}}(t)\|_{m+3+1/2}+\|\nabla\phi^{\mbox{\rm\tiny app}}(t)\|_{m+3}
   \leq M_1, \\[0.5ex]
  1+\eta^{\mbox{\rm\tiny app}}(x,t)-b(x) \geq c_0/2
   \qquad\mbox{for}\quad x\in\mathbf{R}^n, \; 0\leq t\leq T_1.
 \end{array}
\right.
\]
Let $(\eta^{\mbox{\rm\tiny WW}},\phi^{\mbox{\rm\tiny WW}})$ be the solution obtained in 
Theorem \ref{result:theorem 3} and put $T_*=\min\{T_1,T_2\}$, where $T_2$ is that in 
Theorem \ref{result:theorem 3}. 
Then, we have 
\begin{align*}
\sup_{0\leq t\leq T_*}
\bigl(\|\eta^{\mbox{\rm\tiny WW}}(t)-\eta^{\mbox{\rm\tiny app}}(t)\|_{m+2}
&  +\|\phi^{\mbox{\rm\tiny WW}}(t)-\phi^{\mbox{\rm\tiny app}}(t)\|_{m+2}\bigr) \\
& \leq C_2\sup_{0\leq t\leq T_*}(
 \|f_1^{\mbox{\rm\tiny err}}(t)\|_{m+2}+\|\Lambda_0(\delta)^{1/2}f_2^{\mbox{\rm\tiny err}}(t)\|_{m+2}),
\end{align*}
where $\Lambda_0(\delta)=\Lambda(0,0,\delta)$ and 
$C_2$ is a positive constant independent of $\delta\in(0,\delta_*]$. 
\end{theorem}

\noindent
{\bf Proof of Theorem \ref{result:theorem 4}}. \ 
Suppose that the hypotheses in Theorem \ref{result:theorem 4} are satisfied 
for the initial data $(\eta_{(0)},\phi_{(0)})$ and the bottom topography $b$. 
We will construct the initial data 
$\mbox{\boldmath$\phi$}_{(0)}^\delta = (\phi_{0(0)}^\delta,\ldots,\phi_{N(0)}^\delta)^{\rm T}$ 
as a unique solution to 
\[
\mathscr{L}_0^{(N)} \mbox{\boldmath$\phi$}_{(0)}^\delta = \phi_{(0)}, \qquad
\mathscr{L}_i^{(N)} \mbox{\boldmath$\phi$}_{(0)}^\delta = 0
 \quad\mbox{for}\quad i=1,\ldots,N.
\]
We note that the second relation is nothing but the necessary condition \eqref{result:compatibility} 
whereas the first one corresponding the approximate relation \eqref{result:phi}. 
By Lemma \ref{t-derivative:lemma 3} with $m$ replaced by 
$m+4N+8$ in the case (H1) and by $m+4[N/2]+8$ in the case (H2), we see that 
\[
\begin{cases}
 \|\nabla \mbox{\boldmath$\phi$}_{(0)}^\delta\|_{m+4N+7}
  + \delta^{-1} \|\mbox{\boldmath$\phi$}_{(0)}^{\delta \prime}\|_{m+4N+7}
  \leq C & \mbox{in the case (H1)}, \\
 \|\nabla \mbox{\boldmath$\phi$}_{(0)}^\delta\|_{m+4[N/2]+7}
  + \delta^{-1} \|\mbox{\boldmath$\phi$}_{(0)}^{\delta \prime}\|_{m+4[N/2]+7}
  \leq C & \mbox{in the case (H2)},
\end{cases}
\]
where $C=C(c_0,M_0,m,N)$ does not depend on $\delta \in (0,1]$. 
Then, by Lemma \ref{estimate I:lemma 4} we have $\|a(\cdot,0)-1\|_{m+6} \leq C\delta$. 
Therefore, by Sobolev imbedding theorem $a(x,0) \geq 1/2$ for $x \in \mathbf{R}$ if we 
take $\delta_*$ sufficiently small, so that the initial data 
$(\eta_{(0)},\mbox{\boldmath$\phi$}_{(0)}^\delta)$ satisfy the condition \eqref{result:uniform estimate 1} 
in Theorem \ref{result:theorem 1} is satisfied and the solution 
$(\eta^{\mbox{\rm\tiny IK}},\mbox{\boldmath$\phi$}^\delta)$ 
to the Isobe--Kakinuma model exists on the time interval $[0,T_1]$ satisfying 
\[
\begin{cases}
 \|\eta^{\mbox{\rm\tiny IK}}(t)\|_{m+4N+7} + \|\nabla \mbox{\boldmath$\phi$}^\delta (t)\|_{m+4N+7}
  + \delta^{-1} \|\mbox{\boldmath$\phi$}^{\delta \prime} (t)\|_{m+4N+7}
  \leq M & \mbox{in the case (H1)}, \\
 \|\eta^{\mbox{\rm\tiny IK}}(t)\|_{m+4[N/2]+7} + \|\nabla \mbox{\boldmath$\phi$}^\delta (t)\|_{m+4[N/2]+7}
  + \delta^{-1} \|\mbox{\boldmath$\phi$}^{\delta \prime} (t)\|_{m+4[N/2]+7}
  \leq M & \mbox{in the case (H2)}, \\
 1 + \eta^{\mbox{\rm\tiny IK}}(x,t) - b(x) \geq c_0/2
  \quad\mbox{for}\quad x \in \mathbf{R}^n, \; 0 \leq t \leq T_1.
\end{cases}
\]
Now, we define $\phi^{\mbox{\rm\tiny IK}}$ by \eqref{result:phi}, that is, 
$\phi^{\mbox{\rm\tiny IK}} = \sum_{i=0}^N H^{p_j} \phi_j^\delta$. 
Then, we also have 
\[
\begin{cases}
 \|\nabla \phi^{\mbox{\rm\tiny IK}}(t)\|_{m+4N+6} & \mbox{in the case (H1)}, \\
 \|\nabla \phi^{\mbox{\rm\tiny IK}}(t)\|_{m+4[N/2]+6} & \mbox{in the case (H2)}.
\end{cases}
\]
Moreover, by Theorem \ref{result:theorem 2} with $m$ replaced by $m+4N+7$ in the case (H1) 
and by $m+4[N/2]+7$ in the case (H2), we see that 
$(\eta^{\mbox{\rm\tiny IK}},\phi^{\mbox{\rm\tiny IK}})$ satisfy 
\[
\left\{
 \begin{array}{l}
  \partial_t\eta^{\mbox{\rm\tiny IK}}
   -\Lambda(\eta^{\mbox{\rm\tiny IK}},\delta)\phi^{\mbox{\rm\tiny IK}} = \mathfrak{r}_1, \\
  \displaystyle
  \partial_t\phi^{\mbox{\rm\tiny IK}}+\eta^{\mbox{\rm\tiny IK}}
   +\frac12|\nabla\phi^{\mbox{\rm\tiny IK}}|^2
   -\delta^2\frac{(\Lambda(\eta^{\mbox{\rm\tiny IK}},\delta)\phi^{\mbox{\rm\tiny IK}}
    +\nabla\eta^{\mbox{\rm\tiny IK}}\cdot\nabla\phi^{\mbox{\rm\tiny IK}})^2}{2(1
     +\delta^2|\nabla\eta^{\mbox{\rm\tiny IK}}|^2)}
   = \mathfrak{r}_2,
 \end{array}
\right.
\]
where $(\mathfrak{r}_1,\mathfrak{r}_2)$ satisfy 
\[
\begin{cases}
 \|(\mathfrak{r}_1(t),\mathfrak{r}_2(t))\|_{m+3} \leq C\delta^{4N+2}
  & \mbox{in the case {\rm (H1)}}, \\
 \|(\mathfrak{r}_1(t),\mathfrak{r}_2(t))\|_{m+3} \leq C\delta^{4[N/2]+2}
  & \mbox{in the case {\rm (H2)}},
\end{cases}
\]
Therefore, applying Theorem \ref{result:stability} we obtain the desired estimate 
\eqref{result:error estimate 2}. 
\quad$\Box$


\bigskip
Tatsuo Iguchi \par
{\sc Department of Mathematics} \par
{\sc Faculty of Science and Technology, Keio University} \par
{\sc 3-14-1 Hiyoshi, Kohoku-ku, Yokohama, 223-8522, Japan} \par
E-mail: iguchi@math.keio.ac.jp

\end{document}